\UseRawInputEncoding
\documentclass[12pt]
{article}
\usepackage[english]{babel}
\usepackage{amsfonts,amsmath,amssymb,amsthm,graphicx,afterpage,url}
\usepackage{epsfig}
\graphicspath{{morefig/}{figures/}{figures00/}{figures02/}{figures05/}{figlms/}{knots-figures/}}

\def\arf{\mathop{\fam0 arf}}
\def\lk{\mathop{\fam0 lk}}

\def\R{{\mathbb R}} \def\Z{{\mathbb Z}}  

\long\def\comment#1\endcomment{}

\newtheoremstyle{mydefinition}% name
  {3pt}%      Space above
  {3pt}%      Space below
  {\normalfont}%         Body font
  {\parindent}%         Indent amount (empty = no indent, \parindent = para indent)
  {\bfseries}% Thm head font
  {.}%        Punctuation after thm head
  { }%     Space after thm head: " " = normal interword space;
  {}%         Thm head spec (can be left empty, meaning `normal')

\theoremstyle{theorem}
\newtheorem{theorem}{Theorem}[section]
\newtheorem{pr}[theorem]{Problem}
\newtheorem{assertion}[theorem]{Assertion}

\newtheorem{remark}[theorem]{Remark}
\newtheorem{lemma}[theorem]{Lemma}

\begin{document}

\title{A user's guide to basic knot and link theory
\footnote{I am grateful to A. Enne who prepared the figures and a post-production version of \cite{EEF}, and to
S. Chmutov, D. Eliseev, A. Enne, M. Fedorov, A. Glebov, N. Khoroshavkina, E. Morozov,
A. Ryabichev, A. Sossinsky and R. \v Zivaljevi\'c for useful discussions and our work on \cite{EEF}.
This text is based on lectures at Independent University of Moscow (including Math in Moscow Program) and Moscow Institute of Physics and Technology, and on \cite{EEF}.}}

%What is knot theory about
%The lectures were in turn based on books \cite{PS96, CDM} which were very useful
%(in spite of some drawbacks mentioned below).

\author{A. Skopenkov
\footnote{\texttt{https://users.mccme.ru/skopenko}, Moscow Institute of Physics and Technology, Independent University of Moscow.
Supported in part by the Russian Foundation for Basic Research Grant No. 19-01-00169 and by Simons-IUM Fellowship.
}}

\date{}

\maketitle

\begin{abstract}
%{\bf Abstract.}
This is an expository paper.
We define simple invariants of knots or links (linking number, Arf-Casson invariants and Alexander-Conway polynomials) motivated by interesting results whose statements are accessible to a non-specialist or a student.  %(e.g. Theorems \ref{wha-tre} and \ref{t:hopf}).
The simplest invariants naturally appear in an attempt to unknot a knot or unlink a link.
Then we present certain `skein' recursive relations for the simplest invariants, which allow to introduce stronger invariants.
We state the Vassiliev-Kontsevich theorem in a way convenient for calculating the invariants themselves, not only the dimension of the space of the invariants.
%We also present coloring invariants.
% although we cannot explain in an elementary way how they appear.
No prerequisites are required; we give rigorous definitions of the main notions in a way not obstructing intuitive understanding.
\end{abstract}

\tableofcontents

\newpage
{\bf On the style of this text}

%A `theorem' or a `lemma' is a problem considered to be more important.

Usually I formulate a beautiful or important statement before giving a sequence of definitions and results %(lemmas, assertions, etc)
which constitute its proof.
In this case, in order to prove this statement, one may need to read
%solve
some of the subsequent material.
%problems.
I give hints on that after the statements (but I do not want to deprive you of the pleasure of finding the right moment when you finally are ready to prove the statement).
Some theorems are presented without proof, so I give references instead of hints.

%In general, if you are stuck on a certain problem, try looking at the next ones.
%They may turn out to be helpful.
%In the text they are denoted by numbers without a word `Lemma', `Theorem', etc.

In this text assertions are simple parts of a theory (for a reader already familiar with part of the material they are quick reminders).
For the same reason a small number of problems is presented.
% (if a mathematical statement is formulated as a problem, then the objective is to prove this statement).
For assertions and problems hints or solutions are presented in \S\ref{s:app}, together with proofs of theorems and lemmas.
However, a reader is recommended to prove assertions (and to solve problems) himself/herself.
In order to get a thorough understanding of the material a reader can also
consider theorems and lemmas as problems (beware the previous paragraph!).
This is peculiar not only to Zen monasteries but also to serious mathematical education, see \cite[\S1.1]{HC19}, \cite[\S1.2]{Sk20m}.

%Some part of the material is presented as a sequence of problems.
%these assertions are used later, are
%Most problems are assertions.

%More complicated
%problems
%statements are marked by stars.
Remarks are formally not used later.

%Theorems are usually hard to prove when they are stated.
%We either indicate when a reader could prove a theorem,

\section{Main definitions and results on knots}\label{0isot}

\begin{figure}[h]\centering
\includegraphics[scale=3.6]{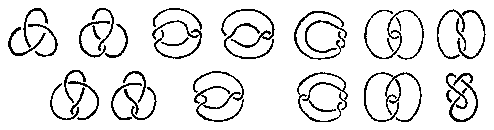}
\caption{Knots isotopic to the trefoil knot (top row) and to the figure eight knot (bottom row)}
\label{f:trefoil}
\end{figure}

We start with informal description of the main notions (rigorous definitions are given after remark \ref{r:why}).
You can imagine a {\it knot} as a thin elastic string whose ends have been glued together, see fig. \ref{f:trefoil}.
As in this figure, knots are usually represented by their `nice' plane projections called {\it knot diagrams}.
Imagine laying down  the rope on a table and carefully recording how it crosses itself
(i.e. which part lies on top of the other).
It should be kept in mind that the projections of the same knot on different planes can look quite dissimilar.

%You can visualize a knot as a thin elastic string whose end points have been glued together
%It should be kept in mind that the projections of the same knot on different planes may look quite different.

A {\bf trivial knot} is the outline (the boundary) of a triangle.

By an {\it isotopy} of a knot we mean its continuous deformation in space as a thin elastic string;
no self-intersections are allowed throughout the deformation.
%Informally, an {\it isotopy} is a continuous deformation without introducing self-intersections.
Two knots are {\it isotopic} if one can be transformed to the other by an isotopy.

\begin{assertion}\label{a:crossing}
(a) All the knots represented in the top row of fig. \ref{f:trefoil} are isotopic to each other.
(For one pair of these knots decompose your isotopy into Reidemeister moves shown in fig. \ref{f:reid}.)

%Some two knots represented in the top row of fig. \ref{f:trefoil} are isotopic to the leftmost knot in this row.

(b) The same is true for the knots represented in the bottom row of fig. \ref{f:trefoil}.

(c) All knots with the same knot diagram are isotopic.
\end{assertion}

\begin{remark}[why a rigorous definition of isotopy is necessary?]\label{r:why}
In  fig. \ref{f:naisotopy} we see an isotopy between the trefoil knot and the trivial knot.

\begin{figure}[h]\centering
\includegraphics[scale=0.6]{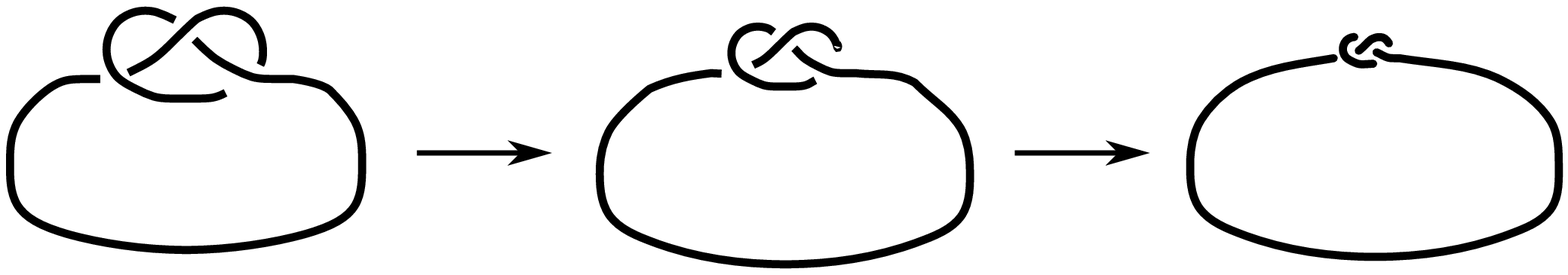}
\caption{A (non-ambient) isotopy between the trefoil knot and the trivial knot}
\label{f:naisotopy}
\end{figure}

Is it indeed an isotopy?
This is the so called `piecewise linear non-ambient isotopy', which is {\it different} from the `piecewise linear ambient isotopy' defined and used later.
(The first notion better reflects the idea of continuous deformation without self-intersections, but is hardly accessible to high school students, cf. \cite{Sk16i}.)
In fact, any two knots are piecewise linear non-ambient isotopic!

%Indeed, by the first definition any two knots are isotopic!
%According to the rigorous definition of `piecewise linear ambient isotopy' presented below, it is not.
%However, according to the rigorous definition of `piecewise linear non-ambient isotopy' it is.

A usual problem with intuitive definitions is not that it is hard to make them rigorous, but that
this can be done in different non-equivalent ways.
\end{remark}

%a reader cannot guess their meaning, but that he/she can do so ,
%and does not have time to check which way is correct.

%{\bf Rigorous definitions.}

A {\bf knot} is a spatial closed non-self-intersecting polygonal line.\footnote{This is not to be confused with {\it oriented knot} defined below in \S\ref{s:ori}.}

%Unfortunately, such a confusion is present even in good textbooks. E.g. in \cite{PS96},
%text around Theorem 1.5 presumably works with oriented knots
%CDM does not give the promised `precise definition' of a knot

A {\bf plane diagram} of a knot is its generic\footnote{A polygonal line in the plane is {\it generic} if there is a polygonal line $L$ with the same union of edges such that no 3
%three
vertices of $L$ belong to any line and no 3
%three
segments joining some vertices of $L$ have a common interior point.}
projection onto a plane\footnote{A university-mathematics terminology is `a generic image under projection onto a plane'.}, together with the information
which part of the knot `goes under' and which part `goes over' at any given crossing.

\begin{assertion}\label{p:diag} For any knot diagram there is a knot projected to this diagram.
(Such a knot need not be unique; see though assertion \ref{a:crossing}.c.)
\end{assertion}

\begin{figure}[h]\centering
\includegraphics[scale=0.5]{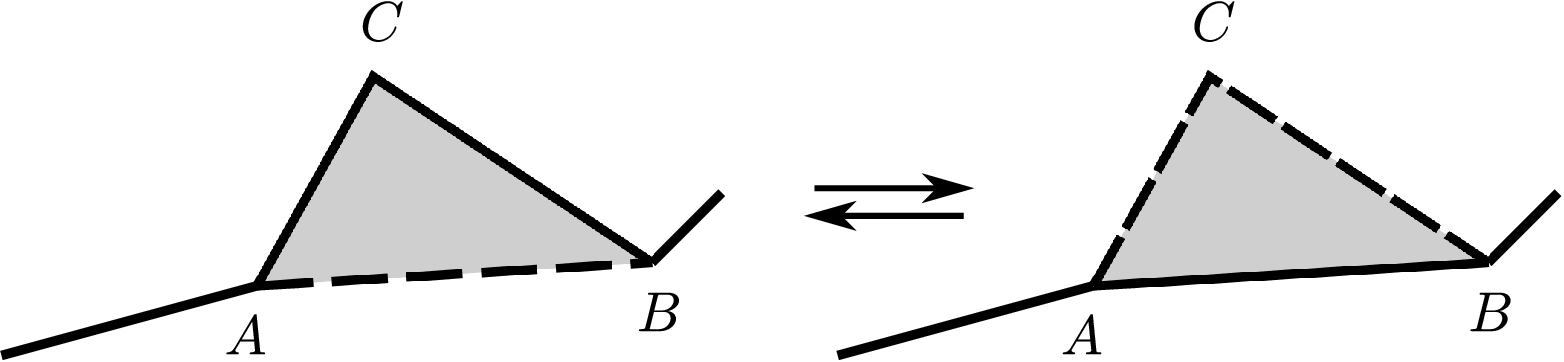}
\caption{Elementary move}
\label{f:elem}
\end{figure}

%(\cite[Fig. 1.4]{PS96})

Suppose that two sides $AC$ and $CB$ of a triangle $ABC$ are edges of a knot.
Moreover, assume that the knot and (the part of the plane bounded by) the triangle $ABC$  do not intersect at any other points.
An {\bf elementary move} $ACB\to AB$ is the replacement of the two edges $AC$ and $CB$ by the edge $AB$, or the inverse operation $AB\to ACB$ (fig. \ref{f:elem}).\footnote{If the triangle $ABC$ is degenerate, then elementary move is either subdivision of an edge or inverse operation.}
Two knots $K$ and $L$ are called (piecewise linearly ambiently) {\bf isotopic} if there is a sequence of knots $K_1,\ldots,K_n$ such that $K_1=K$, $K_n=L$ and every subsequent knot $K_{j+1}$ is obtained from the previous one $K_j$ by an elementary move.
% (of the type described above).
%\footnote{It would be nice to have a reference to a published proof that this definition
%equivalent to the standard one \cite[Definition 1.1]{Skopenkov2016i}.}

%Two plane diagrams of knots are called {\it Reidemeister equivalent} if one can be obtained from the other
%by Reidemeister moves \cite[Fig. 4.1]{Pr95}.

\begin{theorem}\label{wha-tre} (a) The following knots are pairwise not isotopic: the trivial knot, the trefoil knot, the figure eight knot.

%(a) The trivial knot, the trefoil knot and the figure eight knot are pairwise not isotopic.

(b) There is an infinite number of pairwise non-isotopic knots.
\end{theorem}

This is proved using {\it Arf} and {\it Casson invariants}, see \S\ref{0arf} and \S\ref{0cas}, cf. \S\ref{s:reco}.

%or using {\it proper colorings,} see

%(same is true for theorem \ref{arf-k7} below).
%In \S\ref{s:reco} we sketch an alternative proof of a part of theorem \ref{wha-tre}.

%{\it Amphisphericity, or achirality.}

The {\it mirror image} of a knot $K$ is the knot whose diagram is obtained by changing all the crossings (fig. \ref{f:crossing}) in a diagram of $K$.
By assertion \ref{a:crossing}.b the figure eight knot is isotopic to its mirror image.

%(The mirror image of an isotopy class is well defined since Reidemeister moves are preserved under mirror
%images, see below.)

\begin{theorem}\label{wha-nchir} The trefoil knot is not isotopic to its mirror image.
\end{theorem}

Theorem \ref{wha-nchir}
%and \ref{wha-ninv} below are
is proved using {\it the Jones polynomial} \cite[\S3]{PS96}, \cite[\S2.4]{CDM}.
The proof is outside the scope of this text.

\begin{theorem}[Conway--Gordon; cf. {\cite[Theorem~2]{CG83}}]\label{arf-k7} Take any 7 points in 3-space, no four of which belong to any plane.
Take ${7\choose2}=21$
%non-self-intersecting polygonal lines (e.g.
segments joining them.
%Suppose that every two polygonal lines intersect only at their common end, if there is such an end
%(in particular, the polygonal lines are disjoint if they do not have a common end).
Then there is a closed polygonal line formed by taken
%polygonal lines
segments and non-isotopic to the boundary of a triangle.
%More generally, if the graph $K_7$ is} piecewise linearly embedded {\it in space, then there exists a knotted cycle in this graph.}
\end{theorem}

This is proved using {\it Arf invariant}, see \S\ref{0arf}.
The details are outside the scope of this text.

%\newpage
\section{Main definitions and results on links}\label{0isoli}

A {\bf link} is a collection of pairwise disjoint knots, which are called the {\it components} of the link.
Ordered collections are called ordered or colored links, while non-ordered collections are called non-ordered or non-colored links.
In this text we abbreviate `ordered link' to just `link'.

A {\bf trivial link} (with any number of components) is a link formed by triangles in parallel planes.

{\it Plane diagrams} and {\it isotopy} for links are defined analogously to knots.

\begin{figure}[h]\centering
\includegraphics{3-1.eps}\includegraphics[scale=0.6]{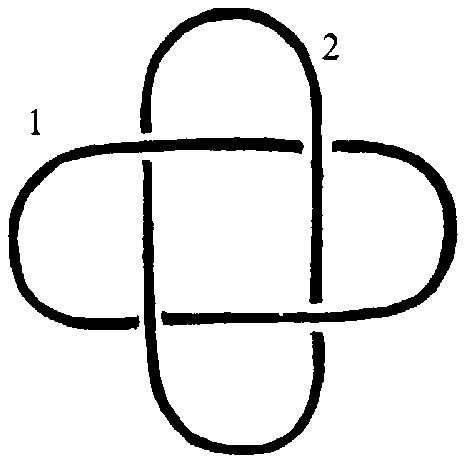}
\includegraphics[scale=0.5]{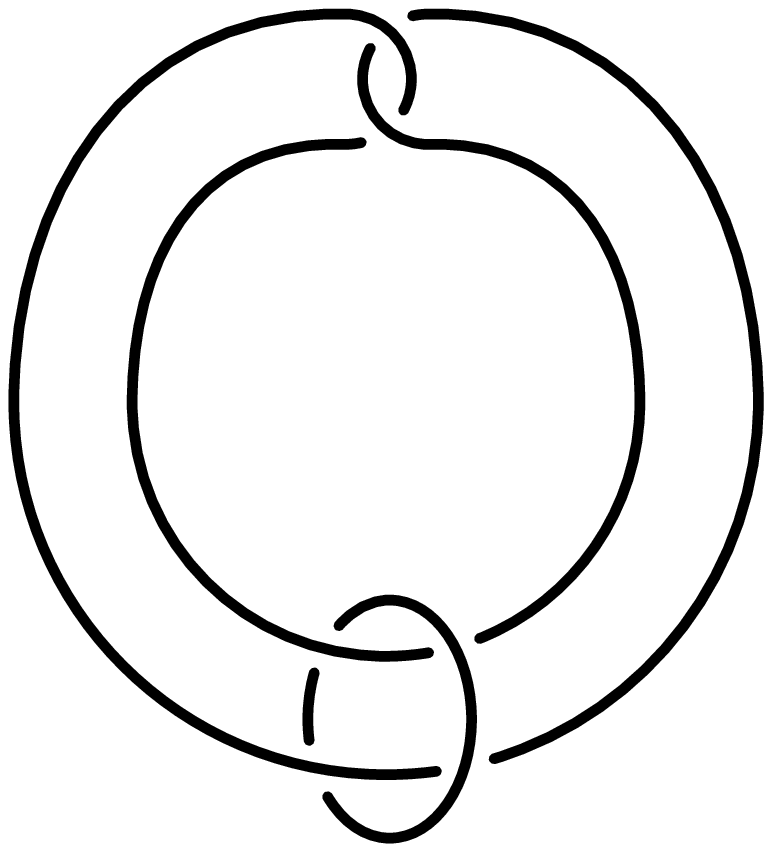}\includegraphics{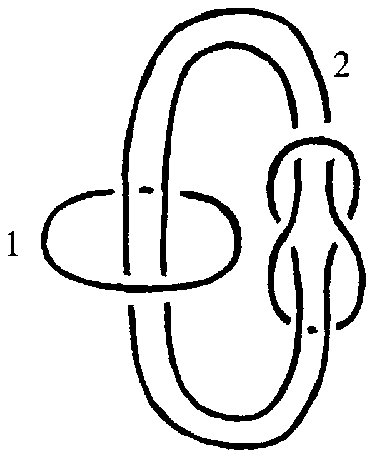}
\caption{The Hopf link, the trivial link and another three links}\label{f:hopf}
\end{figure}

%\cite{Pr95} fig. 2.14a \cite[Fig. 4.5]{Pr95}

\begin{figure}[h]\centering
\includegraphics{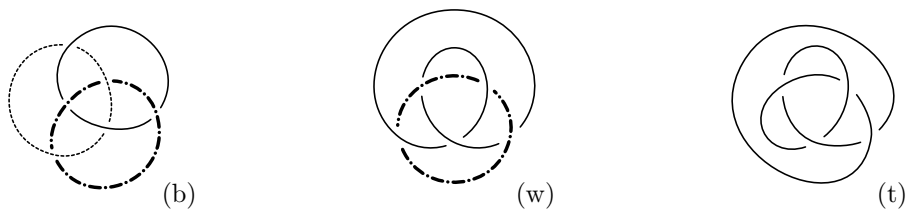}
%na noutbuke ne kompiliruetsja? net, dvipdf
\caption{The Borromean rings, the Whitehead link and the trefoil knot}\label{f:borwhitre}
\end{figure}

%\newpage

%(a) Define (piecewise linearly) isotopic {\it pairs} of disjoint knots.

\begin{assertion}\label{a:excha} (a) The Hopf link is isotopic to the link obtained from the Hopf link by switching the components.

%The following three links are isotopic: the Hopf link, the mirror image of the Hopf link and

(b) The Hopf link is isotopic to some link whose components are symmetric with respect to some straight line.

%Problems  2.4.a*b* from \cite[\S2]{Pr95}.

(c) The fourth link in fig.  \ref{f:hopf} is isotopic to the Whitehead link in fig. \ref{f:borwhitre}.w.

(d,e) The same as in (a,b) for the Whitehead link.

(f) The Borromean rings link is isotopic to a link whose components are permuted in a cyclic way under the rotation by angle $2\pi/3$ with respect to some straight line.
\end{assertion}

\begin{proof} (a) This follows by (b) (or can be proved independently).

%(b) It is not difficult to find one of the two possible lines with such symmetry.

(d) This follows by (e) (or can be proved independently).

%\cite[\S2]{Pr95}

\begin{figure}[h]\centering
\includegraphics[scale=0.9]{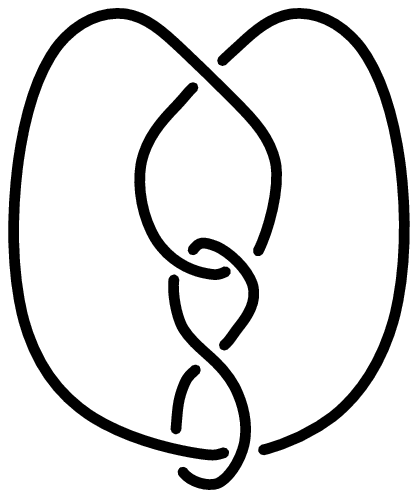}\qquad
\includegraphics[scale=0.9]{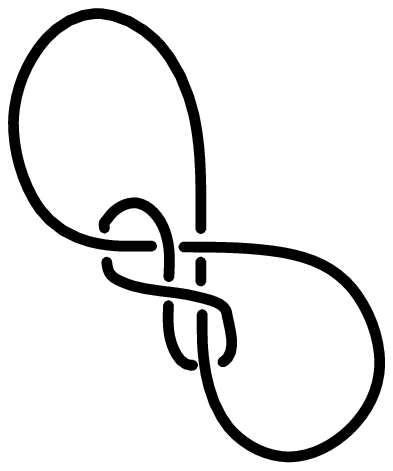}\qquad
\includegraphics[scale=0.9]{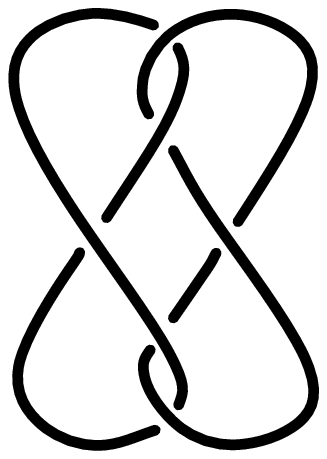}\qquad
\caption{Isotopy of the Whitehead link}
	\label{f:whi}
\end{figure}

%\cite{Pr95} fig. 2.14bcd

(e) See figure \ref{f:whi}.

(f) Take the quadrilaterals from figure \ref{f:borr}, left.
%circumscribed around these ellipses and symmetric w.r.t. the coordinate axes.
Then the straight line is the bisector of any octant formed by the quadrilaterals.

\begin{figure}[h]\centering
\includegraphics[scale=0.5]{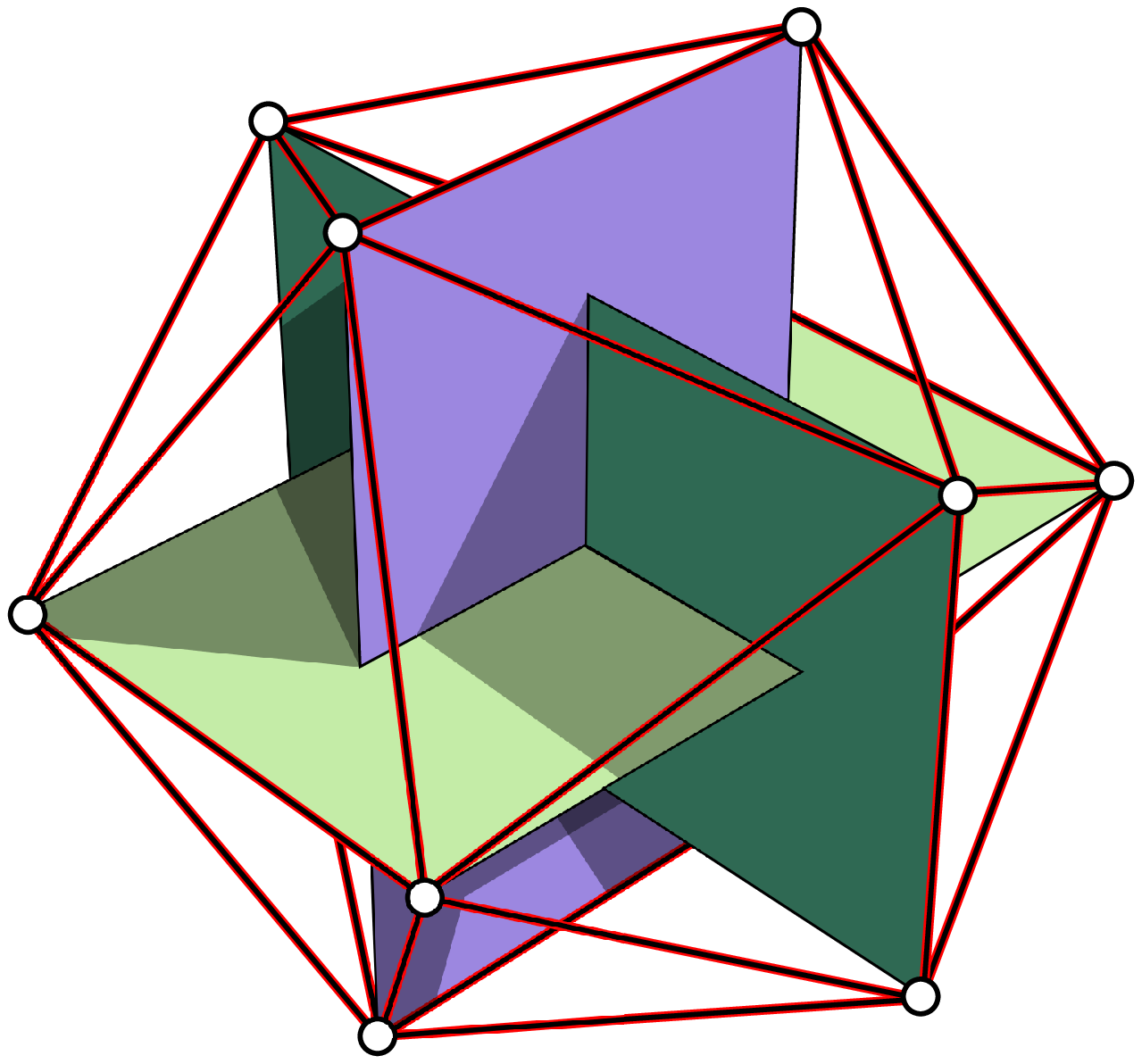}\qquad\includegraphics[scale=0.9]{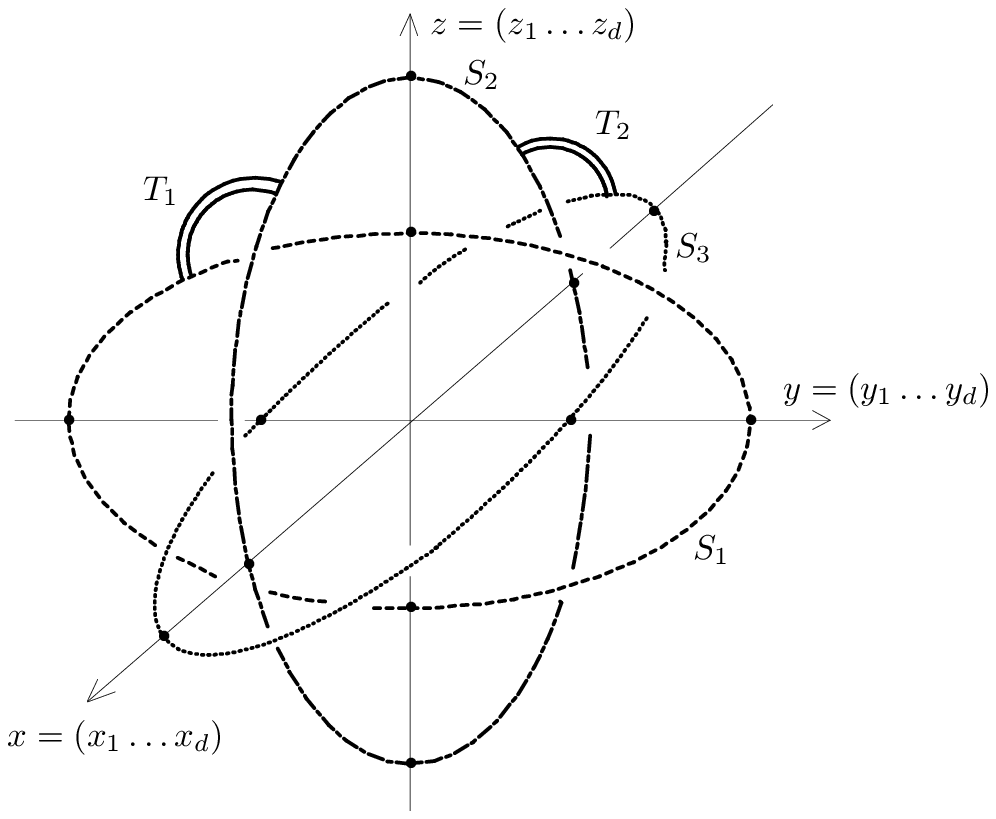}
\caption{Borromean rings}
\label{f:borr}
\end{figure}

There is also the following beautiful curvilinear construction.
Take three ellipses given by the following three systems of equations:
$$\left\{\begin{array}{c} x=0\\ y^2+2z^2=1\end{array}\right., \qquad
\left\{\begin{array}{c} y=0\\ z^2+2x^2=1\end{array}\right. \qquad\text{and}\qquad
\left\{\begin{array}{c} z=0\\ x^2+2y^2=1. \end{array}\right.$$
See figure \ref{f:borr}, right.
Then the straight line is given by $x=y=z$.
\end{proof}

\begin{theorem}\label{t:hopf}
(a) The following links are pairwise non-isotopic: the Hopf link, the trivial link, the Whitehead link.

(b) The Borromean rings link is not isotopic to the trivial link.
%{\it Neither the Hopf link nor the Whitehead link nor the Borromean rings link is isotopic to the trivial link.}
\end{theorem}

%For the Hopf link this is
This is proved using {\it linking number modulo 2}, invent it yourself or see \S\ref{0rapl2},
%The rest of the theorem is
and {\it the Alexander-Conway polynomials}, see \S\ref{0con}.
Alternatively, one can use the `triple linking' (Massey-Milnor) number and `higher linking' (Sato-Levine) number  \cite[\S4.4-\S4.6]{Sk}.
%Part (a) can also be proved using {\it proper colorings}, see \S\ref{s:reco}.

\begin{theorem}[Conway--Gordon, Sachs; {\cite[Theorem 1.1]{Sk14}, cf. \cite[Theorem~1]{CG83}}]\label{rlt-cgslin}  If no 4 of 6 points in 3-space  lie in the same plane, then there are two linked triangles with vertices at these 6 points.
That is, the part of the plane bounded by the first triangle intersects the outline of the second triangle exactly at one point.
%Moreover, the number of linked unordered pairs of triangles with vertices at these 6 points is odd.
\end{theorem}

This is proved using {\it linking number modulo 2}, see \S\ref{0rapl2}.
The details are outside the scope of this text, see \cite{Sk14}.

\section{Some basic tools}\label{s:tools}

\begin{remark}[some accurate arguments]\label{r:accur} In the following paragraph we prove that {\it if a knot lies in a plane, then the knot is isotopic to the trivial knot.}

%Denote by  the elementary moves of fig.~\ref{f:elem}.
Denote the knot in a plane by $M_1M_2\ldots M_n$.
Take a point $Z$ outside the plane.
Then $M_1M_2\ldots M_n$ is transformed to the trivial knot $M_1ZM_n$ by the following sequence of elementary moves:
$$M_1M_2\to M_1ZM_2,\quad ZM_2M_3\to ZM_3,\quad ZM_3M_4\to ZM_4,\quad\ldots,\quad ZM_{n-1}M_n\to ZM_n.$$

The following result shows that intermediate knots of an isotopy from a knot lying in a plane to the trivial knot
 can be chosen also to lie in this plane.

{\it Schoenflies theorem.} Any closed polygonal line without self-intersections in the plane is isotopic (in the plane) to a triangle.

%Any two closed polygonal lines without self-intersections in the plane are isotopic in the plane.

This is a stronger version of the following celebrated result.

{\it Jordan theorem.} Every closed non-self-intersecting polygonal line $L$ in the plane $\R^2$ splits the plane into exactly two parts, i.e. $\R^2-L$ is not connected and is a union of two connected sets.

A subset of the plane is called {\it connected} if every two points of this subset can be connected by a polygonal line lying in this subset.
%Deduce (j) from (s).
%In fact, (j) has a direct proof used in the proof of (s).

For an algorithmic explanation why the Jordan Theorem (and so the Schoenflies Theorem) is non-trivial,
and for a proof of the Jordan Theorem, see \S1.3 `Intersection number for polygonal lines in the plane' of \cite{Sk18}, \cite{Sk}.
\end{remark}

%We do not require as rigorous proof of lemma \ref{p:arfmot} as a formal reduction to the Schoenflies Theorem. %Same holds for other results (like assertions \ref{l:crossing} and \ref{l:reid}).
%{\it Remark.} This idea can be made rigorous using the following result.

\begin{assertion}\label{a:arfmot} Suppose that there is a point on a knot such that if we go around the knot starting from this point, then on some plane diagram we first meet only overcrossings, and then only undercrossings.
Then the knot is isotopic to the trivial knot.\footnote{This assertion would be a motivation for introduction of the Arf invariant (\S\ref{0arf}). The proof illustrates in low dimensions one of the main ideas of the celebrated Zeeman's proof of the higher-dimensional Unknotting Spheres Theorem, see the survey \cite[Theorem 2.3]{Sk16c}.}
\end{assertion}

\begin{figure}[h]\centering
\includegraphics[scale=0.5]{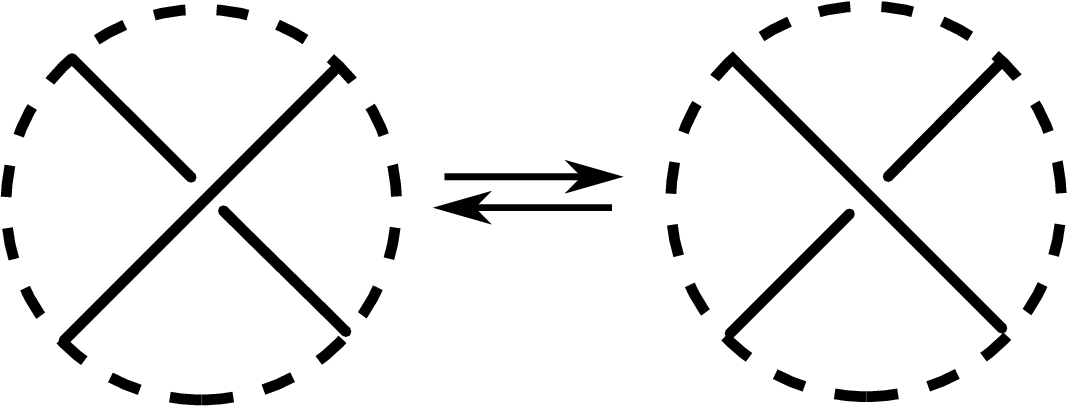}
\caption{Crossing change.
The plane diagrams are identical outside the disks bounded by dashed circles.
%as in figure \ref{f:crossing}, while outside this disk the plane diagram  remains unchanged.
No other sides of the plane diagrams except for the pictured ones intersect the disks.
(Same for fig. \ref{f:reid}, \ref{f:planiso}, \ref{f:skeino}, \ref{f:pass} and \ref{f:skein}.)}
\label{f:crossing}
\end{figure}

%\cite[p. 42]{CDM}

A {\bf crossing change} is change of overcrossing to undercrossing or vise versa, see fig. \ref{f:crossing}.

%altering one of the crossings, reversing which arc goes over and which goes under at one of the crossing points,

Clearly, after any crossing change on the leftmost diagrams of the trefoil knot and the figure eight shown in fig. \ref{f:trefoil}
we obtain a diagram of a knot isotopic to the trivial knot.

\begin{lemma}\label{l:crossing} Every plane diagram of a knot
%(or a link)
can be transformed by  crossing changes to a plane diagram of a knot isotopic to the trivial knot.\footnote{This simple lemma will be used for construction of invariants using recursive ({\it skein}) relations, see \ref{con-lk2},   \ref{con-thma},\ref{con-lk2},\ref{con-thm},\ref{con-ale} и \ref{vas-kon}.}
\end{lemma}

%(In the proof you can use without proof the Schoenfliess theorem \ref{t:schoe}.a below.)

\begin{figure}[h]\centering
\includegraphics[scale=0.45]{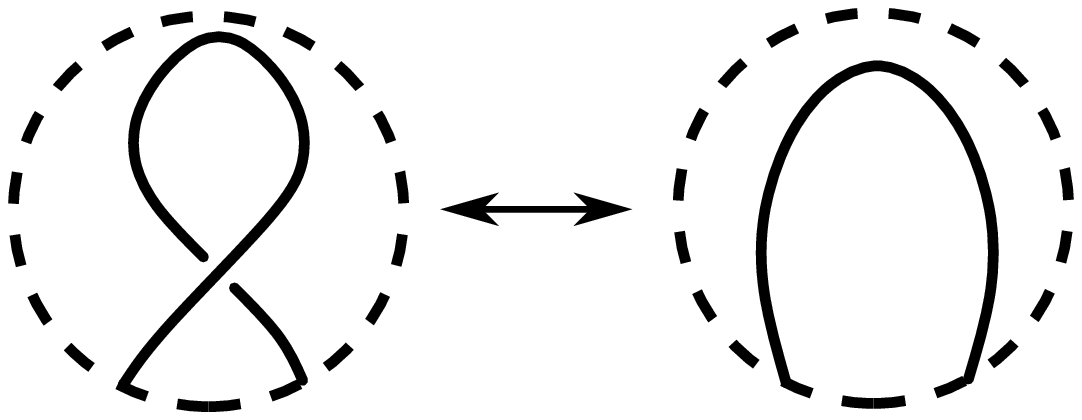}\qquad
\includegraphics[scale=0.45]{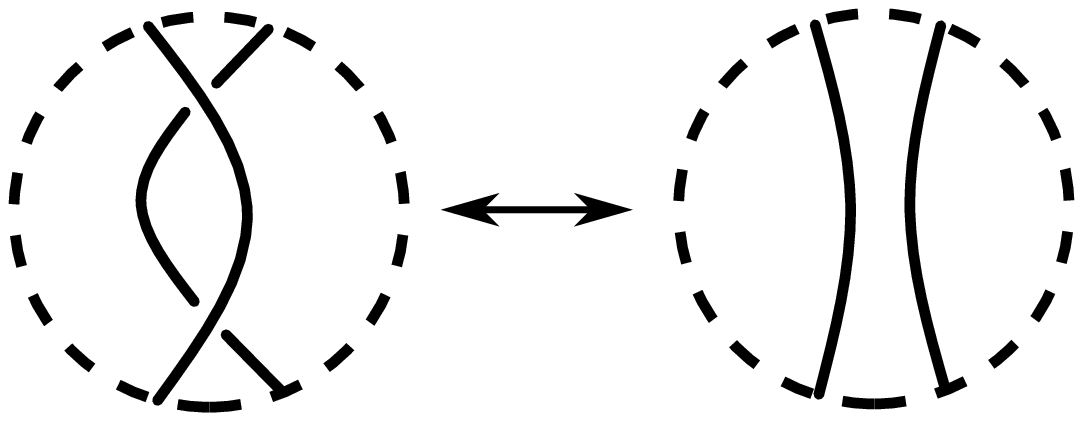}\qquad
\includegraphics[scale=0.45]{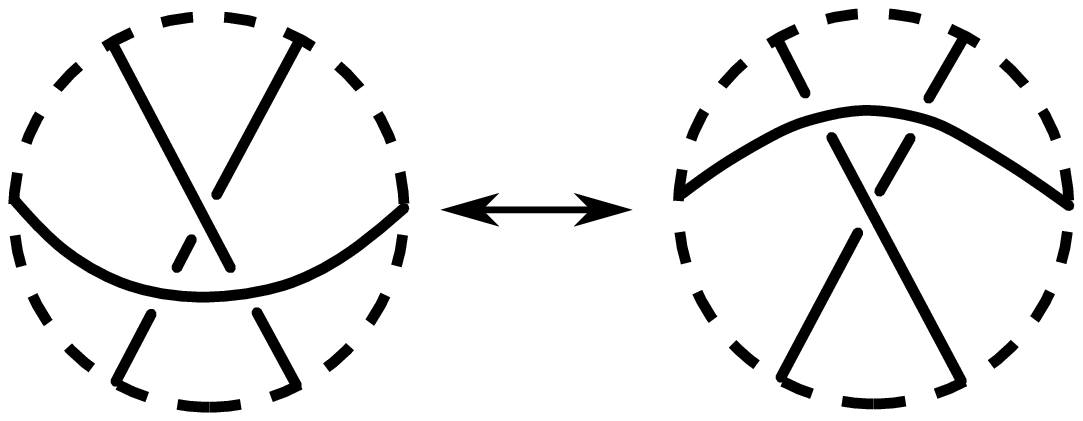}
\caption{Reidemeister moves}
\label{f:reid}
\end{figure}

\begin{figure}[h]\centering
\includegraphics[scale=0.45]{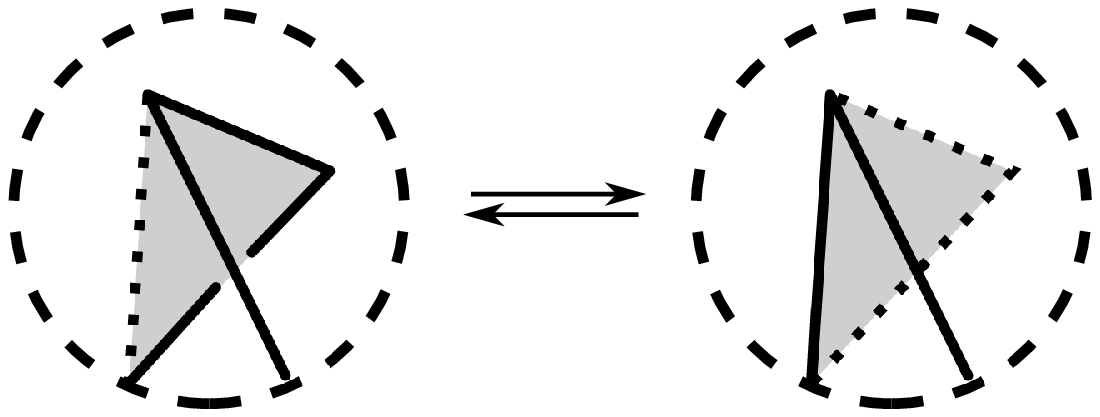}\qquad\qquad
\includegraphics[scale=0.45]{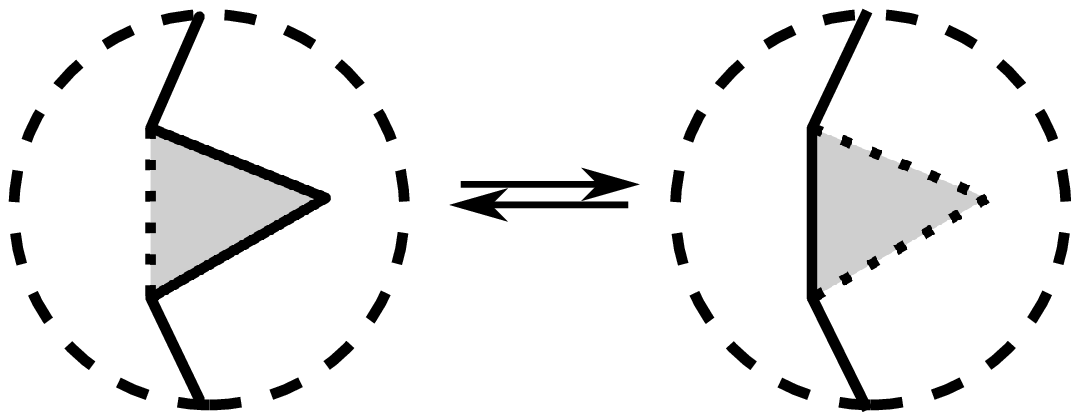}\quad
\includegraphics[scale=0.45]{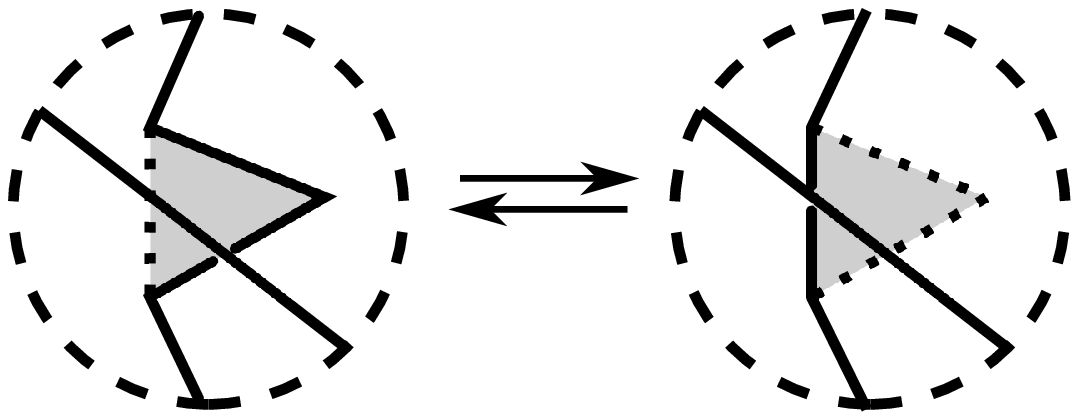}
\caption{Left: to a rigorous definition of the first Reidemeister move
\newline
Middle, right: plane isotopy moves}
\label{f:planiso}
\end{figure}

%\cite[Fig. 1.9]{PS96} \cite[Fig. 1.10]{PS96}

In this text instead of knots up to isotopy we shall study plane diagrams of knots up to (equivalence generated by) {\bf Reidemeister moves} (shown in fig. \ref{f:reid})\footnote{A rigorous definition of the first Reidemeister move is easily given using fig. \ref{f:planiso} (left). The other Reidemeister moves have analogous rigorous definitions.
%The participants are not required to use these rigorous definitions in solutions.
% (except for problems \ref{p:arfmot} and lemma \ref{l:crossing}
%which allow a simple solutions not using Reidemeister moves).
You can use informal description of Reidemeister moves in fig. \ref{f:reid} and so ignore plane isotopy moves.
See footnote \ref{fo:reid}.}
and {\it plane isotopy moves} shown in fig. \ref{f:planiso} (middle, right).
I.e. we shall use without proof the following result.
%Reidemeister theorem \ref{l:reid} below.

\begin{theorem}[Reidemeister]\label{l:reid} Two knots are isotopic if and only if some plane diagram of the first knot can be obtained from some plane diagram of the second one by Reidemeister moves and plane isotopy moves.
\end{theorem}

See \cite[\S1.7]{PS96}.\footnote{\label{fo:reid} Since \cite[\S1.6]{PS96} does not contain as rigorous definition of Reidemeister moves as that of plane isotopies, the argument in \cite[\S1.7]{PS96} does not constitute a rigorous proof.
We believe that a rigorous proof can be recovered using rigorous definition of Reidemeister moves.
\newline
This shows that having plane isotopy in the statement \cite[\S1.7]{PS96} does not make the statement rigorous, and thus should be avoided.
On an intuitive level, plane isotopies should better be ignored.
With the alternative rigorous definition below, plane isotopies can be expressed via Reidemeister moves and so again should better be ignored in the statement.
\newline
Let us present an alternative rigorous definition of the first Reidemeister move.
The other Reidemeister moves have analogous rigorous definitions.
% we present this definition in order to show that rigorous definitions are easy to give, and to justify
%(for a specialist) that our simplified formulation of lemma \ref{l:reid} is equivalent to the standard one.
On the plane take a closed non-self-intersecting polygonal line $L$ whose interior (see the Jordan theorem in remark \ref{r:accur}) intersects a knot diagram $D$ by a non-self-intersecting polygonal line $M$ joining two points on $L$.
Let $N$ be a closed non-self-intersecting polygonal line in the interior of $L$ such that $N\cap L=\emptyset$, $N\cap M$ is one point and $M\cup N$ is a generic (self-intersecting) polygonal line.
{\it The first Reidemeister move} is replacement of $M$ to $M\cup N$ in $D$, with any `information'
at the appearing crossing.}
The analogues of lemma \ref{l:crossing} and theorem \ref{l:reid} for links are correct.

%Alternatively, one can rigorously define Reidemeister moves
%analogously to definition of plane isotopies   \cite[beginning of \S1.6]{PS96}.

%It would be nice to have a reference to a published proof that a simple definition of plane isotopies %\cite[\S1.6]{PS96} is equivalent to the standard one \cite[Definition 1.1]{Sk16i}
%(which is hardly accessible to high-school students).

%The standard formulation \cite[\S1.7]{PS96} involves  (not to be confused with isotopy of knots on the plane),

\section{The Gauss linking number modulo 2 via plane diagrams}\label{0rapl2}

Suppose that there is an isotopy between two 2-component links, and the second component is fixed throughout the isotopy.
Then the trace of the first component is a self-intersecting cylinder disjoint from the second component.
If after the isotopy the components are unlinked, then the cylinder can be completed to
a self-intersecting disk disjoint from the second component.
This observation, together with \cite[the Projection lemma 4.3.2]{Sk}, motivates the following definition.

The {\bf linking number modulo two} $\lk_2$ of the plane diagram of a 2-component link is the number modulo 2 of crossing points on the diagram at which the first component passes above the second component.

%\ref{rap-link} below.

%\begin{figure}[h]\centering
%\includegraphics{3-16.eps}\qquad\includegraphics{3-18.eps}
%\caption{Some links}
%\label{f:links}
%\end{figure}

\begin{pr}\label{a:lk2} Find the linking number modulo 2 for the plane diagrams in fig. \ref{f:hopf} and
for pairs of Borromean rings in fig. \ref{f:borwhitre}.b.
\end{pr}

\begin{lemma}\label{l:lk2} The linking number modulo 2 is preserved under Reidemeister moves.
\end{lemma}

This lemma is easily proved separately for every Reidemeister move.

By lemma \ref{l:lk2} the {\bf linking number modulo 2} of a 2-component link (or even of its isotopy class) is well-defined
by setting it to be the linking number modulo 2 of any plane diagram of the link.

We shall use without proof the following {\it Parity lemma:} any two closed polygonal lines in the plane whose vertices are in general position intersect at an even number of points.
For a discussion and a proof see \S1.3 `Intersection number for polygonal lines in the plane' of \cite{Sk18, Sk}.

\begin{assertion}\label{a:lk2-pr}
(a) Switching the components of a 2-component link preserves the linking number modulo 2.

(b) There is a 2-component link which is not isotopic to the trivial link but whose linking number modulo 2 is  zero.
\end{assertion}

Part (b) is proved using {\it integer-valued} linking coefficient, see \S\ref{0rapl}.

%(b) Taking the mirror image preserves the linking number modulo 2.

\begin{theorem}\label{con-lk2} There is a unique mod2-valued isotopy invariant $\lk_2$ of 2-component links
that assumes value 0 on the trivial link and such that for any links $K_+$ and $K_-$
whose plane diagrams differ by a crossing change of a crossing $A$
$$\lk\phantom{}_2K_+-\lk\phantom{}_2K_- =
\begin{cases}1 &A\text{ is the crossing of different components;} \\
0 &A\text{ is the self-crossing of one component.}\end{cases}$$
\end{theorem}

\begin{assertion}\label{a:lk2comp} If the linking number modulo 2 of two (disjoint outlines of) triangles in space is zero, then the link formed by the triangles is isotopic to the trivial link.
\end{assertion}

The proof is presumably unpublished but not hard.
We encourage a reader to publish the details.
Cf. \cite{Ko19}.

%(a) Any link $L$ of two triangles is isotopic to the link obtained from $L$ by switching the components.
%(a) We do not know any proof except deduction from (b). No, rotation!

\section{The Arf invariant}\label{0arf}

%\cite[1.8.4]{CDM}\cite[3.6.7]{CDM}
%In the lecture it is explained how the following notion of  Arf number arises in attempts to unknot a knot.

Take a plane diagram of a knot and a point $P$ on the diagram different from crossing points.
Call $P$ a {\it basepoint}.
A non-ordered pair of crossing points $A$ and $B$ is called {\bf skew} (or $P$-skew) if going around the diagram in some direction starting from $P$ and marking only crossings at $A$ and $B$, we first mark overcrossing at $A$, then undercrossing at $B$, then undercrossing at $A$, and at last overcrossing at $B$.
The {\it $P$-Arf invariant} of the plane diagram is the parity of the number of all skew pairs of crossing points.

%, see problem \ref{con-cas}.e.

%Hereafter the point $P$ from the definition of the Arf number is called .
%The corresponding skew pairs of crossings are called {\it $P$-skew}.
%The Arf number calculated with $P$ as a basepoint is called {\it $P$-Arf number} $\arf_P$.

%An ordered pair of crossing points, $A$ above $A'$ and $B$ above $B'$, is called {\it linked} (or {\it skew})
%if going around the knot (along the orientation) starting from $P$ we first pass through $A$, then through $B'$, %then through $A'$, and at last through $B$.
%intersecting arrows for which the base point is situated between the ends of the arrows.

\begin{pr}\label{e:con-cas} (a) If the $P$-Arf invariant of a plane diagram is non-zero, then $P$ is not a point as in assertion \ref{a:arfmot}.

(b) Find the $P$-Arf invariant (of some plane diagram) of the trivial, the trefoil and the figure eight knots (for arbitrary choice of a basepoint $P$).
\end{pr}

\begin{lemma}\label{p:con-cas-ab} (a) The $P$-Arf invariant is independent of the choice of a basepoint $P$.

(b) The Arf invariant of a plane diagram is preserved under Reidemeister moves.
\end{lemma}

By (a) the {\it Arf invariant} of a plane diagram is well-defined by setting it to be the $P$-Arf invariant for any basepoint $P$.
So the statement of (b) makes sense.
By (b) the {\bf Arf invariant} (Arf number) $\arf$ of a knot (or even of isotopy class of a knot) is well-defined by setting it to be the Arf invariant of any plane diagram of the knot.

\emph{Hints.} (a) It suffices to show that the Arf invariant remains unchanged when the basepoint moves through one crossing on the plane diagram.

(b) Prove the statement for each Reidemeister move separately.
Cleverly choose a basepoint!

\begin{assertion}\label{a:arf-pr}
There is a knot which is not isotopic to the trivial knot but which has zero Arf invariant.
\end{assertion}

This is proved using \emph{Casson} invariant, see \S\ref{0cas}.

\begin{figure}[h]\centering
\includegraphics[scale=0.6]{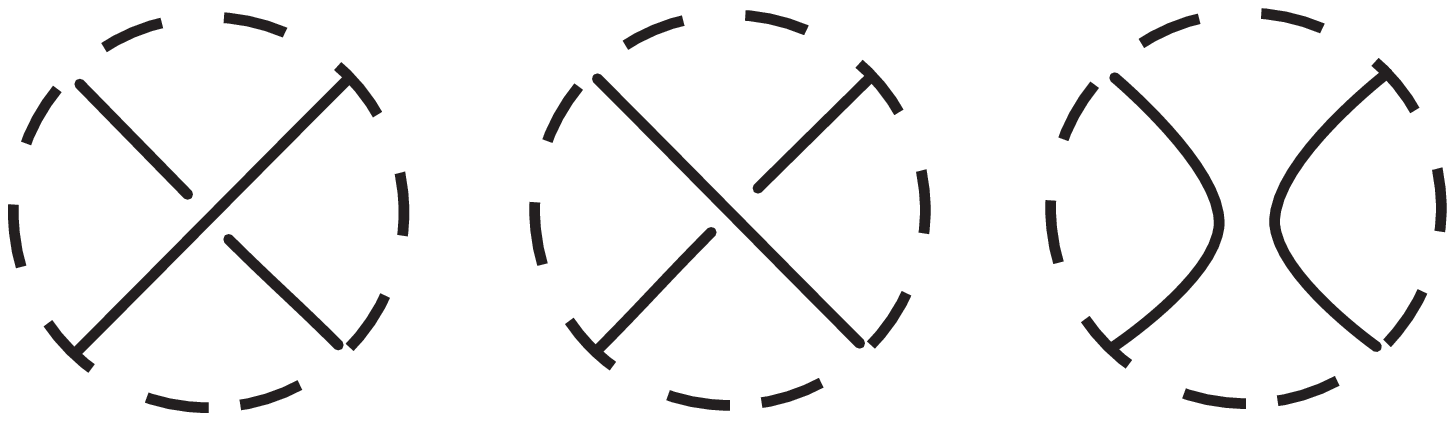}
\caption{Links $K_+,K_-,K_0$}
\label{f:skeino}
\end{figure}

\begin{theorem}\label{con-thma} There is a unique mod2-valued isotopy invariant $\arf$ of knots that assumes value 0 on the trivial knot and such that
$$\arf K_+-\arf K_-=\lk\phantom{}_2K_0.$$
for any knots $K_+$ and $K_-$ whose plane diagrams differ as shown in fig. \ref{f:skeino} so that $K_0$ is a 2-component link. (The latter is equivalent to the existence of the orientation for which fig. \ref{f:skeino} becomes fig. \ref{f:skein}.)
\end{theorem}

\begin{figure}[h]\centering
\includegraphics[scale=0.8]{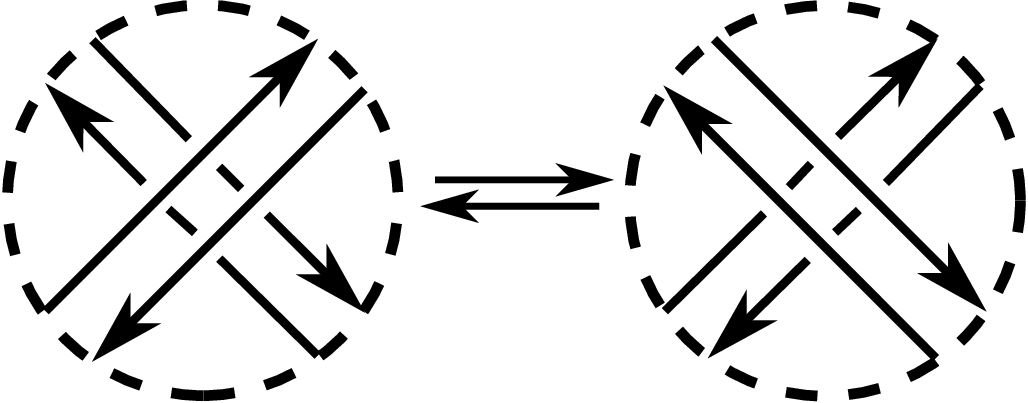}
\caption{Pass move}
\label{f:pass}
\end{figure}

%(\texttt{https://arxiv.org/pdf/1406.5573.pdf}, Definition 3.1)

\begin{assertion}\label{lk-pamo} Two knots  are called \emph{pass equivalent} if some plane diagram of the first knot (with some orientation) can be transformed
to some plane diagram of the second knot (with some orientation) using Reidemeister moves and \emph{pass moves} of fig. \ref{f:pass}.

(a) If two knots are pass equivalent, then their Arf invariants are equal.

(b) The eight figure knot is pass equivalent to the trefoil knot.

(c) \cite[pp. 75--78]{Ka87} If the Arf invariants of two knots are equal, then the knots are pass equivalent.
\end{assertion}

%\newpage
\section{Appendix: proper colorings}\label{s:reco}

%This section is not used in the sequel and so could be omitted.

%This section illustrates relations with colourings in combinatorics \cite{Ra20}.

Section \ref{s:reco}  only uses the material of \S\ref{0isot} and \S\ref{0isoli}.
See more in \cite{Pr98}.

A {\it strand} in a plane diagram (of a knot or link) is a connected piece that goes from one undercrossing to the next.
%The number of strands is the same as the number of crossings.
A {\bf proper coloring} of a plane diagram (of a knot or link) is a coloring of its strands
%from the definition of 3-colorability.
in one of three colors so that at least two colors are used, and at each crossing, either all three colors are present or only one color is present.
A plane diagram (of a knot or link) is {\bf 3-colorable} if it has a proper coloring.

\begin{pr}\label{wha-cole} For each of the following knots or links take any diagram and decide if it is 3-colorable.

(a) the trivial knot.
\quad
(b) the trefoil knot.
\quad
(c) the figure eight knot.

(d-i) links in fig. \ref{f:hopf} and \ref{f:borwhitre}.b.
\end{pr}

\begin{lemma}[{\cite[pp. 29-30, Theorem 4.1]{Pr95}}]\label{wha-colre}
The 3-colorability of a plane diagram is preserved under the Reidemeister moves.
\end{lemma}

%\begin{pr}\label{wha-col} (a-j) For every of the knots and links of problem \ref{wha-cole}
%take any diagram and find the number of its proper colorings.
%\end{pr}

\begin{theorem}\label{wha-colr} (a) Neither of links in fig. \ref{f:hopf} and \ref{f:borwhitre} (except the trivial link) is isotopic to the trivial link.

(b) The $5_1$ knot is not isotopic to the trivial knot.
\end{theorem}

\begin{figure}[h]\centering
\includegraphics[scale=0.9]{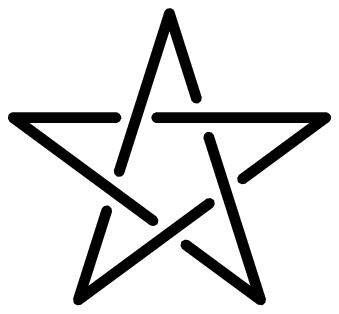}
\caption{The $5_1$ knot}
\label{f:51}
\end{figure}

%\newpage
\section{Oriented knots and links and their connected sums}\label{s:ori}

\begin{figure}[h]\centering
\includegraphics[scale=0.5]{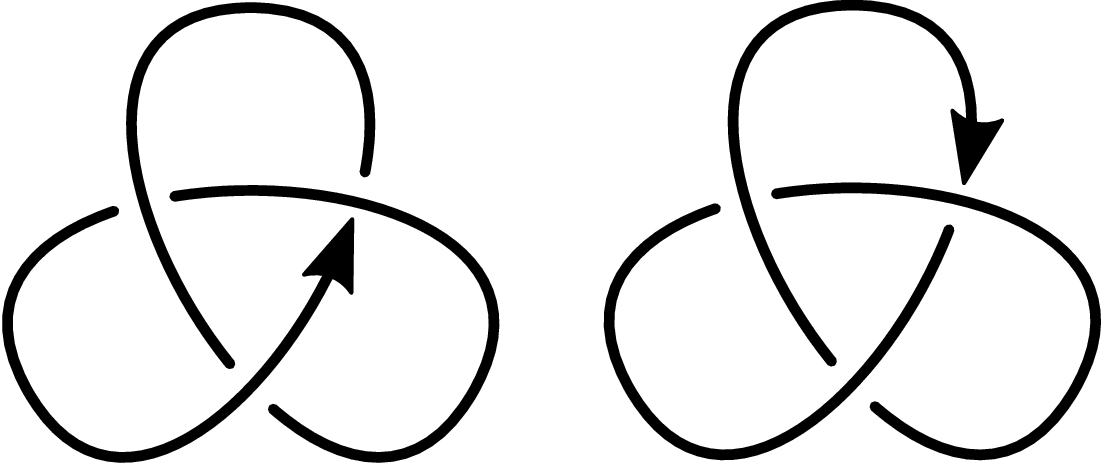}
\caption{Two trefoil knots with the opposite orientations}
\label{f:ori}
\end{figure}

One knows what is oriented polygonal line, so one knows what is oriented knot (fig. \ref{f:ori}).

Both the informal notion and rigorous definition of {\it isotopic} oriented knots
are given analogously to isotopic knots.

\begin{assertion}\label{wha-chirsph} Isotopic oriented polygonal lines without self-intersections
%tangles
on the plane and on the sphere
%and on the torus
are defined analogously to isotopic oriented knots in space.
%(The mathematical name for tangle is `knot in the plane', and `tangle' means something else.)

(a) An oriented spherical triangle is isotopic on the sphere  to the same triangle with the opposite orientation.

(b) The analogue of (a) for the plane is false.
\end{assertion}

%(c)* Any oriented circle on the torus is not isotopic {\it on the torus} to the same circle with
%the opposite orientation.
%(This is proved using homology.)

\begin{assertion}\label{wha-inv} The following pairs of knots with opposite orientations are isotopic: two trivial knots, two trefoil knots, two figure eight knots.
\end{assertion}

\begin{theorem}[Trotter, 1964]\label{wha-ninv}
There exists an oriented knot which is not isotopic to the same knot with the opposite orientation.
\end{theorem}

This is proved using {\it the Jones polynomial} \cite[\S3]{PS96}, \cite[\S2.4]{CDM}; the proof is outside the scope of this text.

\begin{figure}[h]\centering
\includegraphics[scale=0.5]{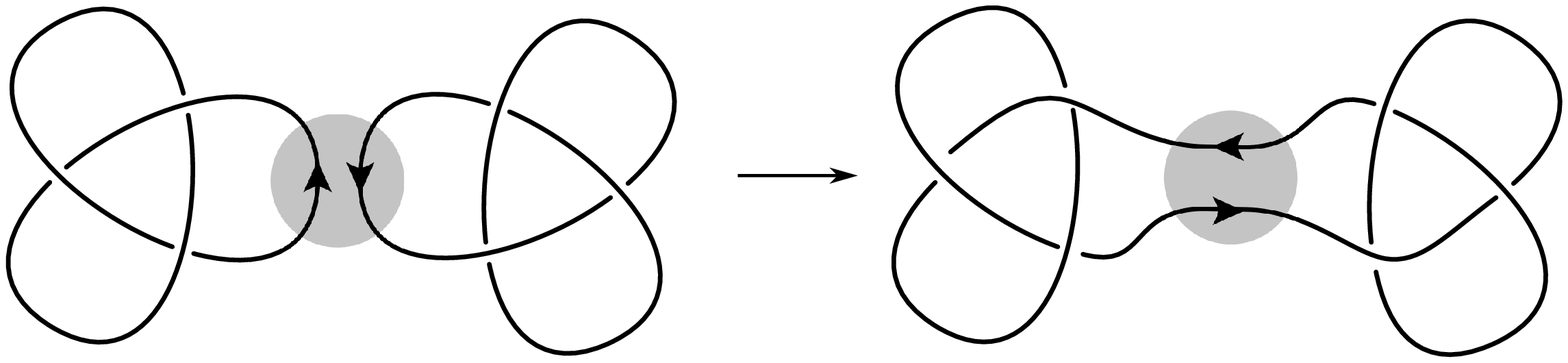}
\includegraphics[scale=0.4]{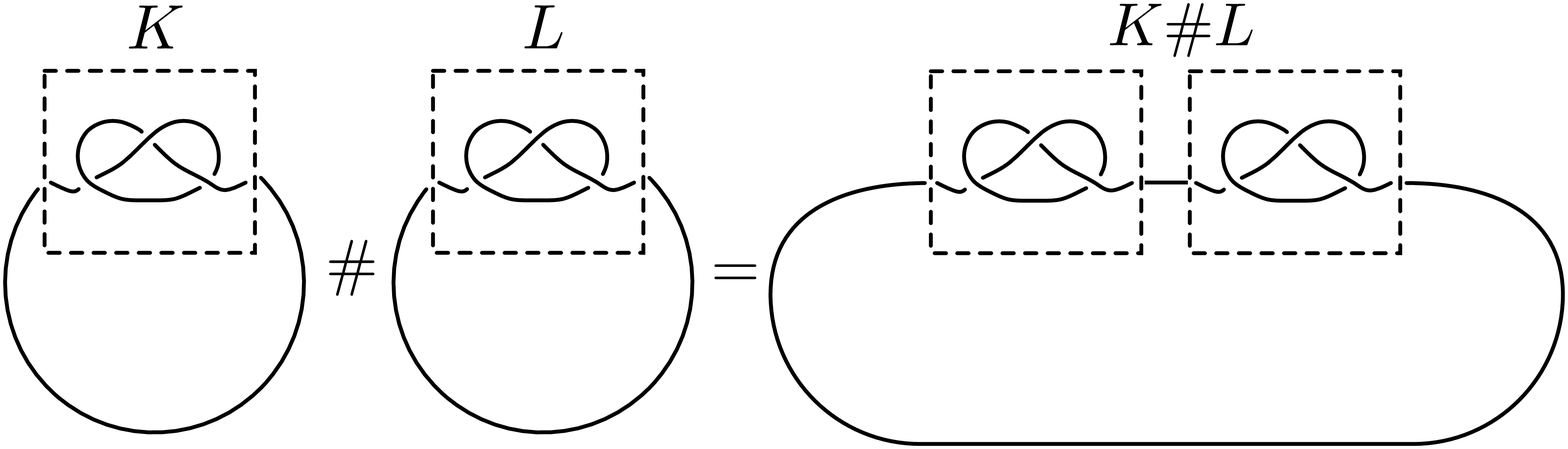}
\caption{Connected sum of knots}
\label{f:consum}
\end{figure}

%\cite[\S1, figures 1.6 and 1.17]{PS96} \cite[p. 25]{CDM}

The {\bf connected sum} $\#$ of {\it oriented} knots is defined in fig. \ref{f:consum}.
The {\bf connected sum} $\#$ of (non-oriented) knots is defined analogously (just ignore the arrows).
Neither operation is well-defined on the set of knots or oriented knots, respectively.
So we denote by $K\# L$ any of the connected sums of knots or oriented knots $K$ and $L$.

%\footnote{More precisely, consider disjoint oriented plane diagrams of the two oriented knots.
%Find a rectangle in the plane where one pair of sides are edges of each knot, but the rectangle is otherwise %disjoint from the knots, and the edges are oriented around the outline of the rectangle in the same direction.
%Now join the two diagrams together by deleting these edges from the knots and adding the edges
%that form the other pair of sides of the rectangle.
%The resulting connected sum diagram inherits an orientation consistent with the orientations
%of the two original diagrams.}

%Не очень понятный текст.
%Не написано, что берутся противоположные стороны прямоугольника.
%Что значит 'rectangle is otherwise disjoint from the knots' ---
%только контур прямоугольника не пересекается с узлами, или его внутренность?
%Почему такой прямоугольник найдется?
%Наконец, как все это связано с рис. \ref{f:consum}?
%Стороны пунктирных прямоугольников на этом рисунке явно не содержат ребер узлов.

%E.g. the knots in \cite[\S1, Figure ]{PS96} are {\it composite}, i.e. are connected sums of some knots.

\begin{assertion}\label{wha-cons} For any
%isotopy classes $K,L,M$ of
oriented knots $K,L,M$ and
%the isotopy class $O$ of
the trivial oriented knot $O$ we have

(a) $K\#O=K$. \quad (b) $K\#L=L\#K$. \quad (c) $(K\#L)\#M=K\#(L\#M)$. \quad

(d) For any non-oriented knots $K,L$ we have $\arf(K\#L)=\arf K+\arf L$.

(The rigorous meaning of (a) is `any connected sum of $K$ and $O$ is isotopic to $K$'.
Analogous rigorous meanings have (b,c) and (d).
See though remark \ref{r:well}.)
\end{assertion}

%associativity

\begin{proof}[Sketch of the proof]
(a) See fig. \ref{f:kok}.

\begin{figure}[h]\centering
\includegraphics[scale=0.7]{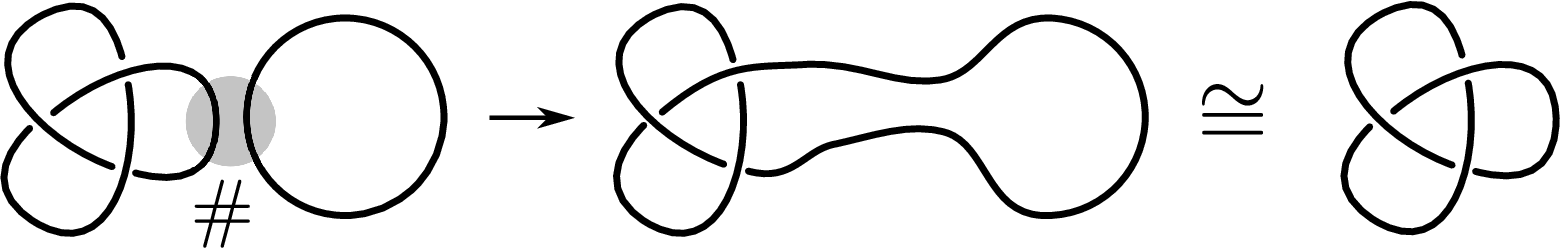}
\caption{Proof of $K\#O = K$}
\label{f:kok}
\end{figure}

\begin{figure}[h]\centering
\begin{tabular}{@{}cc@{}}
\includegraphics[scale=0.4]{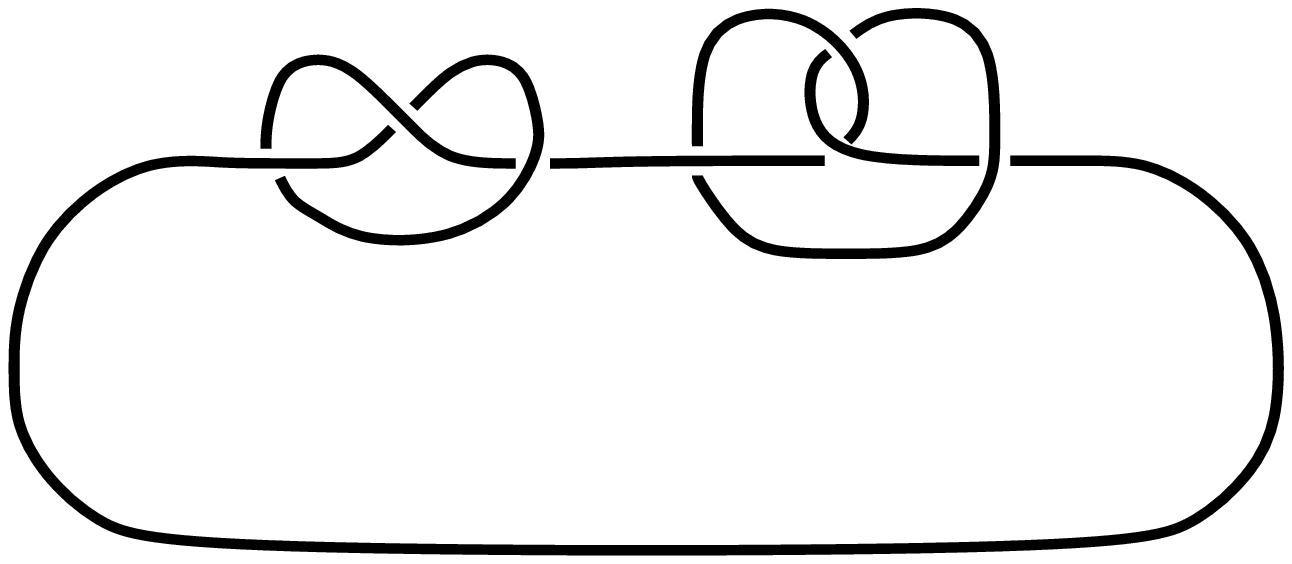}&\\
\includegraphics[scale=0.4]{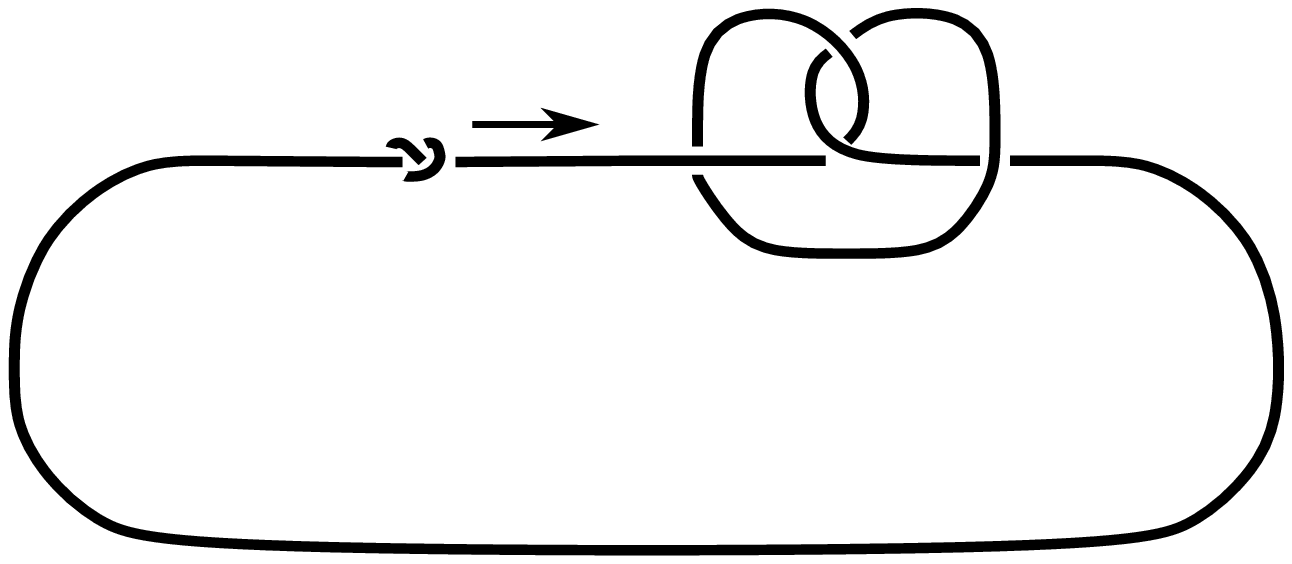}&\\
\includegraphics[scale=0.4]{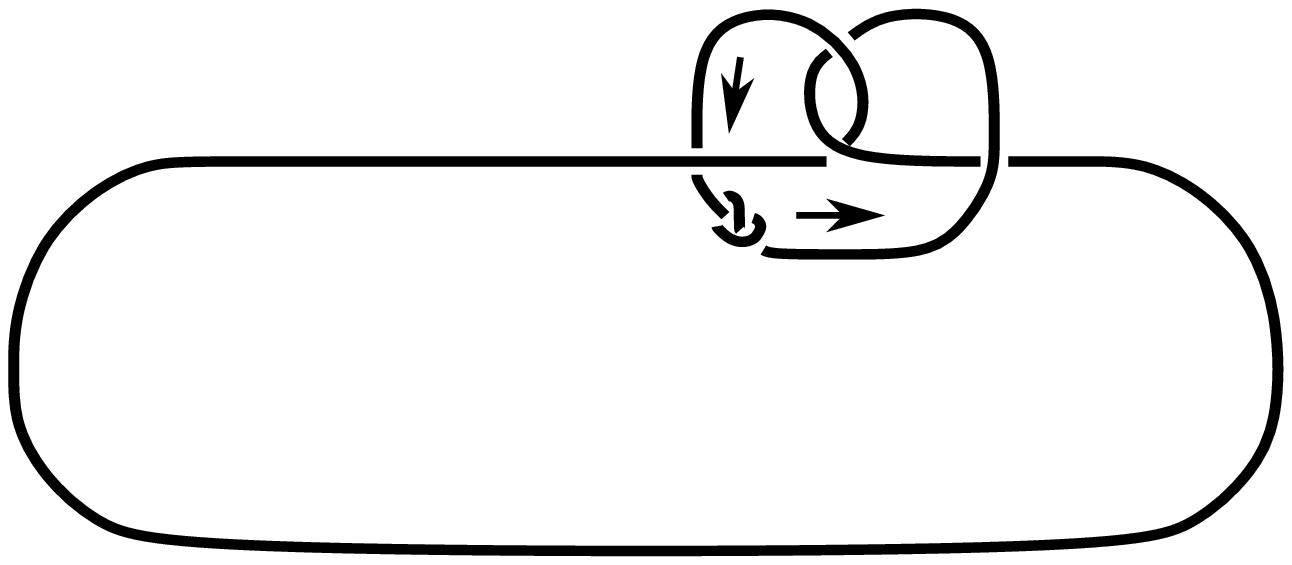}& \qquad \qquad \includegraphics[scale=0.4]{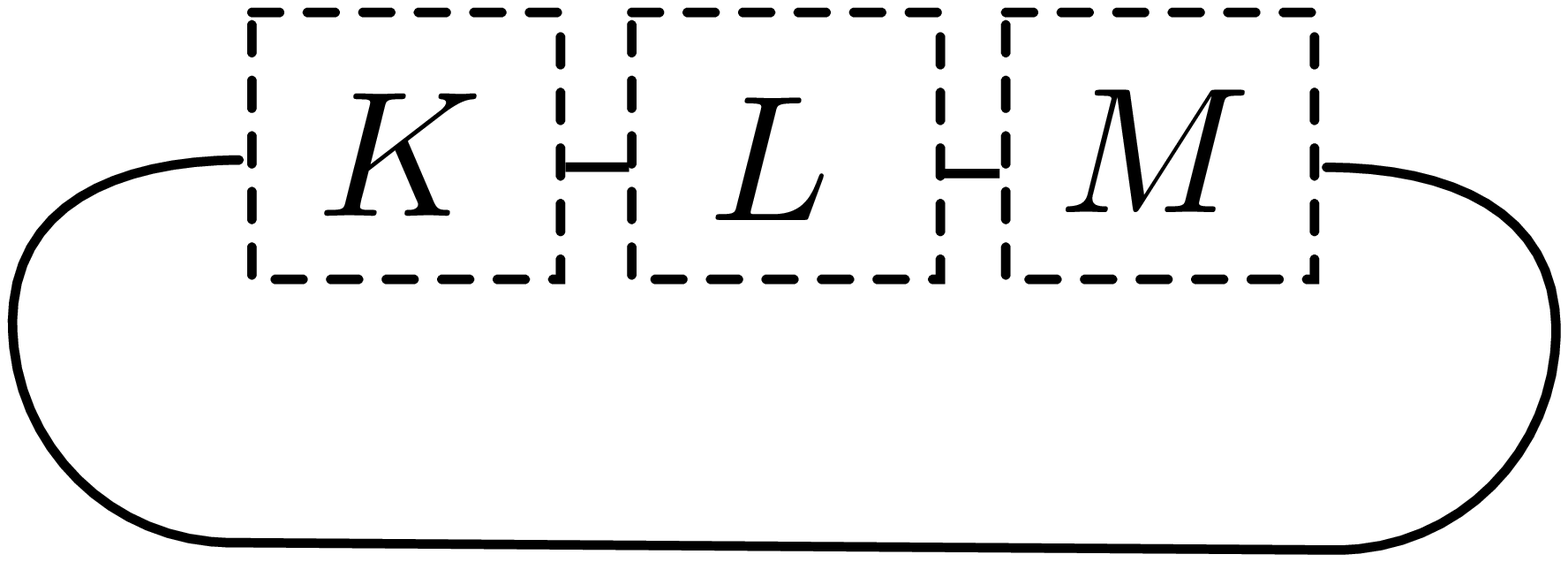}\\
\includegraphics[scale=0.4]{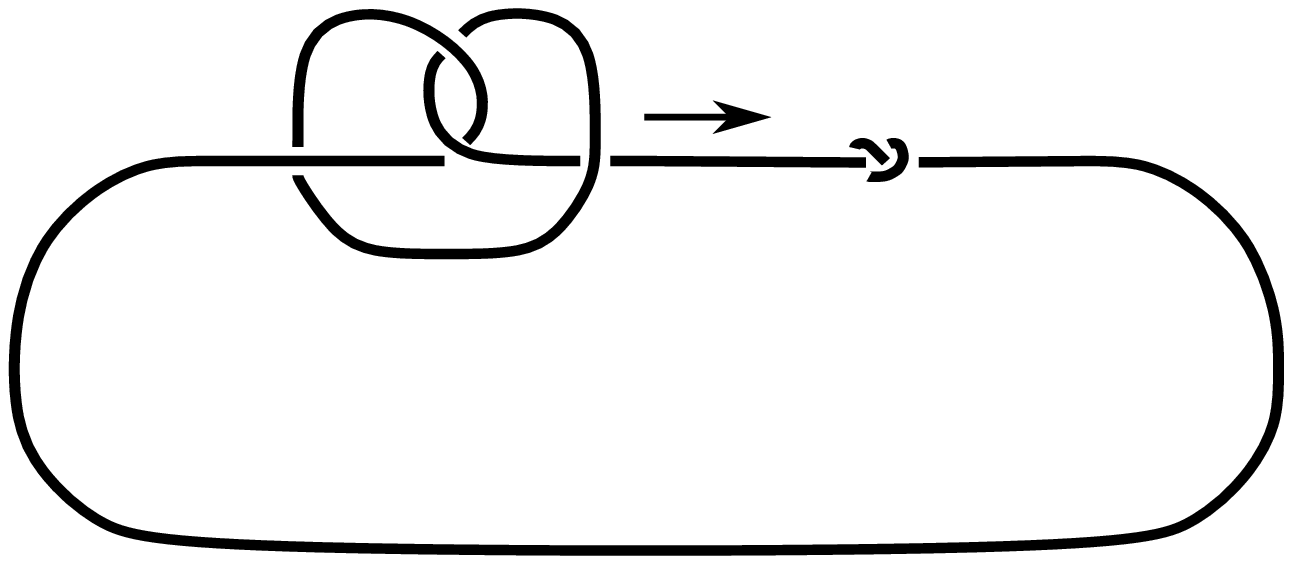}&\\
\includegraphics[scale=0.4]{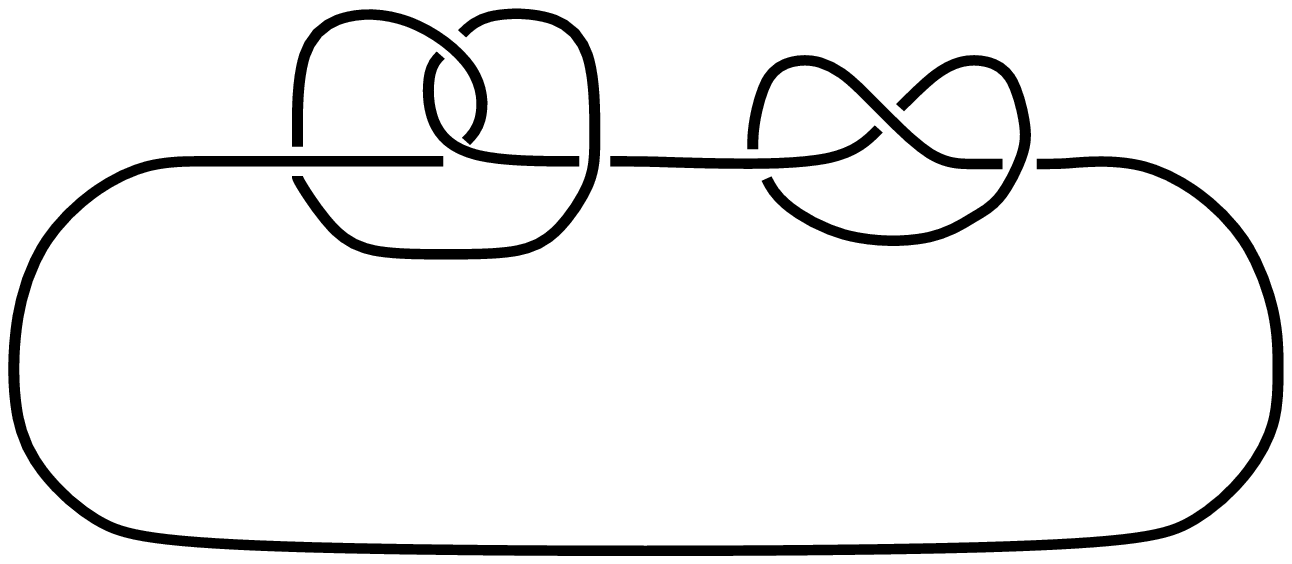}&\\
\end{tabular}
\caption{Proofs of $K\#L = L\#K$ and of $(K\#L)\#M = K\#(L\#M)$}\label{f:comm}
\end{figure}

(b) Take a small knot isotopic to $L$ and push it through $K$, see fig. \ref{f:comm}, left.

(c) The left hand and the right hand side of the equality are isotopic to the knot in fig. \ref{f:comm}, right.

(d) Choose basepoint close to the `place of connection'.
Check that all skew pairs of crossings in $K\#L$ are obtained from the skew pairs of crossings in $K$ and in $L$.
\end{proof}

\begin{remark}\label{r:well} An {\it isotopy class} of a knot is the set of knots isotopic to this knot.
The oriented isotopy class $[K\#L]$ of the connected sum of two oriented isotopy classes $[K],[L]$ of oriented knots $K,L$ is independent of the choices used in the construction, and of the representatives $K,L$ of $[K],[L]$.
Hence the connected sum of oriented isotopy classes of oriented knots is well-defined by $[K]\#[L]:=[K\#L]$, see  \cite[Remark 2.3.a]{Sk15}.
For isotopy classes of non-oriented knots the connected sum is not well-defined \cite{CSK}.
\end{remark}

\begin{theorem}\label{t:consum}
For any isotopy classes $K,L,M$ of knots and the isotopy class $O$ of the trivial oriented knot we have

(a) if $K\#L=O$, then $K=L=O$;

(b) if $K\#L=K\#M$, then $L=M$.
\end{theorem}

The proof is outside the scope of this text, see \cite[Theorem 1.5]{PS96}.
(In this part of \cite{PS96} one needs to replace `knot' by `oriented knot' because of remark \ref{r:well}.)

\begin{figure}[h]\centering
\includegraphics[scale=0.9]{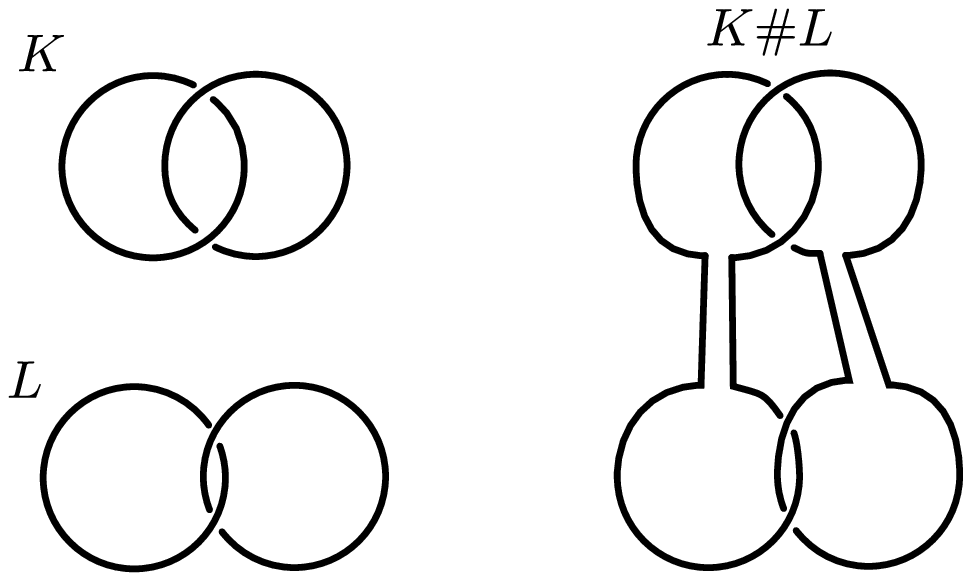}
\caption{Connected sum of links}
\label{f:consuml}
\end{figure}

%fig. 1 of https://arxiv.org/pdf/1704.06501.pdf

The connected sum $\#$ of links (ordered or not, oriented or not) is defined analogously to the connected sum of knots, see fig. \ref{f:consuml}.
This is not a well-defined operation on links, and assertion \ref{wha-csli} shows that this does not
give a well-defined operation on their isotopy classes.
So we denote by $K\# L$ any of the connected sums of $K$ and $L$.

\begin{assertion}\label{wha-consli} (a,b,c) The analogues of assertions \ref{wha-cons}.a,b,c for links are true.

(b) For any non-oriented 2-component links $K,L$ we have $\lk_2(K\#L)=\lk_2 K+\lk_2 L$.
\end{assertion}

\begin{figure}[h]
\centering
\includegraphics[scale=0.5]{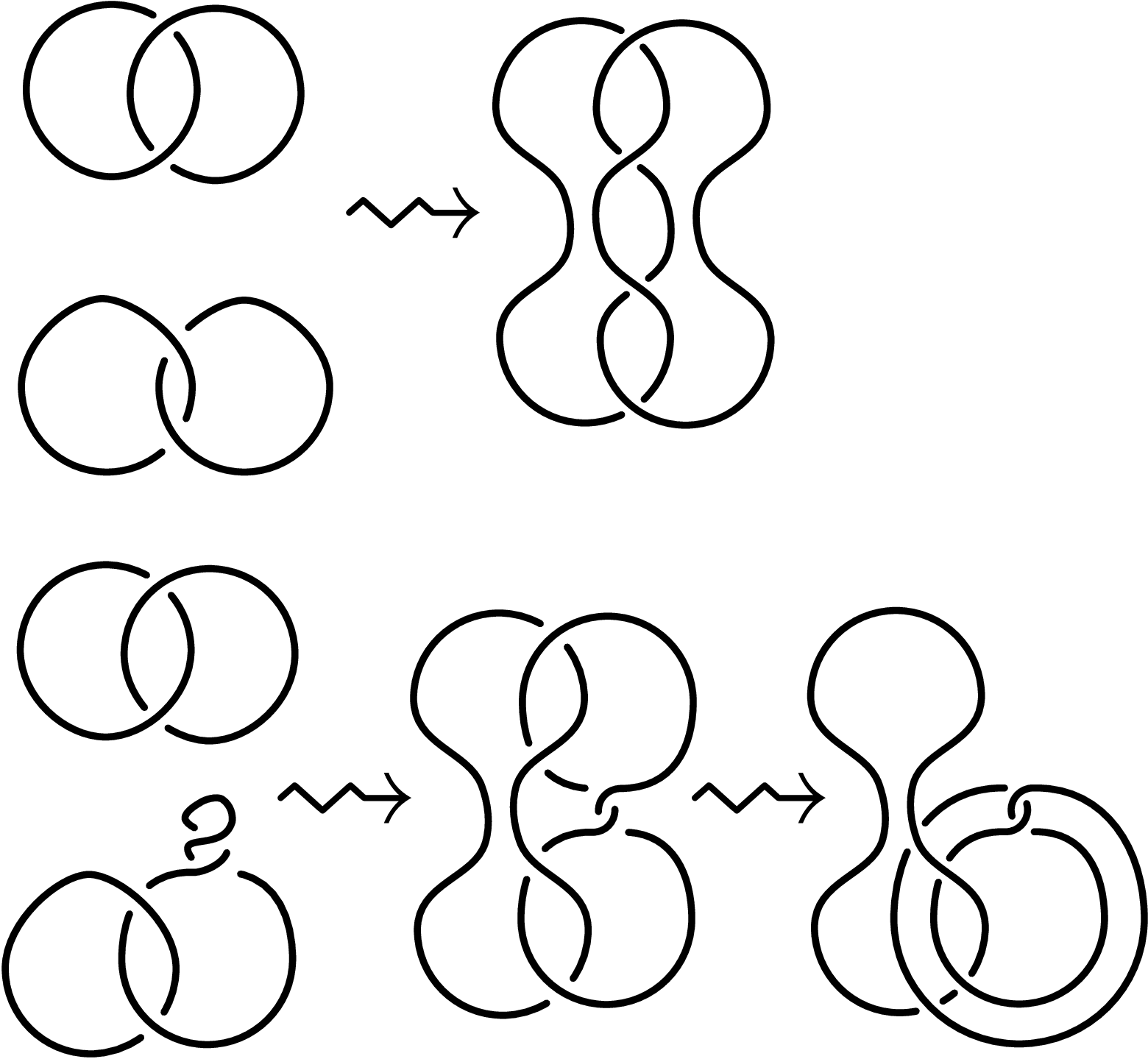}
\caption{Connected sum of isotopy classes of ordered links is not well-defined}
\label{f:ryabichev}
\end{figure}

\begin{remark}\label{wha-csli} There are two isotopic pairs $(K,L)$ and $(K',L')$ of 2-component links such that some connected sums  $K\#L$ and $K'\#L'$ are not isotopic.
(The links could be ordered or non-ordered, oriented or non-oriented, so that we have 4 statements.)
As an example of non-ordered pairs we can take equal links consisting of a trefoil and an unknot in disjoint cubes, cf. \cite[Figure 3.16]{PS96}.
For an example of ordered pairs see \cite{As}.
Fig. \ref{f:ryabichev} presents an alternative example suggested by A. Ryabichev.
\end{remark}

\comment

\begin{figure}[h]\centering
\includegraphics[scale=0.4]{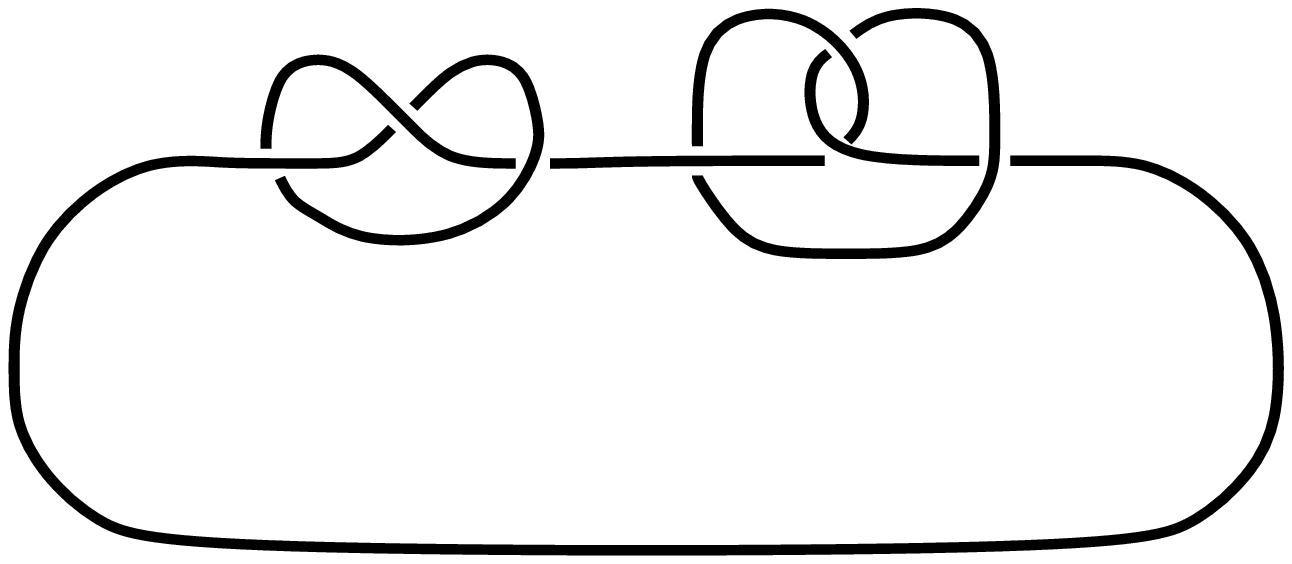}\\
\includegraphics[scale=0.4]{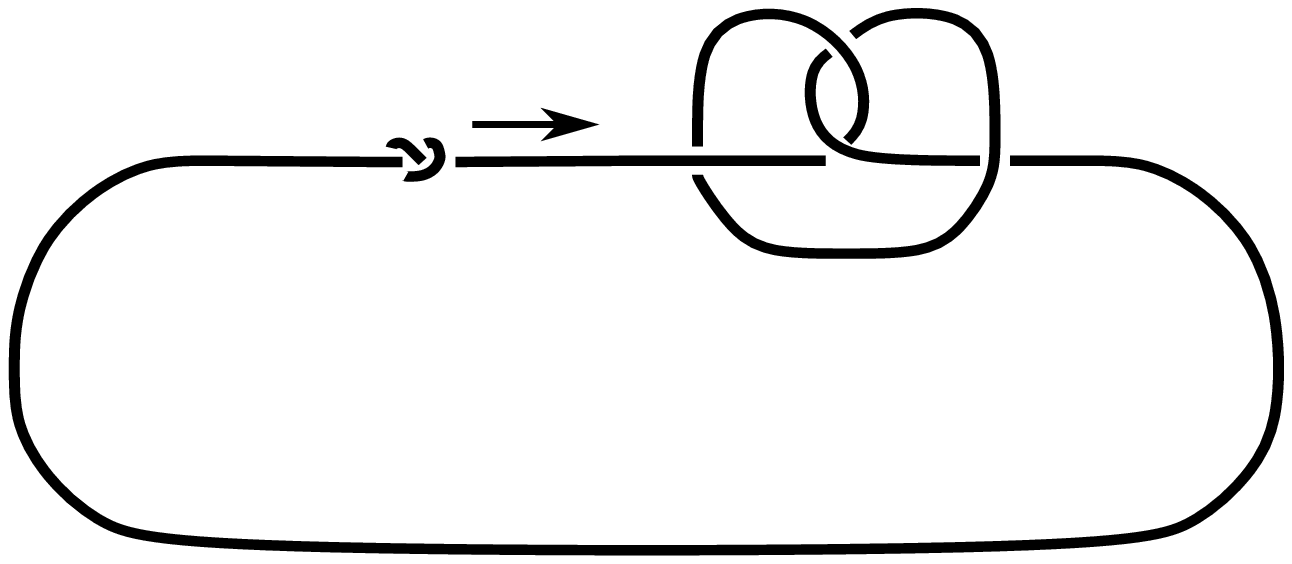}\\
\includegraphics[scale=0.4]{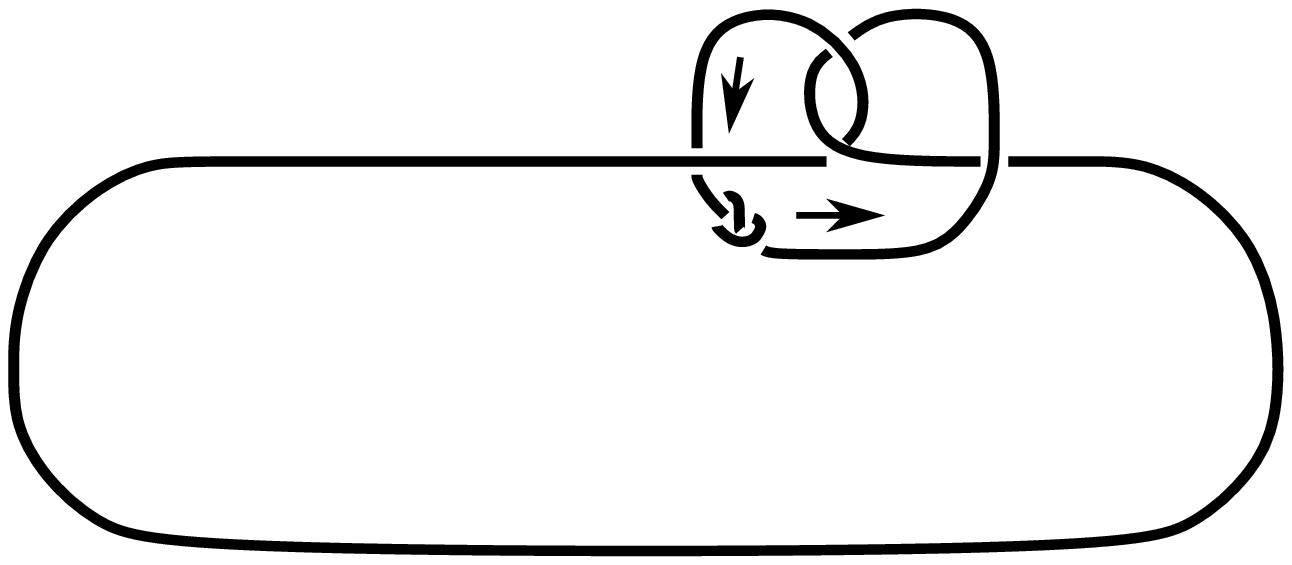}\\
\includegraphics[scale=0.4]{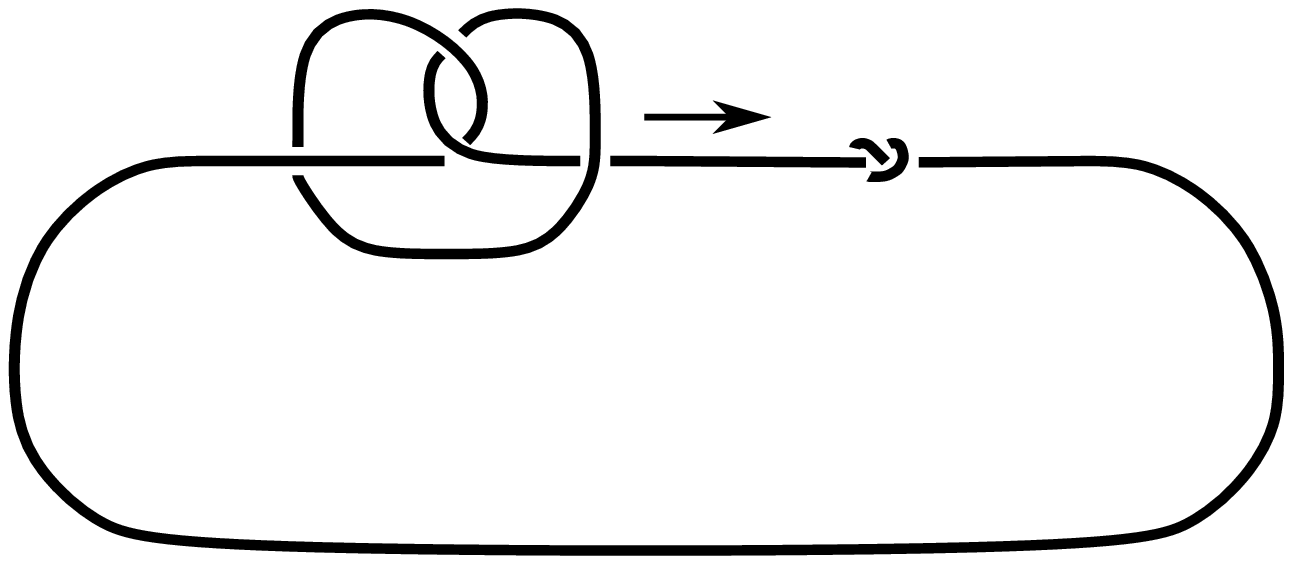}\\
\includegraphics[scale=0.4]{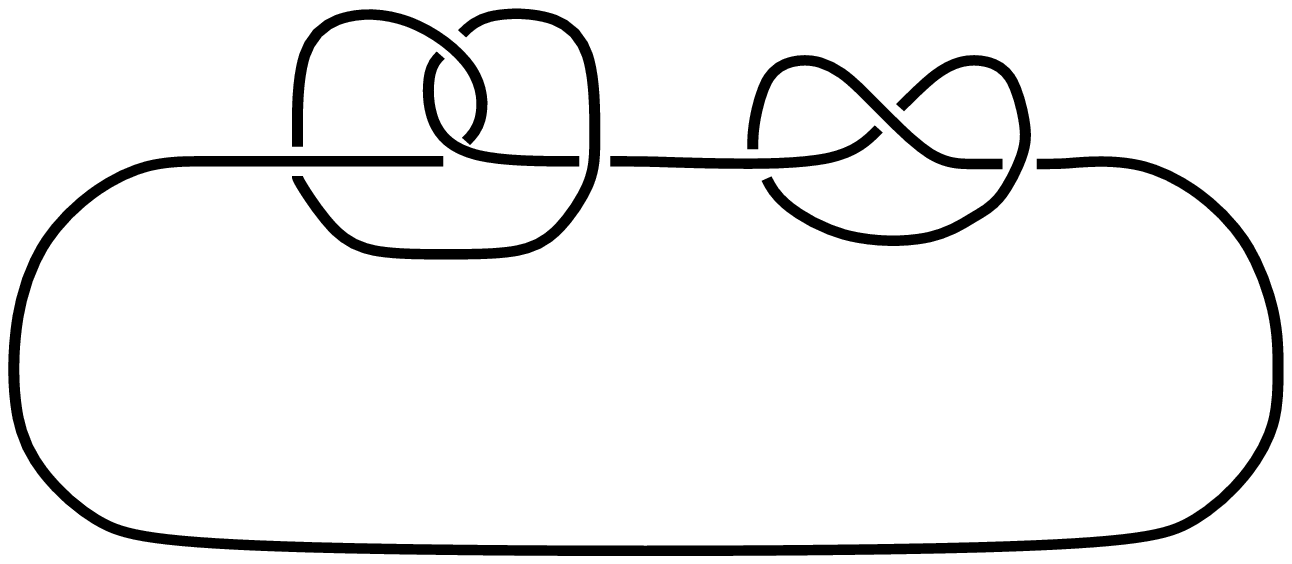}\\
\end{figure}

\begin{figure}[h]\centering
\includegraphics[scale=0.4]{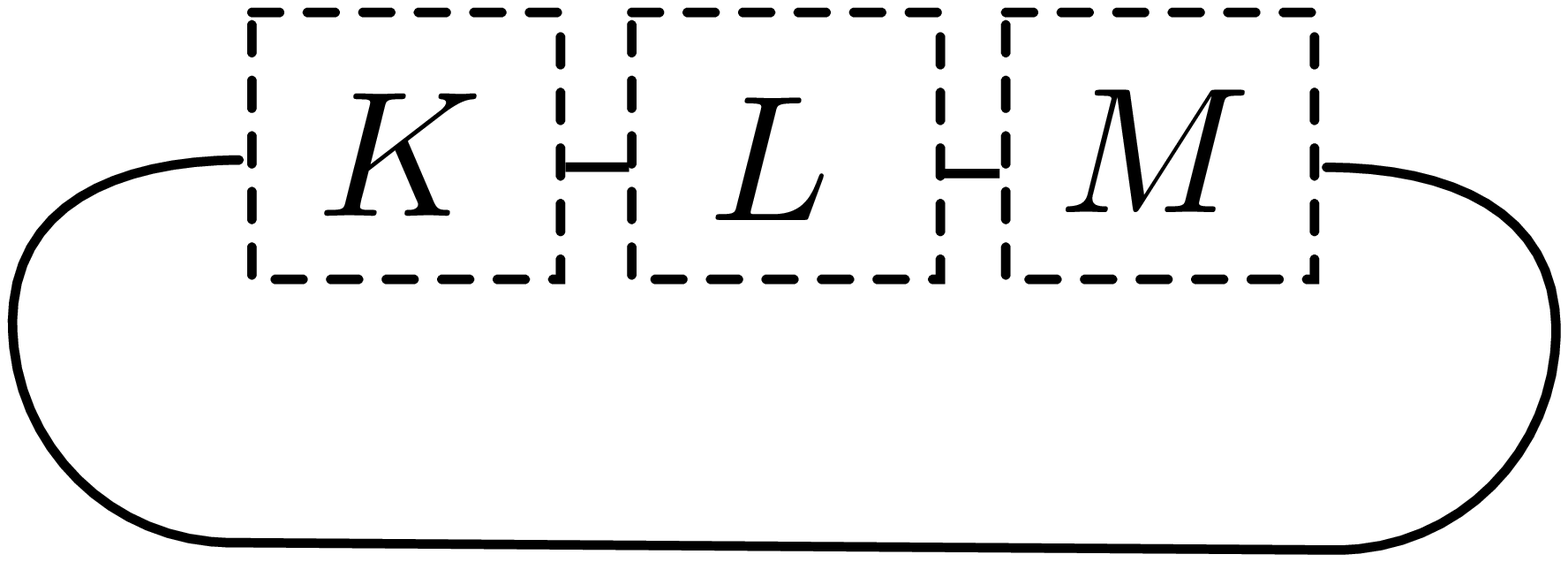}
\caption{Proof of $(K\#L)\#M = K\#(L\#M)$}
\label{f:asso}
\end{figure}

\endcomment

\section{The Gauss linking number via plane diagrams}\label{0rapl}

\begin{figure}[h]\centering
\includegraphics[scale=0.7]{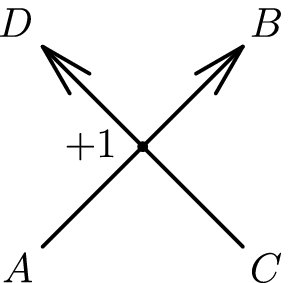}\quad \includegraphics[scale=0.7]{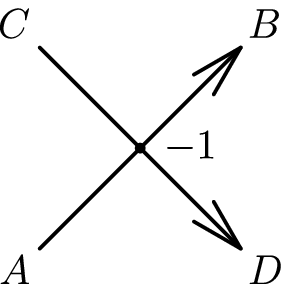}
\caption{The sign of intersection point}
\label{f:sign}
\end{figure}

Let $(\overrightarrow{AB},\overrightarrow{CD})$ be an ordered pair of vectors (oriented segments) in the plane intersecting at a point $P$.
Define {\bf the sign} of the pair to be $+1$ if $ABC$ is oriented clockwise and to be $-1$ otherwise (fig. \ref{f:sign}).

The {\bf linking number} $\lk$ of the plane diagram of an oriented 2-component link is the sum of signs at all those crossing points on the diagram at which the first component passes above the second component.
At every crossing point the {\it first} (the {\it second}) vector is the oriented edge of the first (the second) component.

\begin{pr}\label{e:lk} Find the linking number for (some plane diagram of) the Hopf link and pairs of Borromean rings, for your choice of orientation on the components.
\end{pr}

\begin{lemma}\label{p:lk} The linking number is preserved under Reidemeister moves.
\end{lemma}

The proof is analogous  to lemma  \ref{l:lk2}.
It suffices to check that the signs of all crossing points do not change.

By Lemma \ref{p:lk} the {\bf linking number} of an oriented 2-component link (or of its isotopy class)
is well-defined by setting it to be the linking number of any plane diagram of the link.

The {\it absolute value of the linking number} of a (non-oriented) 2-component link (or of its isotopy class)
is well-defined by taking any orientations on the components.

We shall use without proof the following {\it Triviality lemma:} for any two closed oriented polygonal lines in the plane whose vertices are in general position the sum of signs of their intersection points is zero.
For a discussion and a proof see \S1.3 `Intersection number for polygonal lines in the plane' of \cite{Sk18}, \cite{Sk}.

%\begin{lemma}[Triviality]\label{l:triv}  \end{lemma}

\begin{assertion}\label{p:lk-pr} (a) Switching the components of a link negates the linking number.

(b) Reversing the orientation of either of the components negates the linking number.

(c) There is an oriented 2-component link whose linking number is $-5$.

(d) For any of the connected sums $K\#L$ of oriented 2-component links $K,L$ we have $\lk(K\#L)=\lk K+\lk L$.

(e) There is a 2-component link which is not isotopic to the trivial link but which has zero linking number.
\end{assertion}

Part (e) is proved using \emph{Alexander-Conway polynomial}, see \S\ref{0con}.

%For any integer $n$ there is an oriented link of linking number $n$.

\begin{figure}[h]\centering
\includegraphics[scale=0.6]{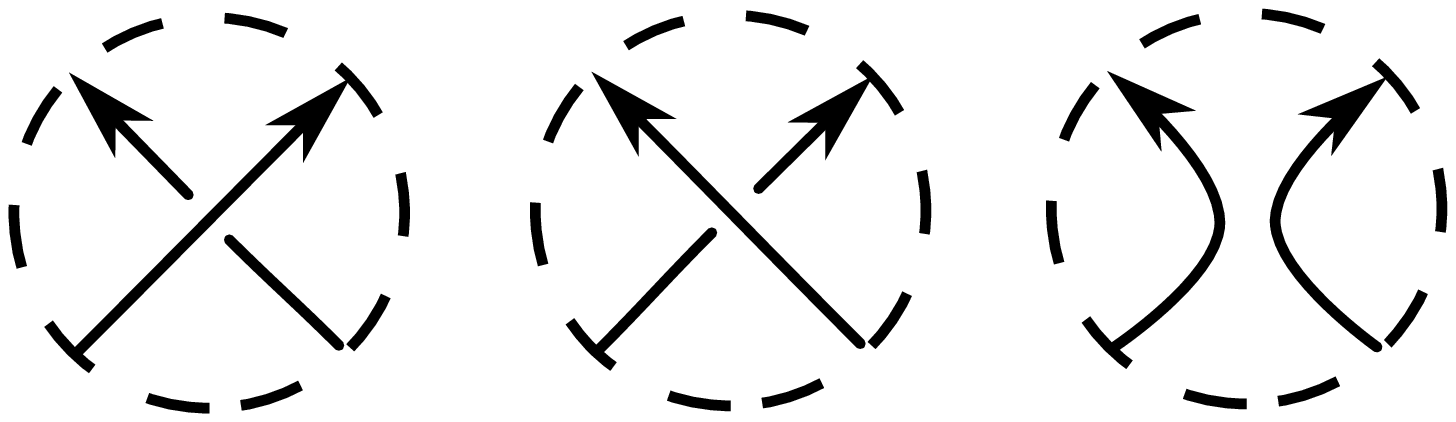}
\caption{Oriented links $K_+,K_-,K_0$}
\label{f:skein}
\end{figure}

\begin{theorem}\label{con-lk}
There is a unique integer-valued isotopy invariant $\lk$ of oriented 2-component links that assumes value 0 on the trivial link and such that for any links $K_+$ and $K_-$
whose plane diagrams differ by a crossing change of a crossing $A$ as shown in fig. \ref{f:skein}
$$\lk K_+-\lk K_- =
\begin{cases}1 &A\text{ is the crossing of different components;} \\
0 &A\text{ is the self-crossing of one component.}\end{cases}$$
\end{theorem}

The proof is analogous to theorem \ref{con-lk2}.

%\cite[2.3.1]{CDM}

\section{The Casson invariant}\label{0cas}

The {\bf sign} of a crossing point of an oriented plane diagram of a knot is defined after figure \ref{f:sign};
the first (the second) vector is the vector of overcrossing (of undercrossing).
Clearly, the sign is independent of the orientation of the diagram, and so is defined for non-oriented diagram.

The {\bf sign} of a $P$-skew pair of crossing points in a plane diagram of a knot (for any basepoint $P$) is the product of the signs of the two crossing points.

The {\it $P$-Casson invariant} of a plane diagram is the sum of signs over all $P$-skew pairs of crossing points.

\begin{pr}\label{e:con-casca} (a) Same as problem \ref{e:con-cas}.b for the Casson invariant.

(b) Draw a plane diagram of a knot and a basepoint $P$ such that $P$-Casson invariant is $-5$.
\end{pr}

\begin{lemma}\label{p:con-casca-ef} (a,b) Same as lemma \ref{p:con-cas-ab}.a,b for the Casson invariant.
\end{lemma}

Hence the {\bf Casson invariant} (number) $c_2$ of a plane diagram, of a knot, or even of isotopy class of a knot, is well-defined by setting it to be the $P$-Casson invariant of any plane diagram of the knot for any basepoint $P$.

\begin{assertion}\label{con-casca} (a,b) Same as assertions \ref{wha-cons}.d and \ref{a:arf-pr} for the Casson invariant.
\end{assertion}

Part (b) is proved using \emph{Alexander-Conway polynomial}, see \S\ref{0con}.

%Find and prove the analogue of lemma \ref{p:arf-pr} for the Casson invariant.
%(a) The Casson invariants of inverse knots are the same.
%(c) The Casson invariant is uniquely defined by the conditions of Theorem \ref{con-thm}.
%(a) Are the Casson invariants of mirror-symmetric knots the same or the opposite?

Denote by $D_+,D_-,D_0$ any three diagrams of oriented (knots or) links differing as shown in fig. \ref{f:skein} (for a convention on figures see caption to fig. \ref{f:crossing}).
We also denote by $K_+,K_-,K_0$ any three links who have diagrams $D_+,D_-,D_0$.

\begin{theorem}\label{con-thm} There is a unique integer-valued isotopy invariant $c_2$ of (non-oriented) knots that assumes value 0 on the trivial knot and such that for any knots $K_+$ and $K_-$ whose plane diagrams differ as shown in fig. \ref{f:skein}
$$c_2(K_+)-c_2(K_-)=\lk K_0.$$
%Here $K_+,K_-$ are knots from fig. \ref{f:crossing}
%and $K_0$ is the oriented link obtained by using some compatible orientations on $K_+,K_-$;
(Observe that $K_0$ has to be a 2-component  link; the number $\lk K_0$ is well-defined because change of the orientation on both components of an oriented link does not change the linking number.)
\end{theorem}

%This invariant is called {\it Casson invariant} and its reduction modulo 2 is called {\it Arf invariant}.
%In this section the existence part of Theorem \ref{con-thm} can be used without proof.

%\newpage
\section{Alexander-Conway polynomial}\label{0con}

Section \ref{0con} only uses the material of \S\ref{0isot}, \S\ref{0isoli} and \S\ref{s:ori}
(except that assertions \ref{p:conzer}.bc use  \S\ref{0rapl} and \S\ref{0cas}).

\begin{pr}\label{t:a3} (a) There is a unique mod2-valued isotopy invariant $\arf$ of oriented 3-component links that assumes value 0 on the trivial link and for which
$$\arf K_+ -\arf K_- =\begin{cases}\lk\phantom{}_2 K_0 &\text{at the crossing point different components cross each other;} \\
0 &\text{at the crossing point one component crosses itself.}\end{cases}$$
(Here $\lk_2 K_0$ is defined because $K_0$ is a 2-component link.)\footnote{This assertion is a particular case of mod2 version of theorem~\ref{con-ale}; it would be interesting to obtain a direct proof because such a proof could illuminate an idea of proof of theorem~\ref{con-ale} by presenting the idea in the simplest non-trivial situation.
\newline
Theorem \ref{con-thma} is the analogue of this assertion for 1-component links (knots).
The definition of $\arf$ given in \S\ref{0arf} applies to knots only and here the point is to extend it to 3-component links.}

(b) Assuming the existence of the invariant $\arf$ from (a), calculate (for your choice of orientation on the components) the $\arf$ invariant of the Borromean rings.
\end{pr}

%The reader is recommended to prove first mod2-analogues of the following results and problems.

\begin{theorem}\label{con-ale} There is a unique infinite sequence $c_{-1}=0,c_0,c_1,c_2,\ldots$
of $\Z$-valued isotopy invariants of oriented non-ordered links such that

$\bullet$ $c_0=1$ and $c_1=c_2=\ldots=0$ for the trivial knot;

$\bullet$ for any $n\ge0$ and links $K_+,K_-,K_0$ from fig.~\ref{f:skein} we have
$$c_n(K_+)-c_n(K_-)=c_{n-1}(K_0).$$
%where $K_0$ is $K_0$ from fig.~\ref{f:skein} with some ordering??? of the components.
\end{theorem}

%Actually two versions of theorem~\ref{con-ale} are equivalent.
Proofs of the existence in assertion \ref{t:a3}.a and in theorem \ref{con-ale} are outside the scope of this text.
We shall use the existence without proof.\footnote{It is not clear whether the statement in  \cite[\S2.3.1]{CDM} involves ordered or non-ordered links.
We deduce the stronger version (for non-ordered links) from the weaker version (for ordered links) in \S\ref{s:app}.}
See a proof in \cite{Al28}, \cite[\S3-\S5]{Ka06'}, \cite{Ka06}, \cite[\S5.5 and \S5.6]{Ma18}, \cite{Ga20}.
For a relation to proper colorings see \cite[\S6]{Ka06'}.

%You can earn a plus-mark (plus-sign) for proving the uniqueness, and solve other problems assuming the existence.
%his book "Knots and Physics".
%There he mimics the R-matrix and Yang-Baxter equation, but it runs quite elementary.

The polynomial $C(K)(t):=c_0(K)+c_1(K)t+c_2(K)t^2+\ldots$ is called the {\it Conway polynomial}, see assertion \ref{p:conzer}.e.
Introduction of this polynomial allows to calculate all the invariants $c_n$ as quickly as one of them.
The formula in theorem~\ref{con-ale} is equivalent to
$$
C(K_+)-C(K_-)=tC(K_0).
$$

%(c) The Conway polynomial is uniquely defined by the conditions of Theorem \ref{con-ale}.
%{\it There is a unique $\Z[t]$-valued invariant $C$
%of oriented links  that assumes value 1 on the trivial knot and for which $C(K_+)-C(K_-)=tC(K)$.}

%It is not too complicated to calculate each invariant directly using the definition.
%However it is more effective to calculate Conway polynomial which gives all invariants $c_n$
%from theorem~\ref{con-ale} (see solution of problem~\ref{con-cal}).

%In this section Theorem \ref{con-ale} can be used without proof.

\begin{pr}\label{con-cal}
%Assuming the existence of the Conway polynomial,
Calculate the Conway polynomial of the following links (for your choice of orientation on the components).

(a) the trivial link with 2 components; \quad
(b) the trivial link with $n$ components; \quad

(c) the Hopf link; \quad
(d) the trefoil knot; \quad
(e) the figure eight knot; \quad

(f) the Whitehead link; \quad
(g) the Borromean rings; \quad
(h) the $5_1$ knot.
\end{pr}

\begin{assertion}\label{p:conzer} (a) We have $c_0(K)=1$ if $K$ is a knot and $c_0(K)=0$ otherwise (i.e. if $K$ has more than one component).

(b) For a knot $K$ we have $c_{2j+1}(K)=0$ and $c_2$ is the Casson invariant.

(c) For a 2-component link $K$ we have $c_{2j}(K)=0$ and $c_1$ is the linking coefficient.

(d) For a $k$-component link $K$ we have $c_j(K)=0$ if either $j\le k-2$ or $j-k$ is even.

(e) For every knot or link all but a finitely many of the invariants $c_n$ are zeroes.
\end{assertion}

\begin{assertion}\label{con-pro} (a) Change of the orientations of all components of a link (in particular, change of the orientation of a knot) preserves the Conway polynomial.

%the attached copy of

(b) \cite[2.3.4]{CDM} There is a 2-component link such that change of the orientation of its one component changes the degree of the Conway polynomial (so this change neither preserves nor negates the Conway polynomial).

(c) For any of the connected sums $K\#L$ of knots $K,L$ we have $C(K\#L)=C(K)C(L)$.
\end{assertion}

%(b) Does change of the orientation of any link preserve the Conway polynomial?

%(b) The Conway polynomials of mirror-symmetric knots are the same.

%\begin{figure}[h]\centering
%\includegraphics[width=2.7cm]{pict1.33.eps}
%\caption{A split link (Fall 2013-2-11)}
%\label{f:split}
%\end{figure}

\begin{assertion}\label{con-spl}
A link is \emph{split} if it is isotopic to a link whose components are contained in disjoint balls.

%3-dimensional

%(a) The link of fig. \ref{f:split} is split.
(a) Neither Hopf link nor Whitehead link nor Borromean rings link is split.

(b) The linking coefficient of a split link is zero.

(c) The Conway polynomial of a split link is zero.
%\cite[Ex. 4 in p. 64]{CDM}.
\end{assertion}

%This is proved using {\it the Conway poynomial} in \S\ref{0con}.

%\begin{pr}\label{con-motiv} If $C$ is a $\Z[t]$-valued invariant of oriented links that
%assumes value 1 on the trivial knot and $a,b,c\in\Z[t]$ are such that $aC(K_+)+bC(K_-)=cC(K)$,
%then $a=-b=1$ and $c=t$.

%Hint: analogously to \cite[\S3.1]{PS96}.
%\end{pr}

%See also \cite{A}.

\comment

{t:a3}
(b) There is a unique mod2-valued isotopy invariant $a_3$ of oriented 2-component links that assumes value 0 on the trivial link and for which
$$a_3(K_+)-a_3(K_-)=\arf K_0.$$
(Here $\arf K_0$ is defined because $K_0$ is either a knot or a 3-component link.)

(c) There is a unique mod2-valued isotopy invariant $a_3$ of oriented 4-component links that assumes value 0 on the trivial link and for which
$$a_3(K_+)-a_3(K_-)=\begin{cases}\arf K_0 &\text{at the crossing point different components cross each other;} \\
0 &\text{at the crossing point one component crosses itself.}\end{cases}$$
(Here $\arf K_0$ is defined because $K_0$ is a 3-component link.)

\begin{pr}\label{con-cal2}
%Assuming the existence of the invariants $\arf,a_3$ from problem \ref{t:a3},
Calculate (for your choice of orientation on the components)

(a) the $\arf$ invariant of the Borromean rings;

(b,c,d) the $a_3$ invariant of the Hopf link, of the Whitehead link, and of {\it 4-Borromean rings},
i.e. of any link of your choice for which every 3-component sublinks are isotopic to the trivial link,
but the entire link is not isotopic to the trivial link.
\end{pr}

 %{\bf \ref{con-cal2}.} {\it Answers:} (a, b) 0; (c) 1 (independently of the choice of orientation).

\begin{figure}[h]\centering
\includegraphics[scale=0.6]{4-borromean.eps}
\caption{4-Borromean rings}
\label{f:4-borromean}
\end{figure}

\endcomment

%\newpage
\section{Vassiliev-Goussarov invariants}\label{0vas}

Section \ref{0vas} only uses the material of \S\ref{0isot}, \S\ref{0isoli} and \S\ref{s:ori}
(except that problem \ref{vas-pro}.2 uses \S\ref{0cas}).

An (oriented) {\it singular knot} is a closed oriented polygonal line in $\R^3$ whose only self-intersections are double points which are not vertices.
%PL map  $f:S^1\to\R^3$ whose only self-intersections are {\it transversal} double points \cite[4.1]{PS96}.
Two singular knots are {\it isotopic} if there is an orientation preserving PL homeomorphism $h:\R^3\to\R^3$
carrying the first singular knot to the second one, and the orientation on the first singular knot to the orientation on the second one.
Denote by $\Sigma$ the set of all isotopy classes of singular knots.

A {\it chord diagram} is a cyclic word of length $2n$ having $n$ letters, each letter appearing twice.
A chord diagram is depicted as an oriented circle with a collection of chords, cf. \cite[\S1.5]{Sk20}.
For a singular knot $K$ denote by $\sigma(K)$ the following chord diagram.
Move uniformly along the oriented circle and for any point $A$ on the circle take the `corresponding' point $f(A)$ on $K$.
Join by a chord each pair of points on the circle corresponding to the intersection point of $K$ \cite[4.8]{PS96}, \cite[3.4.1]{CDM}.\footnote{In other words, take a PL map $f:S^1\to\R^3$ of the circle whose image is $K$. Take  a chord $XY$ for each pair of points $X,Y$ such that $f(X)=f(Y)$.
A chord diagram should not be confused with the {\it Gauss diagram} (of a projection) of a (non-singular) knot $g:S^1\to\R^3$ which is the (somehow oriented) chord diagram of the composition of the projection $\R^3\to\R^2$ and $g$ \cite[4.8]{PS96} \cite[1.8.4]{CDM}.}

\begin{figure}[h]\centering
\includegraphics[scale=1]{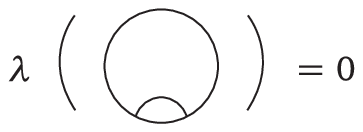}\qquad\qquad\includegraphics[scale=1]{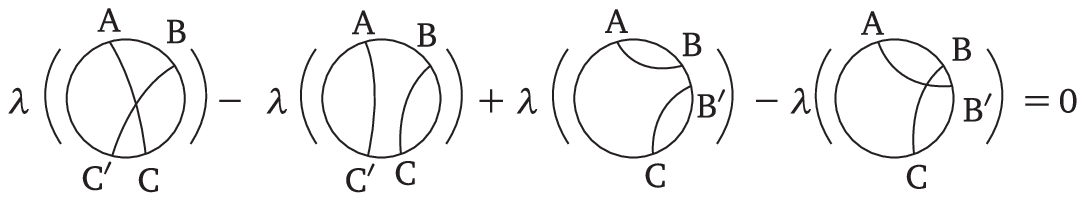}
\caption{The 1-term and 4-term relations}
\label{f:term}
\end{figure}

%\cite[(4.5),(4.6)]{PS96}

\begin{theorem}\label{vas-kon} Assume that $n\ge0$ is an integer and $\lambda:\delta_n\to\R$ a map from the set $\delta_n$ of all chord diagrams that have $n$ chords.
% assuming only a finite number of non-zero values.
The map $\lambda$ satisfies the 1-term and the 4-term relations from fig. \ref{f:term} if and only if
there exists a map $v:\Sigma\to\R$ (i.e. an invariant of singular knots) such that

\begin{figure}[h]\centering
\includegraphics[scale=1]{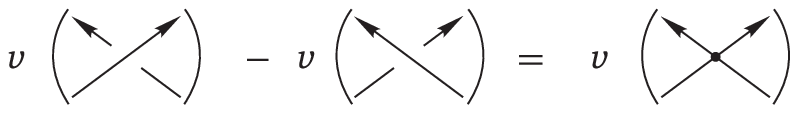}
\caption{The Vassiliev skein relation, notice the difference with fig. \ref{f:skein}}
\label{f:vskein}
\end{figure}

(1) The Vassiliev skein relation from fig. \ref{f:vskein} holds,

%For any singular knots $K_+,K_-$ and $K^0$ from fig. \ref{f:vskein} we have
%$$v(K_+)-v(K_-)=v(K^0).$$ \cite[(4.1)]{PS96}

($2_n$) $v(K)=0$ for each singular knot that has more than $n$ self-intersection points, and

(3) $v(K)=\lambda(\sigma(K))$ for each singular knot $K$ that has exactly $n$ self-intersection points.
\end{theorem}

The proof is outside the scope of this text.

As far as I know, the Vassiliev-Kontsevich theorem \ref{vas-kon} was never stated in this form, which is short and convenient for calculation of the invariants (although this form was implicitly used when the invariants were calculated).
I am  grateful to S. Chmutov for confirmation that theorem \ref{vas-kon} is correct and is indeed equivalent to the standard formulation of the Vassiliev-Kontsevich theorem, see e.g. \cite[Theorem 4.2.1]{CDM}, cf. \cite[Theorem 4.12]{PS96}.

A map $v:\Sigma\to\R$ such that (1) holds is called a {\it Vassiliev-Goussarov invariant}.
If additionally ($2_n$) holds, then $v$ is called a {\it invariant of order at most $n$}.
The following assertion shows that the map $v$ of theorem \ref{vas-kon} is unique up to Vassiliev-Goussarov invariant of order at most $n-1$.

\begin{assertion}[{\cite[Proposition 3.4.2]{CDM}}] The difference between maps $v,v':\Sigma\to\R$ satisfying to (1), ($2_n$) and (3), satisfies to (1) and ($2_{n-1}$).
\end{assertion}

\begin{pr}\label{vas-pro} (a) Prove the `if' part of theorem \ref{vas-kon}.

(0),(1),(2) Prove the `only if' part of theorem \ref{vas-kon} for $n=0,1,2$.

\emph{Hint.} For $n=2$ use theorem \ref{con-thm}.

(The `only if' part of theorem \ref{vas-kon} for $n=3$ could be proved using the coefficient of $h^3$ in $J(e^h)$, where $J$ is the Jones polynomial in $t$-parametrization \cite[2.4.2, 2.4.3]{CDM} \cite[(4.6)]{PS96}.)
\end{pr}

In the remaining problems use (the `only if' part of) theorem \ref{vas-kon} without proof.
Assertion `$v(K)=x$ for any singular knot $K$ whose chord diagram is $A$' is shortened to `$v(A)=x$'.

\begin{pr}\label{vas-cal2} (a) There exists a unique Vassiliev-Goussarov invariant $v_2:\Sigma\to\R$ of order at most 2 such that $v_2(O)=0$ for the trivial knot $O$ and $v_2(1212)=1$.
(Here (1212) is the `non-trivial' chord diagram with 2 chords, see \cite[Figure 4.4]{PS96}, 3rd diagram of the first line.)

\emph{Hint.} This follows from theorem \ref{con-thm}, but try to deduce this from theorem \ref{vas-kon}.

(b,b',c,d) Calculate $v_2$ for the (arbitrary oriented) right trefoil, left trefoil, figure eight knot and the $5_1$ knot.
\end{pr}

\begin{pr}\label{vas-cal3} (a) There exists a unique Vassiliev-Goussarov invariant $v_3:\Sigma\to\R$ of order at most 3 such that $v_3(O)=0$ for the trivial knot $O$ and for the left trefoil $O$, and $v_3(123123)=1$.
(Here (123123) is the `non-trivial most symmetric chord diagram with 3 chords', see \cite[Figure 4.4]{PS96}, 5th diagram of the second line.)

(b,b',c,d) Same as problem \ref{vas-cal2} for $v_3$.

\emph{Hints.} See Problems 2, 3, 4ab, Results/Theorems 11, 13, 14 from \cite[\S4]{PS96}.
\end{pr}

\begin{pr}\label{vas-cal4} (a) \cite[Problem 4.4.b]{PS96} There exists a unique Vassiliev-Goussarov invariant $v_4:\Sigma\to\R$ of order at most 4 such that

$\bullet$ $v_4(O)=0$ for the trivial knot $O$, for the left trefoil $O$, and for the right trefoil $O$,

$\bullet$ $v_4(12341234)=2$, $v_4(12341432)=3$ and $v_4(12341423)=5$.

(b,b',c,d) Same as problem \ref{vas-cal2} for $v_4$.
\end{pr}

%\newpage
\section{Appendix: some details}\label{s:app}

\begin{figure}[h]\centering
\includegraphics[scale=0.25]{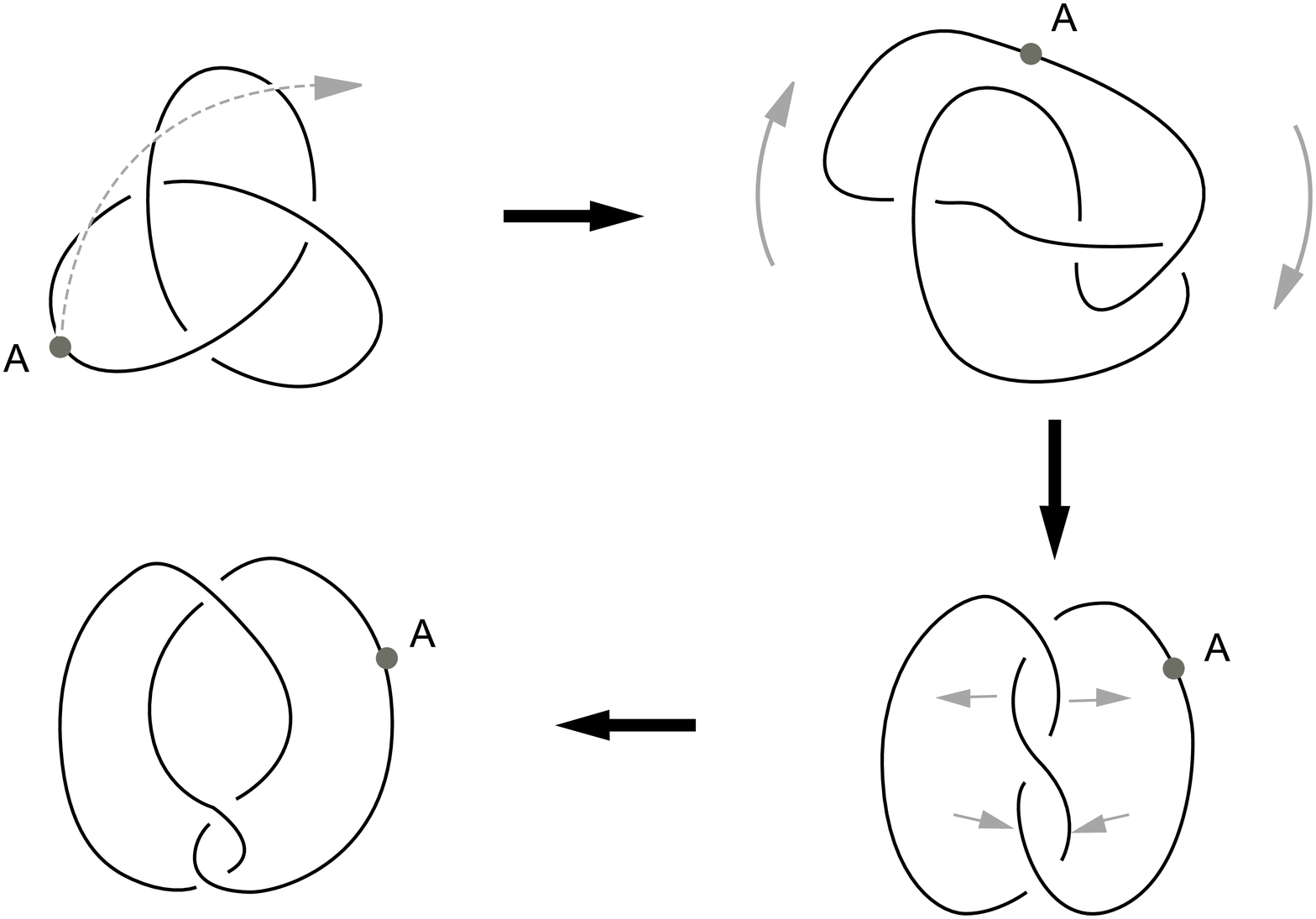}\qquad\includegraphics[scale=0.7]{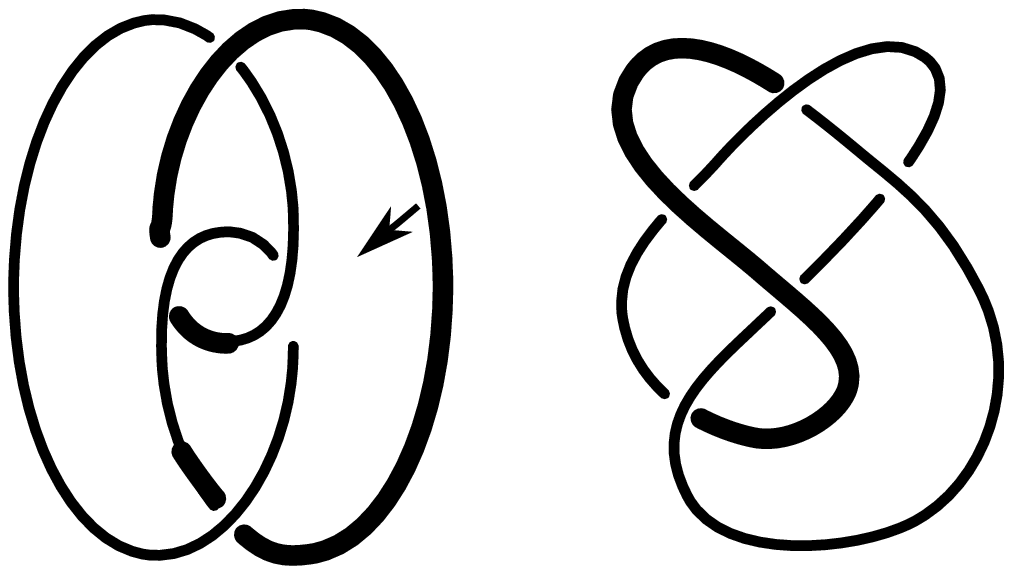}
\caption{Isotopy the trefoil and of the figure eight knot}
\label{f:itref}
\end{figure}

%\smallskip
{\bf \ref{a:crossing}.} (a,b) `Probably the best way of solving this problem is to make a model of the trefoil knot and the figure eight knot by using a shoelace and then move it around from one position to the other.
Fig. \ref{f:itref}
%\ref{f:eight}
gives some hints concerning transformations of the trefoil and the figure eight knot.' \cite[\S2]{Pr95}
(Fig. \ref{f:itref}, left, is prepared by D. Kroo.)

%\begin{figure}[h]\centering
%\caption{Isotopy of the figure eight knot}
%\label{f:eight}
%\end{figure}

%\cite{Pr95} fig. 2.13

(e) Consider two knots with coinciding plane diagrams in a `horizontal' plane $\pi$.
For each point $X$ in the space let $p(X)$ be the line containing $X$, perpendicular to $\pi$.
Let $h(X)$ be the height of $X$ relative to $\pi$, that is positive ($h(X)>0$) if $X$ is in the upper half-space,  and is negative ($h(X)<0$) if $X$ is in the lower half-space.
To each point $A$ of the first knot associate a point $A'$ of the second knot by the following procedure.
%There are two cases:

\emph{Case 1: The projection of the point $A$ on $\pi$ is not a crossing point on the plane diagram. }
In this case $p(A)$ intersects the first knot only at the point $A$.
Since the plane diagrams coincide, the line $p(A)$ intersects the second knot also at a single point.
Define $A'$ to be this point.

\emph{Case 2: The projection of the point $A$ on $\pi$ is a crossing point of the plane diagram.}
In this case the line $p(A)$ intersects the first knot in an additional point $B$.
Since the plane diagrams coincide, the line $p(A)$ intersects the second knot in two points $C$ and $D$,
where we assume that $h(C)>h(D)$. If $h(A)>h(B)$, we define $A'=C$, and in the opposite case $A'=D$.

For each point $A$ of the first knot and each number $t\in[0,1]$ let $A(t)$ be the point on the line $p(A)$ with the height $h(A(t))=(1-t)h(A)+th(A')$.
By construction $A(0)=A$, $A(1)=A'$ and the transformation of the first knot, which moves $A(0)$ in the direction of $A(1)$ with constant speed, so that at the time $t$ it occupies the position $A(t)$, is the required isotopy.

\begin{figure}[h]\centering
\includegraphics[scale=0.5]{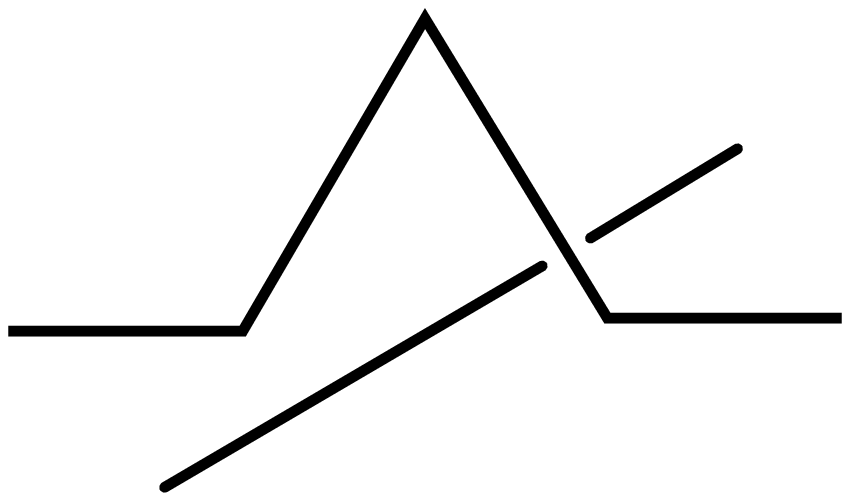}
\caption{The bridge over some crossing point}
\label{f:bridge}
\end{figure}

\smallskip
{\bf \ref{p:diag}.} See fig. \ref{f:bridge}.
For each crossing point of the plane diagram, on the upper edge of the
crossing, choose two points, close to the intersection and on the
opposite sides of the intersection.
Replace the line segment between the two chosen points by a `bridge' rising above the plane diagram, which
connects these two points.
After replacing all crossing points by the corresponding bridges, we obtain the required knot.

\smallskip
{\bf \ref{wha-tre}.} (a) Use the results of problems \ref{e:con-cas}.b, \ref{e:con-casca}.a and lemmas \ref{p:con-cas-ab}.ab, \ref{p:con-casca-ef}.ab.
Alternatively, use the results of problem \ref{wha-cole}.ab  and lemma \ref{wha-colre}.

(b) Take any of the connected sums of $n$ trefoil knots.
By the results of problems \ref{e:con-casca}.a and~\ref{con-casca}.a the Casson invariant of this knot is $n$.
Hence by lemma \ref{p:con-casca-ef}.ab these knots for different values of $n$ are not isotopic.

\smallskip
{\bf \ref{t:hopf}.} (a) In order to distinguish the Hopf link from the other two use the result of problem \ref{a:lk2} and lemma \ref{l:lk2}.
In order to distinguish the Whitehead link from the trivial link use the result of problem \ref{wha-cole} (or \ref{con-cal}) and lemma \ref{wha-colre} (or theorem \ref{con-ale}).

(b) Use the result of problem \ref{con-cal} and theorem \ref{con-ale}.

\smallskip
{\bf \ref{a:arfmot}.}
Choose a knot projected to the given plane diagram in the same way as in assertion~\ref{p:diag}.
Suppose that all the `bridges' lie in the upper half-space w.r.t. the projection plane.
By the assumption there are points $X$ and $Y$ on the knot which divide
the knot into two polygonal lines $p$ and $q$ such that

$\bullet$ $q$ lies in the projection plane and passes only through undercrossings;

$\bullet$ $p$ is projected to polygonal line $p'$ which passes only through overcrossings.

Take a point $Z$ in the upper half-space, and a point $T$ in the lower half-space.
Let us construct an isotopy between the given knot and the closed polygonal line $XZYT$,
which is isotopic to the trivial knot.
The isotopy consists of 3 steps, all of them keeping $X,Y$ fixed.

{\it Step 1. An isotopy between $q$ and $XTY$.}
Suppose that $q=A_0A_1\ldots A_n$, where $A_0=X$ and $A_n=Y$.
Then the isotopy is given by
$$
A_0A_1\to A_0TA_1,\quad TA_1A_2\to TA_2,\quad TA_2A_3\to TA_3,\quad\ldots TA_{n-1}A_n\to TA_n.
$$

{\it Step 2. An isotopy between $p$ and $p'$.} Remove all the `bridges' by elementary moves.

{\it Step 3. An isotopy between $p'$ and $XZY$.} This is done analogously to step~1.

\smallskip
{\bf \ref{l:crossing}.} Follows by assertion \ref{a:arfmot}.

{\it Another idea of the proof (cf.~\cite[Theorem 3.8]{PS96}).}
Denote by $\pi$ the horizontal plane containing the plane diagram.
For each point $X$ in the space, $p(X)$ and $h(X)$ are defined in the solution of the problem~\ref{a:crossing}.c.
Let $l$ be a line in the plane, which passes through a vertex $A_0$ of the plane diagram, while the whole diagram is contained in one of the two half-planes determined by $l$.
Let $A_0,A_1,\ldots A_n$ be all vertices of the plane diagram, in the order of their appearance, while we move along the diagram in some direction.
Choose points $B_0,\ldots,B_n$ so that $A_i\in p(B_i)$ for $i=1,\ldots,n$, and $h(B_i)<h(B_j)$ for $i<j$.
Let $B_{n+1}$ be a point, whose projection on $\pi$ is close to $A_0$ and $h(B_{n+1})>h(B_n)$.
We claim that the knot $B_0\ldots B_nB_{n+1}$ is isotopic to the trivial knot.
Indeed, by the choice of the line $l$, the projection of the knot onto any plane, perpendicular to
the line $l$, is a closed polygonal line without self-intersections.
It remains to modify crossing in the plane diagram so that they are in
agreement with the projection of the constructed knot to the plane $\pi$.

\smallskip
{\bf \ref{a:lk2}.} \textit{Answer:} 1 for the Hopf link and 0 for other links.

\smallskip
{\bf \ref{l:lk2}.} For moves I and III the number of crossing points where the first component passes above the second one does not change.
For move II this number changes by \(0\) or \(\pm2\).
	
\smallskip
{\bf \ref{a:lk2-pr}.} (a) Take a plane diagram of a link. 	
%At any crossing point either the first component passes above the second one,
%or the second component passes above the first one.
By the Parity lemma stated before assertion \ref{a:lk2-pr} the number of crossing points where the first component passes above the second one has the same parity as the number of crossing points where the second component passes above the first one.
This is the required statement.

(b) An example is the third link in fig. \ref{f:hopf}.
This link is not isotopic to the trivial link because they have distinct linking numbers, see \S\ref{0rapl}.

\smallskip
{\bf \ref{con-lk2}.} {\it Existence.} By lemma \ref{l:lk2} the linking number modulo 2 is an isotopy invariant.
The skein relation is easy to check.

{\it Uniqueness.} Suppose that \(\mathrm f\) is another invariant aside from \(\lk_2\) satisfying the assumptions.
Then \(\mathrm f-\lk_2\) is an isotopy invariant assuming zero value on the trivial link and invariant under crossing changes.
The analogue of lemma \ref{l:crossing} for links states that any plane diagram of a link can be obtained from the diagram of a link isotopic to the trivial link by some crossing changes.
Hence $\mathrm f-\lk_2=0$.

\smallskip
{\bf \ref{e:con-cas}.} (a) If $P$ is a point on the plane diagram as in assertion~\ref{a:arfmot}, then there are no $P$-skew pairs of crossings.
Hence the $P$-Arf invariant is zero.

\smallskip
{\bf \ref{p:con-cas-ab}.}
(a) Let $P_1$ and $P_2$ be two basepoints such that the segment $P_1P_2$ contains exactly one crossing point $X$.

{\it Case 1: $P_1P_2$ passes through undercrossing}.
Then $X$ does not form either $P_1$-skew or $P_2$-skew pair with any other crossing.
Hence $P_1$- and $P_2$-Arf invariants of the diagram are equal.

{\it Case 2: $P_1P_2$ passes through overcrossing}.
Then $X$ divides the diagram into two closed polygonal lines $q_1$ and $q_2$ such that
$P_1$ lies on $q_1$ and $P_2$ lies on $q_2$.
Denote by $n_1$ (respectively, $n_2$) the number of intersections of $q_1$ and $q_2$ for which
$q_1$ passes above $q_2$ (respectively, $q_2$ passes above $q_1$).
Denote by $N_1$ the number of $P_1$-skew pairs formed by $X$ and some intersection of $q_1$ and $q_2$.
Denote by $\arf_{P_1}D$ the $P_1$-arf invariant of $D$.
Use analogous notation with $P_1$ replaced by $P_2$.
Then
$$
{\arf}_{P_1} D-{\arf}_{P_2} D=N_1-N_2=n_1-n_2\underset{2}\equiv n_1+n_2\underset{2}\equiv0,
$$
where $D$ is the given plane diagram. Here

$\bullet$ the first equality holds because a pair of crossings is either $P_1$-skew or $P_2$-skew (but not both) if and only if the pair is formed
by $X$ and some intersection of $q_1$ and $q_2$;

$\bullet$ the second equality holds because $N_1=n_1$ and $N_2=n_2$; indeed, an intersection of $q_1$ and $q_2$ forms a $P_1$-skew (respectively, $P_2$-skew) pair with $X$ if and only if at this intersection $q_1$ passes above (respectively, below) $q_2$;

$\bullet$ $\underset{2}\equiv$ are congruences modulo 2;

$\bullet$ the last congruence follows by the Parity lemma for $q_1$ and $q_2$.

\begin{figure}[h]\centering
\includegraphics[scale=0.45]{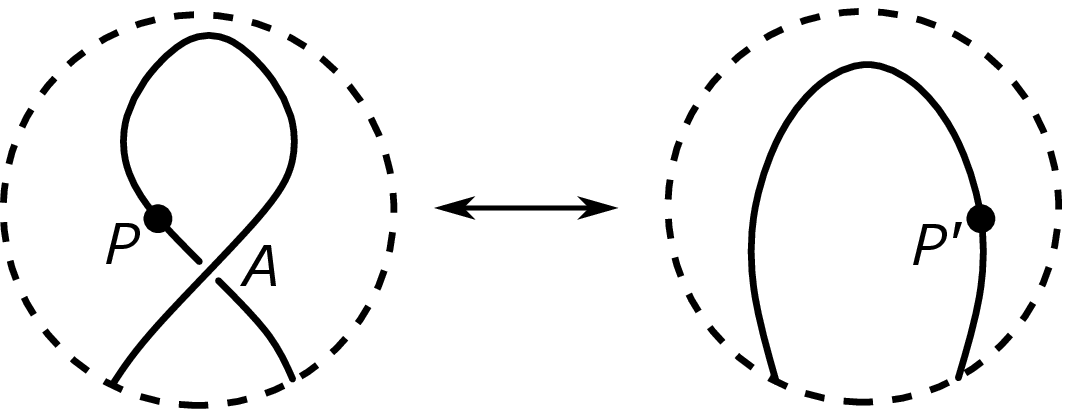}\qquad
\includegraphics[scale=0.45]{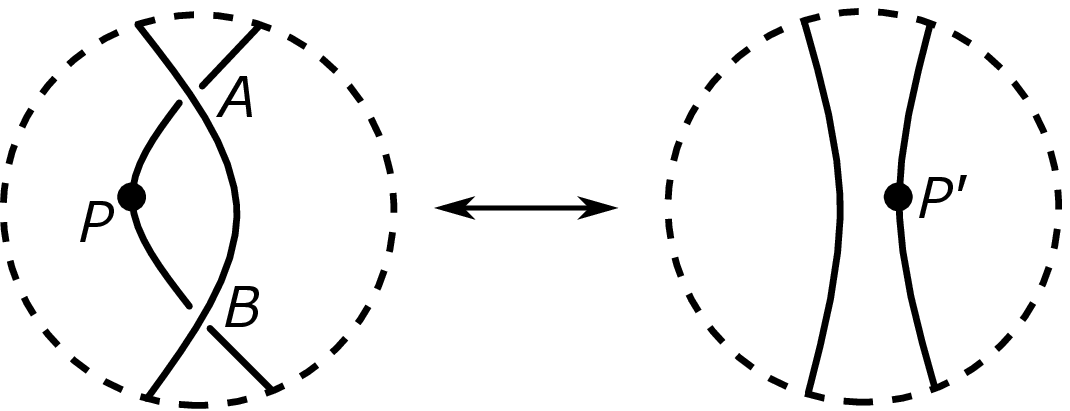}\qquad
\includegraphics[scale=0.45]{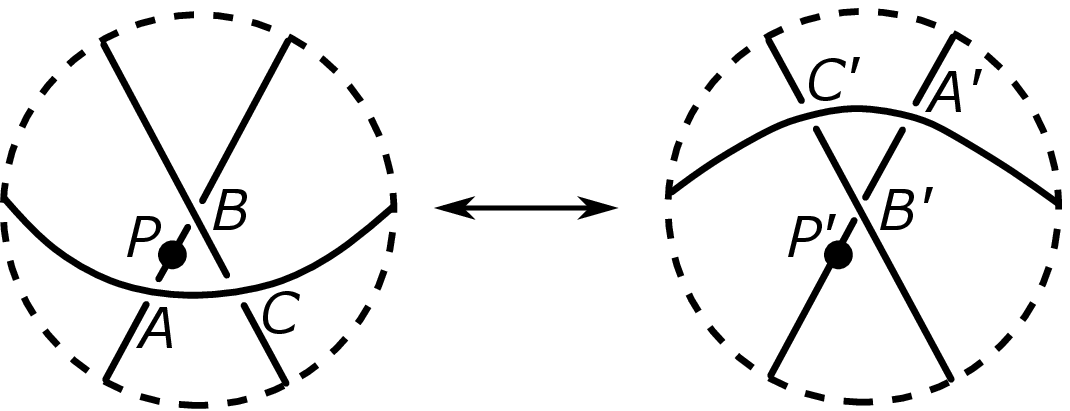}
\caption{Arf-invariant does not change under Reidemeister moves}
%To the proof of \ref{con-cas}.f
\label{f:arf-reid}
\end{figure}

(b) {\it Type $I$ move.} Take basepoints before and after the move as in fig.~\ref{f:arf-reid} (left).
Check that the crossing $A$ does not form a $P$-skew pair with any other crossing.

{\it Type $II$ move.} Take basepoints before and after the move as in fig.~\ref{f:arf-reid} (middle).
Check that neither of the crossings $A$ and $B$ forms a $P$-skew pair with any other crossing.

{\it Type $III$ move.} Take basepoints before and after the move as in fig.~\ref{f:arf-reid} (right).
Check that neither of the crossings $A$, $B$ forms a $P$-skew pair with any other crossing and that neither of the crossings $A'$, $B'$ forms a $P'$-skew pair with any other crossing.
Then check that a crossing $X$ distinct from $A$, $B$, $C$ forms a $P$-skew pair with $C$ if and only if $X$ forms a $P'$-skew pair with $C'$.

\smallskip
{\bf \ref{a:arf-pr}.} Take any connected sum $K$ of the two trefoil knots.
By assertion~\ref{wha-cons}.d $\arf K=0$.
By the result of problem~\ref{e:con-casca}.b and assertion \ref{con-casca}.a $c_2(K)\ne0$.
%is 2.
Hence $K$ is not isotopic to the trivial knot.

\begin{figure}[h]\centering
\includegraphics[scale=0.6]{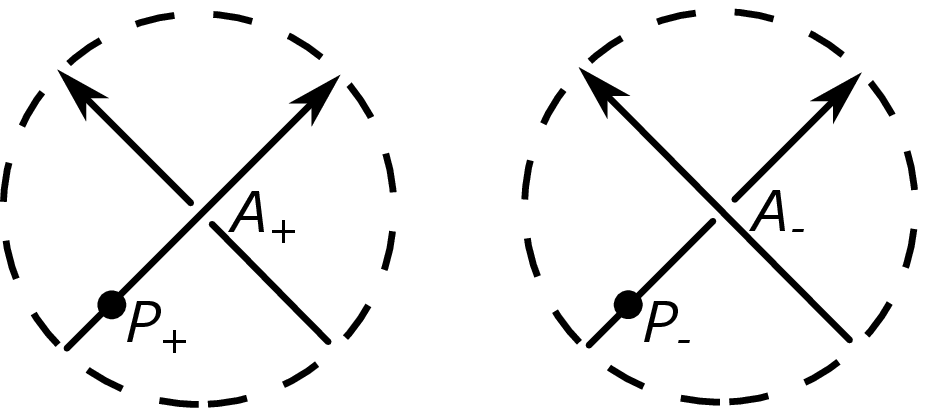}
\caption{To the proof of skein relation for Arf invariant}
%\ref{p:arf-pr}
\label{f:arf-skein}
\end{figure}

\smallskip
{\bf \ref{con-thma}.} {\it Existence.} By lemma \ref{p:con-cas-ab}, the arf invariant is an isotopy invariant.
Here are hints for checking the skein relation.
Take basepoints $P_+$, $P_-$ as in fig.~\ref{f:arf-skein}.
Check that the crossing $A_-$ does not form a $P_-$-skew pair with any other crossing in $K_-$.
Then check that the number of crossings which form a $P_+$-skew pair with $A_+$ in $K_+$ equals $\lk\phantom{}_2K_0$ modulo 2.

{\it Uniqueness.} The proof is analogous to the proof of theorem~\ref{con-lk2}.
Use lemma~\ref{l:crossing} itself instead of its analogue for links.

\begin{figure}[h]\centering
\includegraphics[scale=0.7]{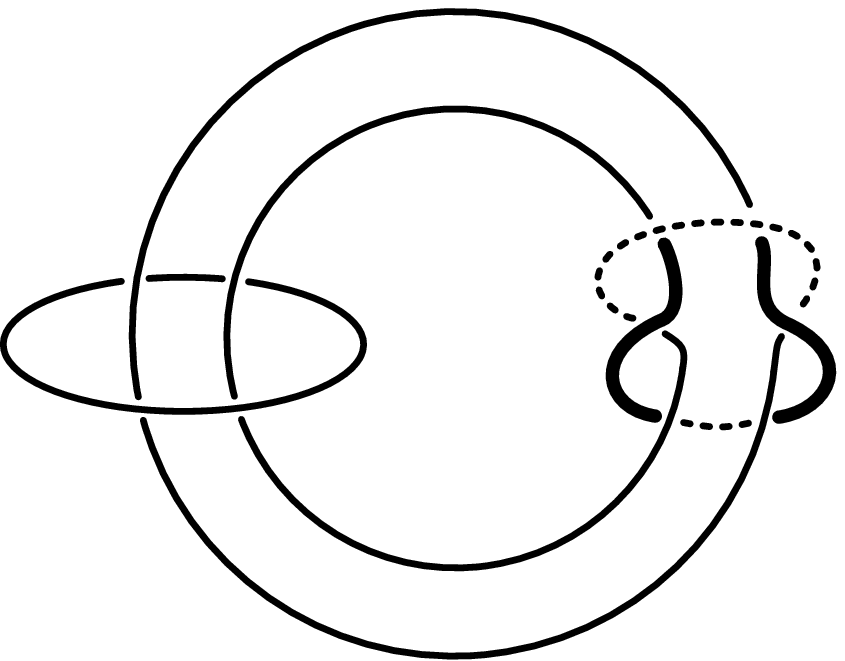}\qquad\includegraphics[scale=0.9]{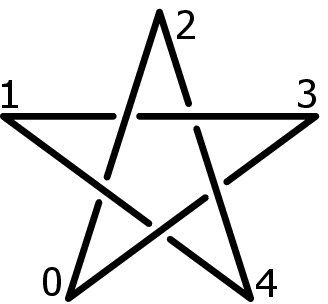}
\caption{A 3-coloring of a link and a 5-coloring of the $5_1$ knot}
\label{f:51-colored}
\end{figure}

\smallskip
{\bf \ref{wha-cole}.} {\it Answers:} b,e,h --- 3-colorable, a,c,d,f,g,i --- not 3-colorable.
For a proper coloring of a diagram of trefoil knot see
%the attached copy of~
\cite[p. 30, figure 4.3]{Pr95}.
For a proper coloring of the last diagram from fig. \ref{f:hopf} see  fig. \ref{f:51-colored} left.
(This diagram was erroneously stated to be not 3-colorable in \cite[\S4]{Pr95}. This minor mistake was found by L.M. Bann\"ohr, S. Zotova and L. Kravtsova.)

\smallskip
{\bf \ref{wha-colr}.}
(a) Most part of (a) follows by lemma \ref{wha-colre} and results of problems~\ref{wha-cole}.d-h
(see \cite[p. 30]{Pr95}).
The last diagram from fig. \ref{f:hopf} is distinguished from the trivial link by the number of proper colorings of a plane diagram.
Prove that this number is preserved under the Reidemeister moves.

(b) A plane diagram is {\it 5-colorable} if there exists a coloring of its strands in five colors $0,1,2,3,4$ so that at least two colors are used, and at each crossing if the upper strand has color $a$ and two lower strands have colors $b$ and $c$, then $2a\equiv b+c\pmod{5}$.
Similarly to lemma~\ref{wha-colre} the 5-colorability of a plane diagram is preserved under Redemeister moves.
The $5_1$ knot is 5-colorable, see fig.~\ref{f:51-colored}, right.
The trivial knot is not.
Hence they are not isotopic.

\smallskip
{\bf \ref{wha-chirsph}.} (b) {\it First solution.}
An oriented polygonal line is called {\it positive} if the bounded part
of the plane is always on the left side of each of its oriented segments
(see the Jordan theorem in remark \ref{r:accur}).
Prove that the positivity of an oriented polygonal line is preserved by elementary moves.

{\it Hint to the second solution.} The positivity can be equivalently defined as follows.
We say that an oriented polygonal line $A_1\ldots A_n$ is {\it positive} if for each of its inner (interior) points $O$ the sum of oriented angles
$\angle A_1OA_2+\angle A_2OA_3+\ldots+\angle A_{n-1}OA_n+\angle A_nOA_1$
is always positive (i.e. the {\it winding number} of the oriented polygonal line around any interior point is positive).

\smallskip
{\bf \ref{wha-inv}.} Each of the three indicated oriented knots is transformed into the oriented  knot with the opposite orientation by the rotation through the angle $\pi$ around the `vertical' axis passing through the `upper' point of the knot (see the leftmost diagram in fig. \ref{f:trefoil}, the first and the second row for the trefoil and the figure eight knot, respectively).
This rotation is included into a continuous family of rotations through the angle $\pi t$, $t\in[0,1]$, with respect to the same line. This is the required isotopy.

\smallskip
{\bf \ref{wha-consli}.} (d) Check that all crossings of different components in $K\#L$ are
obtained from such crossings in $K$ and in $L$.

% For linking number via plane diagrams.

\smallskip
{\bf \ref{e:lk}.} 	\textit{Answers:} $\pm1$; 0.

\smallskip
{\bf \ref{p:lk-pr}.} (a) The proof is analogous  to assertion \ref{a:lk2-pr}.a.
Take a plane diagram of a link. 	
By the Triviality lemma (stated before assertion \ref{p:lk-pr}) the sum of signs of crossing points where the first component passes above the second one has opposite sign to the sum of signs of crossing points where the second component passes above the first one.
Switching the components negates the sign of every crossing point.
This completes the proof.

(b) Reversing the orientation of either of the components negates the sign of every crossing point.

(c) Take the connected sum of 5 Hopf links oriented so that their linking numbers equal to $-1$.

(d) The proof is analogous to assertion \ref{wha-cons}.d.
The signed set of crossing points of plane diagram of \(K\#L\) is the union of the signed sets of crossing points of plane diagrams of links \(K\) and \(L\).

(e) An example is the Whitehead link.
The Whitehead link is not isotopic to the trivial link by theorem \ref{t:hopf}.a.

\smallskip
{\bf \ref{e:con-casca}.} (a) \textit{Answers:} 0, 1 and $-1$.

The trivial knot has no crossings, and so no skew pairs of crossings.
Therefore the Casson invariant of this knot is 0.

All three crossings of the trefoil knot have the same sign.
Since the trefoil knot has exactly one linked pair of crossings (regardless the choice of the base-point),
we obtain that the Casson invariant of this knot is 1.

(b) Take any connected sum of five figure eight knots.
By (a) and assertion~\ref{con-casca}.a below the Casson invariant of this knot is $-5$.

%(Actually, such knots with the smallest number of self-intersections have 10 self-intersections.)

\smallskip
{\bf \ref{p:con-casca-ef}, \ref{con-casca}.a, \ref{con-thm}.} The proofs are analogous to lemma~\ref{p:con-cas-ab}, assertion~\ref{a:arf-pr} and theorem~\ref{con-thma}, respectively.
Take care of the signs of intersection points.
For lemma \ref{p:con-casca-ef}.a use the Triviality lemma stated before assertion~\ref{p:lk-pr}.

% (it suffices to take base-points in the same way).

\smallskip
{\bf \ref{con-casca}.} (b) Take any connected sum of the trefoil knot and the figure eight knot.
By (a) and the answer to problem \ref{e:con-casca}.a the Casson invariant of this knot is 0.
However, by the answers to problems~\ref{con-cal}.d,e and assertion~\ref{con-pro}.c the Conway polynomial of this knot is $(1+t^2)(1-t^2)\ne1$.
Hence this knot is not isotopic to the trivial knot.

%(d) Fix any base-point and check all pairs of crossings.

%\cite[\S4]{Pr95}

\smallskip
{\bf \ref{t:a3}.} (b) {\it Answer:} 0.

{\it Remark.} The invariant $\arf=c_2\mod2$
%$,a_3$
for links may depend on the orientation on the components (for $c_3\mod2$ see \cite[2.3.4]{CDM}).

%However the answers to problems (a), (b), (c) do not depend on the choice of orientation.

\smallskip
Let $D$ be a plane diagram of a link.
Denote by $\mathrm{cr}D$ the number of crossings in $D$.
Denote by $\mathrm{ch}D$ the minimal number of crossing changes required to obtain from $D$ a diagram of a link which is isotopic to the trivial one (such sequence of crossing changes exists by the analogue of lemma \ref{l:crossing} for links).

\smallskip
{\bf \ref{con-ale}.} The uniqueness is analogous to theorems \ref{con-lk},\ref{con-thm}; solve first problem \ref{con-cal}.
%See important comments in the internet version of this project.

{\it Deduction of the stronger version (for non-ordered links) from the weaker version (for ordered links).}
It suffices to show that all invariants $c_n$ defined for ordered links are preserved under changes of the order of the components.

Let $D$ be a plane diagram of a link with two or more components and let $D'$ be a plane diagram obtained from $D$ by a change of the components' order.
The proof is by induction on $\mathrm{cr}D$.
If $\mathrm{cr}D=0$, then $D$ is a diagram of a link which is isotopic to the trivial one and by the answer to
problem~\ref{con-cal}.b we have $C(D)=0$ for any ordering of the components.
Suppose that $\mathrm{cr}\,D>0$; then continue the proof by induction on $\mathrm{ch}D$.
If $\mathrm{ch}D=0$, then $D$ is a diagram of a link which is isotopic to the trivial one; this case is considered above.
Suppose that $\mathrm{ch}D>0$.
Let $D_*$ be a link obtained from $D$ by a crossing change and such that $\mathrm{ch}D_*<\mathrm{ch}D$.
Denote by $D'_*$ is a link obtained from $D'$ by the change of the same crossing.
Then
$$
\pm(C(D)-C(D_*))=C(D_0)\quad\mbox{and}\quad\pm(C(D')-C(D'_*))=C(D'_0),
$$
where $D_0$ is a diagram of a link $K_0$ (with some ordering of the components) from fig.~\ref{f:skein} for $D$, $D_*$ being $D_+$, $D_-$ in some order, and $D'_0$ is the same for   $D'$, $D'_*$.
Note that the diagrams $D_*$ and $D'_*$ coincide up to the order of the components.
The same is true for the diagrams $D_0$ and $D'_0$.
Since $\mathrm{ch}D_*<\mathrm{ch}D$ and $\mathrm{cr}\,D_0<\mathrm{cr}\,D$, by the inductive hypotheses we have $C(D_*)=C(D'_*)$ and $C(D_0)=C(D'_0)$.
Then $C(D)=C(D')$.

%Use assertion \ref{l:crossing}.b
% is given in \cite[2.3.1]{CDM} (the `missing' property is that the Conway polynomials of plane diagrams
%differing by Reidemeister moves are the same). no existence

\smallskip
{\bf \ref{con-cal}.} {\it Answers:} (a, b) 0; (c) $\pm t$; (d) $1+t^2$; (e) $1-t^2$; (f) $\pm t^3$; (g) $\pm t^4$; (h) $1+3t^2+t^4$.

{\it Remark.} The signs in the answers to (c), (f), (g) depend on the orientation on the components.

%Use the Conway skein relation.

{\it Hint.} For examples of such calculations for (a), (c), and (d) see
%the attached copy of~
\cite[2.3.2]{CDM}.

\smallskip
{\bf \ref{p:conzer}.} Let $D$ be a plane diagram of the given link $K$.

(a) For any diagram $D_*$ obtained from $D$ by a crossing change we have $c_0(D)-c_0(D_*)=0$, i.~e. $c_0$ does not change under crossing changes.
By the analogue of lemma~\ref{l:crossing} for links the diagram $D$ can be obtained by crossing changes from a diagram of a link isotopic to the trivial one.
The assertion follows from the definition of $c_0$ on the trivial knot and the answer to problem~\ref{con-cal}.b.

(b,c) The first parts are particular cases of (d).
The second parts follow from the definition of $c_1,c_2$ and theorems~\ref{con-lk}, \ref{con-thm}.

(d) The proof is by induction on $\mathrm{cr}\,D$.
If $\mathrm{cr}\,D=0$, then $K$ is isotopic to the trivial link.
If $K$ is a knot, then $C(D)=1$.
Otherwise $C(D)=0$ by assertion~\ref{con-cal}.b.
Suppose that $\mathrm{cr}\,D>0$; then continue the proof by induction on $\mathrm{ch}D$.
If $\mathrm{ch}D=0$, then $K$ is isotopic to the trivial link; this case is considered above.
Suppose that $\mathrm{ch}D>0$.
Let $D_*$ be a link obtained from $D$ by a crossing change and such that $\mathrm{ch}D_*<\mathrm{ch}D$.
Then $\pm(c_j(D)-c_j(D_*))=c_{j-1}(D_0)$, where $D_0$ is the diagram from fig.~\ref{f:skein} corresponding to $D$, $D_*$ being $D_+$, $D_-$ in some order.
Note that the link $D_*$ consists of $k$ components and the link $D_0$ consists of $k'=k\pm1$ components. Therefore if $j\le k-2$, then $j-1\le k'$ and if $j-k$ is even, then $(j-1)-k'$ is even.
Since $\mathrm{ch}D_*<\mathrm{ch}D$ and $\mathrm{cr}\,D_0<\mathrm{cr}\,D$, by the inductive hypothesis we have $c_j(D_*)=c_{j-1}(D_0)=0$.
Then $c_j(D)=0$.

(e) Prove analogously to (d) that {\it $c_j(D)=0$ for any plane diagram $D$ and $j>\mathrm{cr}\,D$}.

\smallskip
{\bf \ref{con-pro}.} (a) The proof is analogous to assertion~\ref{p:conzer}.d.

(c) Let $D$ and $E$ be plane diagrams of $K$ and $L$.
Analogously to assertion~\ref{p:conzer}.d prove that $C(D\#E)=C(D)C(E)$ by induction on $\mathrm{cr}\,D$ for fixed $E$.

\begin{figure}[h]\centering
\includegraphics[scale=0.6]{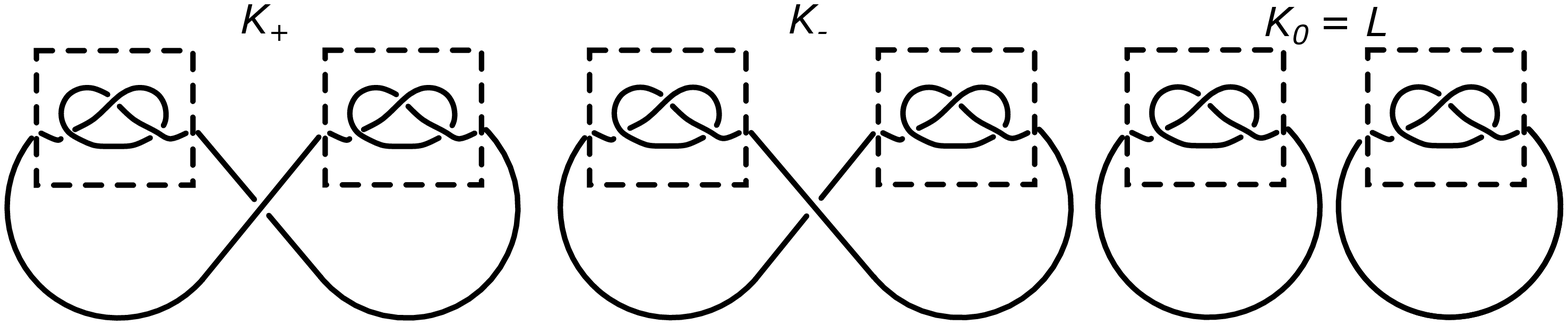}
\caption{Proof that $C(\mbox{split link})=0$}
\label{f:con-split}
\end{figure}

\smallskip
{\bf \ref{con-spl}.} (a) Follows from answers to problems \ref{con-cal}.c,f,g and (b,c).

(c) If $L$ is a split link, then there exist links $K_+$, $K_-$, $K_0$ such that

$\bullet$ their plane diagrams differ like in fig.~\ref{f:skein};

$\bullet$ the links $K_+$ and $K_-$ are isotopic;

$\bullet$ the link $K_0$ is isotopic to $L$.

We have $C(L)=C(K_0)=\frac{1}{t}(C(K_+)-C(K_-))=0$, see fig. \ref{f:con-split}.

%\newpage

%Colin C. Adams, The Knot Book: An elementary introduction to the mathematical theory of knots., 2004.

\end{document}

\title{Основы теории узлов и зацеплений для пользователя\footnote{Благодарю А. Еннэ за подготовку рисунков и работу над \cite{EEF} после Летней конференции,
Е. Морозова за перевод обновленной мной версии на русский язык и А. Глебова, Д. Елисеева, А. Еннэ, Р. Живалевича,  Е. Морозова, А. Рябичева, А. Сосинского, М. Федорова, Н. Хорошавкину, С. Чмутова за полезные обсуждения и нашу работу над \cite{EEF}.
Этот текст основан на лекциях автора в Независимом Московском Университете (включая Math in Moscow Program) и
в Московском Физико-Техническом Институте, а также на \cite{EEF}.}}

\author{А. Скопенков\footnote{\texttt{https://users.mccme.ru/skopenko},
Московский Физико-Технический Институт, Независимый Московский Университет.
Частично поддержан грантом РФФИ No. 19-01-00169 и грантом фонда Саймонса-НМУ.}}

\date{}

\maketitle

\begin{abstract}
Это методическая статья.
Мы определяем простые инварианты узлов и зацеплений (коэффициент зацепления, инварианты Арфа-Кэссона и многочлен Александера-Конвея) для применения к интересным результатам, формулировки которых доступны неспециалисту.
Простейшие инварианты естетственно появляются при попытках распутать узел или зацепление.
Затем мы представляем рекурсивные скейн-соотношения для простейших инвариантов, которые позволяют ввести более сложные инварианты.
Мы формулируем теорему Васильева-Концевича в виде, удобном для вычислления самих инвариантов, а не только размерности пространства инвариантов.
Предварительных знаний не требуется; мы даем строгие определения основных объектов на языке, не препятствуующем их интуитивному пониманию.
\end{abstract}

\tableofcontents

\newpage
{\bf О стиле этого текста}

Как правило, я привожу \emph{формулировку} красивого или важного утверждения \emph{перед} последовательностью определений и результатов, составляющих его \emph{доказательство}.
В~таких случаях для доказательства утверждения требуется часть дальнейшего материала.
Об этом указано после формулировок (но я не хочу лишать Вас удовольствия найти момент, когда вы наконец-то сможете доказать такое утверждение).
Некоторые теоремы приведены без доказательств, тогда даются ссылки вместо указаний.

{\it Утверждения} в этом тексте --- простые элементы теории
(для читателя, частично знакомого с материалом, они являются краткими напоминаниями).
С той же целью приведено небольшое количество задач.
Указания, решения и доказательства к утверждениям, задачам, теоремам и леммам приведены в \S\ref{s:app}.
Однако рекомендую читателю подумать над ними самостоятельно.
%In order to get a thorough understanding of the material a reader can also
%consider theorems and lemmas as problems (beware the previous paragraph!).
%[Следующие фраза и сноска взяты из опубликованной в МП статьи:]
Изучение путем решения задач не только характерно для серьезного изучения математики \cite[\S1.1]{HC19}, \cite[\S1.2]{ZSS}, но и продолжает древнюю культурную традицию.\footnote{Например, послушники дзенских монастырей обучаются, размышляя над загадками, данными им наставниками.
Впрочем, эти загадки являются скорее наводящими на размышления парадоксами, а не задачами.
См. подробнее \cite{S}.}

%Some more complicated statements are marked by stars.
{\it Замечания} формально не используются в дальнейшем.

\section{Основные определения и результаты об узлах}\label{0isot}

\begin{figure}[h]\centering
\includegraphics[scale=3.6]{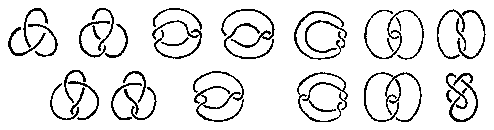}
\caption{Узлы, изотопные трилистнику (верхний ряд) и восьмерке (нижний ряд)}
\label{f:trefoil}
\end{figure}

Начнем с неформального описания основных понятий (строгие определения даны после замечания  \ref{r:why}).
{\it Узел} можно представить себе как тонкую эластичную веревку, концы которой склеены ---
см. рис. \ref{f:trefoil}.
Как и на этом рисунке, узлы обычно изображаются с помощью их <<удобной>> проекции на плоскость, которая называется {\it диаграммой узла}.
Представьте, что вы раскладываете веревку на столе, тщательно записывая как она пересекает себя (т.е. какая часть веревки лежит сверху).
Необходимо помнить, что проекции одного и того же узла на разные плоскости могут выглядеть по-разному.

{\bf Тривиальный узел} --- это контур (граница) треугольника.

Под {\it изотопией} узла понимается его непрерывная деформация в пространстве как тонкой эластичной нити;
любые самопересечения нити во время деформации запрещены.
Два узла называются {\it изотопными}, если один из них можно преобразовать в другой с помощью изотопии.

\begin{assertion}\label{a:crossing}
%(a) Некоторые два узла из верхнего ряда рис. \ref{f:trefoil} изотопны самому левому узлу.
(a) Все узлы, представленные в верхнем ряду на рис. \ref{f:trefoil}, изотопны друг другу.
(Для одной пары из этих узлов представьте изотопию с помощью движений Рейдемейстера, показанных на рис. \ref{f:reid}.)

(b) То же верно для узлов в нижнем ряду на рис. \ref{f:trefoil}.

(c) Все узлы с одинаковыми диаграммами изотопны.
\end{assertion}

\begin{remark}[почему строгое определение изотопии необходимо?]\label{r:why}
На рис. \ref{f:naisotopy} показана изотопия между трилистником и тривиальным узлом.

\begin{figure}[h]\centering
\includegraphics[scale=0.6]{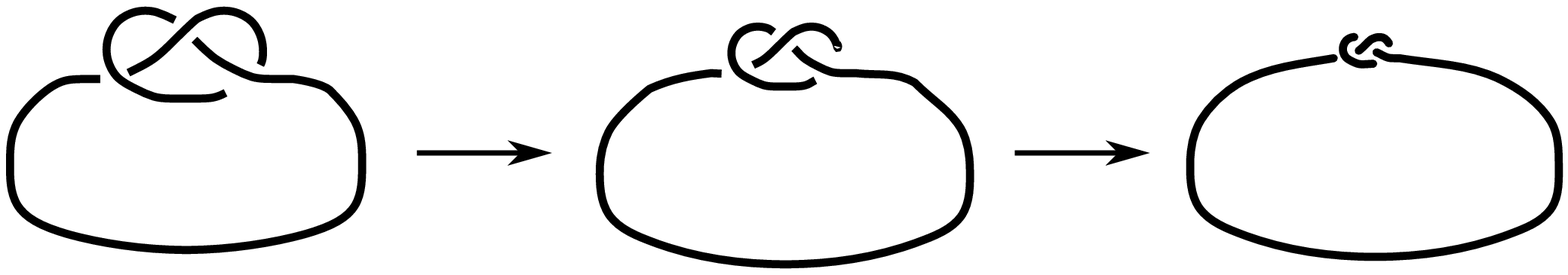}
\caption{Необъемлемая изотопия между трилистником и тривиальным узлом}
\label{f:naisotopy}
\end{figure}

Действительно ли это изотопия?
Это так называемая <<кусочно-линейная необъемлемая изотопия>>, которая \emph{отличается} от <<кусочно-линейной объемлемой изотопии>>, определенной и используемой далее.
(Первое понятие лучше отражает идею непрерывной деформации без самопересечений, но малодоступно для старшеклассников, ср. \cite{Sk16i}.)
На самом деле, любые два узла кусочно-линейно необъемлемо изотопны!

Обычная проблема при работе с интуитивными определениями не в том, что их трудно сделать строгими,
а в том, что это можно сделать разными неэквивалентными способами.
\end{remark}

{\bf Узел} --- это пространственная замкнутая несамопересекающаяся ломаная.\footnote{Его не следует путать с {\it ориентированным узлом}, определенным далее в \S\ref{s:ori}.}

{\bf Плоская диаграмма} узла --- это его проекция общего положения\footnote{Ломаная на плоскости называется ломаной {\it общего положения}, если существует ломаная $L$ с тем же объединением ребер такая, что никакие 3 вершины ломаной $L$ не лежат на одной прямой, и никакие три отрезка, соединяющие некоторые вершины ломаной $L$, не имеют общей внутренней точки.}
на плоскость\footnote{Терминология университетской математики --- <<образ общего положения при проекции на плоскость>>.}, вместе с информацией, какая часть узла <<проходит под>> (проход), а какая --- <<проходит над>> (переход) в каждом конкретном перекрестке.

\begin{assertion}\label{p:diag} Для любой диаграммы узла существует узел, проецирующийся на эту диаграмму.
(Такой узел не обязательно единственный; см. впрочем утверждение \ref{a:crossing}.c.)
\end{assertion}

\begin{figure}[h]\centering
\includegraphics[scale=0.5]{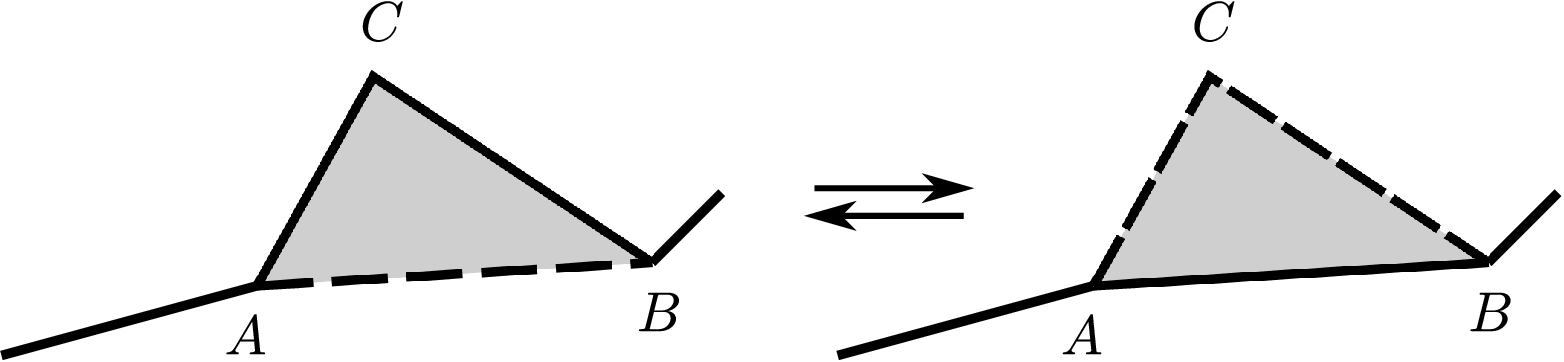}
\caption{Элементарное движение}
\label{f:elem}
\end{figure}

Предположим, что стороны $AC$ и $CB$ треугольника $ABC$ --- ребра некоторого узла.
Кроме того, предположим, что узел и часть плоскости, ограниченная треугольником $ABC$, ни в каких других точках не пересекаются.
{\bf Элементарное движение} $ACB\to AB$ --- это замена двух ребер $AC$ и $CB$ ребром $AB$, либо обратная операция $AB\to ACB$ (рис. \ref{f:elem}).\footnote{Если треугольник $ABC$ вырожденный, то элементарное движение является либо разбиением ребра, либо обратной операцией.}
Два узла $K,L$ называются (кусочно-линейно объемлемо) {\bf изотопными}, если существует последовательность узлов $K_1,\ldots,K_n$ такая, что $K_1=K$, $K_n=L$, и каждый узел последовательности $K_{j+1}$ получен из предыдущего узла $K_j$ элементарным движением.

\begin{theorem}\label{wha-tre} (a) Следующие узлы попарно не изотопны: тривиальный узел, трилистник, восьмерка.

(b) Существует бесконечное множество попарно не изотопных узлов.
\end{theorem}

Эта теорема доказывается с помощью {\it инвариантов Арфа} и {\it Кэссона},
см. \S\ref{0arf} и \S\ref{0cas}, ср. \S\ref{s:reco}.

{\it Зеркальный образ} узла $K$ --- это узел, плоская диаграмма которого получается изменением всех перекрестков (рис. \ref{f:crossing}) диаграммы узла $K$.
Из утверждения \ref{a:crossing}.b следует, что восьмерка изотопна своему зеркальному образу.

\begin{theorem}\label{wha-nchir} Трилистник не изотопен своему зеркальному образу.
\end{theorem}

Теорема \ref{wha-nchir} доказывается с помощью {\it полинома Джонса} \cite[\S3]{PS96}, \cite[\S2.4]{CDM}.
Доказательство выходит за рамки данного текста.

\begin{theorem}\label{arf-k7}
Если никакие 4 из 7 точек пространства не лежат в одной плоскости, то существует замкнутая ломаная, составленная из некоторых
%${7\choose2}=21$
отрезков, соединяющих эти точки, не изотопная контуру треугольника.
\end{theorem}

Эта теорема доказывается с помощью {\it Арф-инварианта}, см. \S\ref{0arf}.
Детали выходят за рамки данного текста, см. \cite[теорема~2]{CG83}.
%Конвей--Гордон

\section{Основные определения и результаты о зацеплениях}\label{0isoli}

{\bf Зацепление} --- это набор попарно непересекающихся узлов, которые называются {\it компонентами} зацепления.
Упорядоченные наборы называются упорядоченными или раскрашенными зацеплениями, в то время как неупорядоченные наборы называются неупорядоченными или нераскрашенными зацеплениями.
В этом тексте под словом <<зацепление>> понимается <<упорядоченное зацепление>>.

{\bf Тривиальное зацепление} (с любым числом компонент) --- это зацепление, состоящее из треугольников, лежащих в параллельных плоскостях.

{\it Плоские диаграммы} и {\it изотопии} для зацеплений определяются так же, как и для узлов.

\begin{figure}[h]\centering
\includegraphics{3-1.eps}\includegraphics[scale=0.6]{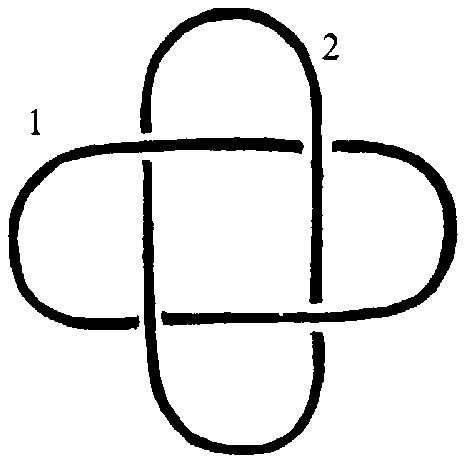}
\includegraphics[scale=0.5]{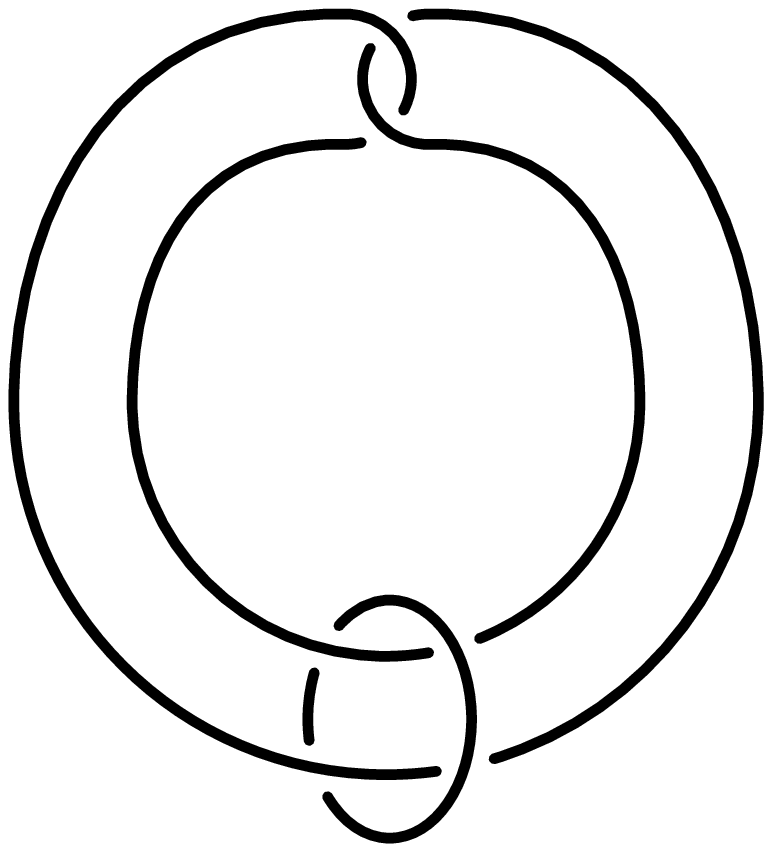}\includegraphics{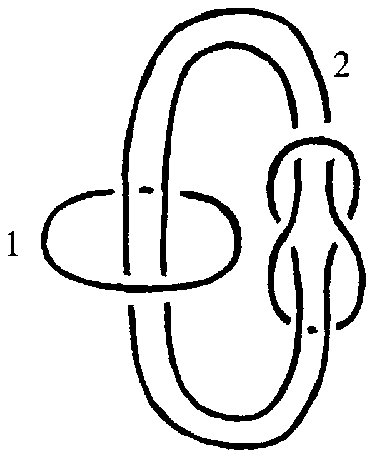}
\caption{Зацепление Хопфа, тривиальное зацепление и еще три зацепления}\label{f:hopf}
\end{figure}

\begin{figure}[h]\centering
\includegraphics{3-7.eps}
\caption{Кольца Борромео, зацепление Уайтхеда и трилистник}\label{f:borwhitre}
\end{figure}

\begin{assertion}\label{a:excha} (a) Зацепление Хопфа изотопно зацеплению, полученному из него перестановкой компонент.

(b) Зацепление Хопфа изотопно некоторому зацеплению, компоненты которого симметричны относительно некоторой прямой.

(c) Зацепление на рис. \ref{f:hopf} внизу слева изотопно зацеплению Уайтхеда на рис. \ref{f:borwhitre}.w.

(d,e) То же, что и в пп.~(a,b) для зацепления Уайтхеда.

(f) Кольца Борромео изотопны зацеплению, компоненты которого переставляются по циклу при повороте на угол $2\pi/3$ относительно некоторой прямой.
\end{assertion}

\begin{proof} (a) Это следует из п.~(b) (или может быть доказано независимо).

(d) Это следует из п.~(e) (или может быть доказано независимо).

\begin{figure}[h]\centering
\includegraphics[scale=0.9]{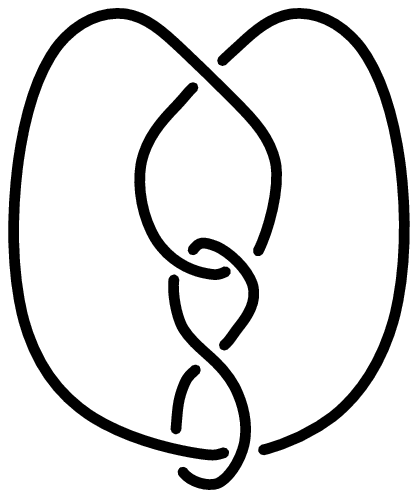}\qquad
\includegraphics[scale=0.9]{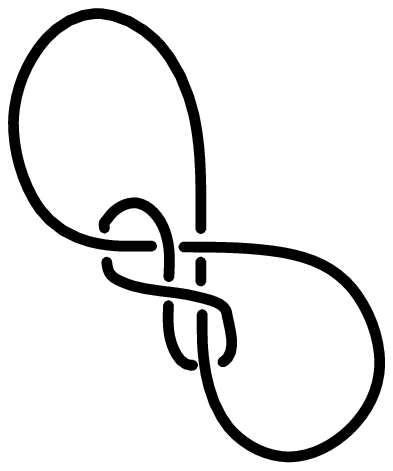}\qquad
\includegraphics[scale=0.9]{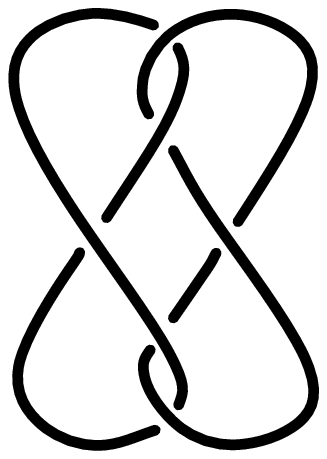}\qquad
\caption{Изотопия зацепления Уайтхеда}
	\label{f:whi}
\end{figure}

(e) См. рис. \ref{f:whi}.

(f) Рассмотрим цветные четырехугольники на рис. \ref{f:borr} слева.
Искомая прямая является биссектрисой любого октанта, образованного четырехугольниками.

\begin{figure}[h]\centering
\includegraphics[scale=0.5]{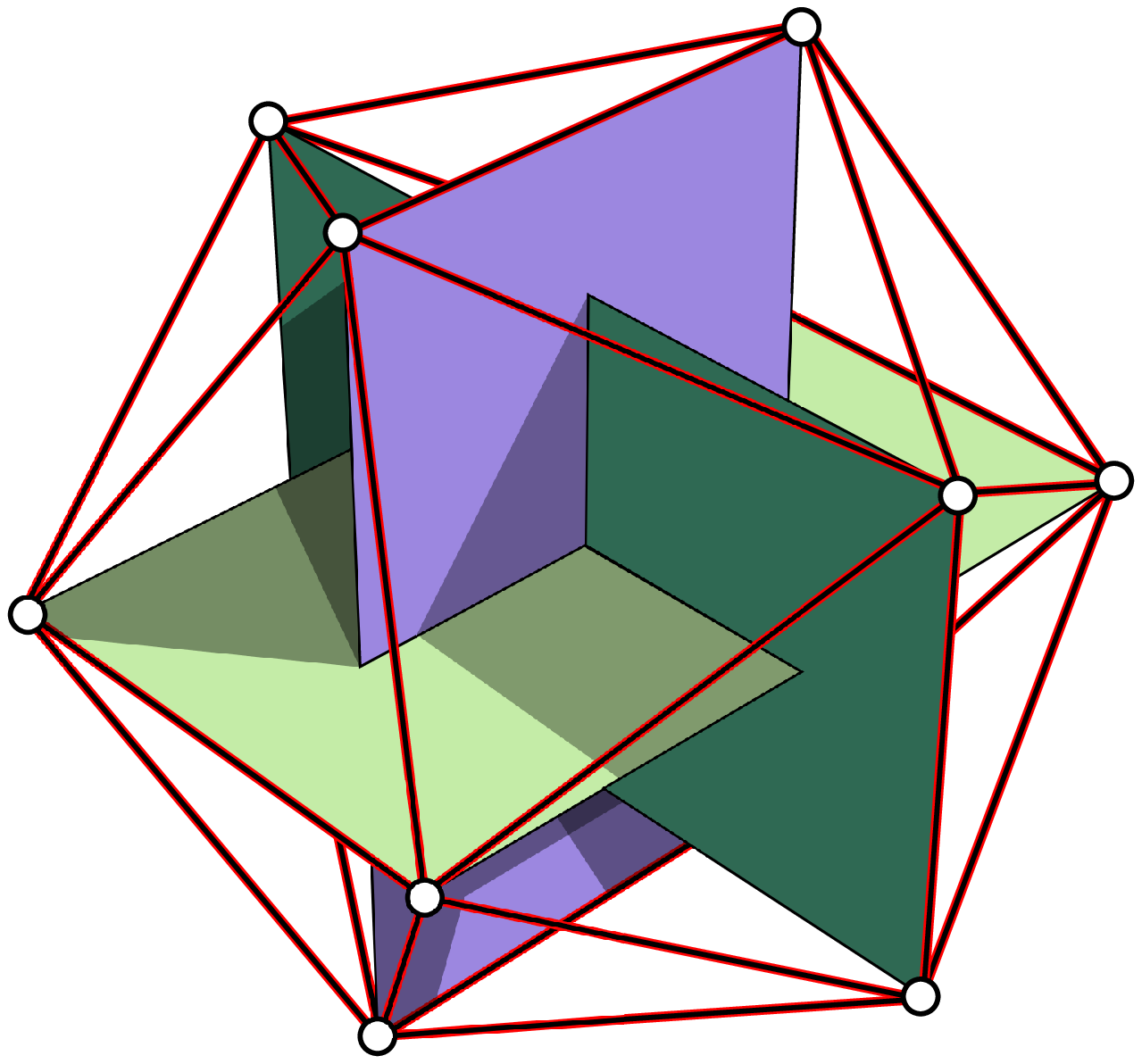}\qquad\includegraphics[scale=0.9]{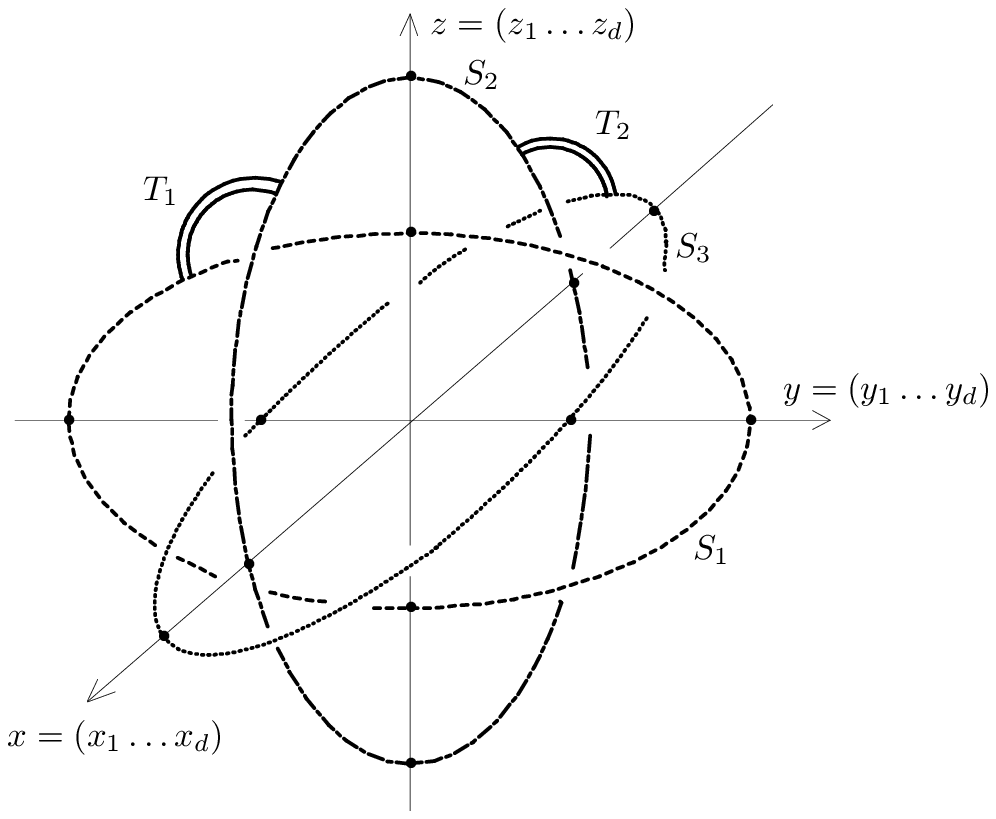}
\caption{Кольца Борромео}
\label{f:borr}
\end{figure}

Существует также следующая изящная криволинейная конструкция.
Возьмем три эллипса, заданные системами уравнений
$$\left\{\begin{array}{c} x=0\\ y^2+2z^2=1\end{array}\right., \qquad
\left\{\begin{array}{c} y=0\\ z^2+2x^2=1\end{array}\right. \qquad\text{и}\qquad
\left\{\begin{array}{c} z=0\\ x^2+2y^2=1. \end{array}\right.$$
(См. рис. \ref{f:borr} справа.)
Искомая прямая задается уравнением $x=y=z$.
\end{proof}

\begin{theorem}\label{t:hopf}
(a) Следующие зацепления попарно не изотопны: зацепление Хопфа, тривиальное зацепление, зацепление Уайтхеда.

(b) Кольца Борромео не изотопны тривиальному зацеплению.
\end{theorem}

Эта теорема доказывается с помощью {\it коэффициента зацепления по модулю 2}, придумайте его самостоятельно или см. \S\ref{0rapl2},
и {\it полиномов Александера-Конвея}, см. \S\ref{0con}.
Также можно использовать <<тройной коэффициент зацепления>> (инвариант Масси-Милнора) и <<высший коэффициент зацепления>> (инвариант Сато-Левина)  \cite[\S4.4-\S4.6]{Sk}.

\begin{theorem}\label{rlt-cgslin}  Если никакие 4 из 6 точек пространства не лежат в одной плоскости, то существует пара зацепленных треугольников с вершинами в этих 6 точках. (Это означает, что внутренность первого треугольника пересекает контур второго треугольника ровно в одной точке.)
\end{theorem}

Эта теорема доказывается с помощью {\it коэффициента зацепления по модулю 2}, см. \S\ref{0rapl2}.
Детали выходят за рамки данного текста, см. \cite[теорема~1]{CG83} и простое изложение в \cite[теорема 1.1]{Sk14}.

%Конвей--Гордон, Закс;

\section{Основные инструменты}\label{s:tools}

\begin{remark}[некоторые строгие рассуждения]\label{r:accur} В следующем абзаце мы доказываем, что \emph{если узел лежит в плоскости, то он изотопен тривиальному узлу.}

Обозначим узел, лежащий в плоскости, через $M_1M_2\ldots M_n$.
Возьмем точку $Z$ вне этой плоскости.
Тогда $M_1M_2\ldots M_n$ преобразуется в тривиальный узел $M_1ZM_n$ следующей последовательностью элементарных движений:
$$M_1M_2\to M_1ZM_2,\quad ZM_2M_3\to ZM_3,\quad ZM_3M_4\to ZM_4,\quad\ldots,\quad ZM_{n-1}M_n\to ZM_n.$$

Следуюший результат показывает, что промежуточные узлы изотопии, переводящей узел, лежащий в плоскости, в тривиальный узел, можно выбрать также лежащими в плоскости.

\emph{Теорема Ш\"енфлиса.} Любая ломаная без самопересечений в плоскости изотопна (в этой плоскости) треугольнику.

Это усиленная версия следующего выдающегося результата.

\emph{Теорема Жордана.} Любая замкнутая несамопересекающаяся ломаная $L$ в плоскости $\R^2$ разбивает эту плоскость в точности на две части, т.е. $\R^2-L$ --- несвязное множество, являющееся объединением двух связных множеств.

Подмножество плоскости называется \emph{связным}, если любые две точки этого подмножества можно соединить ломаной, в нем лежащей.

Алгоритмическое объяснение почему теорема Жордана (а значит и теорема Ш\"енфлиса) нетривиальна, и доказательство теоремы Жордана см. в \S1.3 <<Число пересечений для ломаных на плоскости>> \cite{Sk18}, \cite{Sk}.
\end{remark}

\begin{assertion}\label{a:arfmot} Предположим, что существует такая точка узла, что если мы идем вдоль узла, начиная с этой точки, то на некоторой диаграмме мы сначала встречаем только переходы, а затем только проходы.
Тогда этот узел изотопен тривиальному узлу.\footnote{Это утверждение будет мотивировкой для введения Арф-инварианта (\S\ref{0arf}).
Доказательство иллюстрирует в малых размерностях основные идеи выдающегося доказательства Зимана многомерной теоремы о незаузленности сфер, см. обзор \cite[теорема 2.3]{Sk16c}.}
\end{assertion}

\begin{figure}[h]\centering
\includegraphics[scale=0.5]{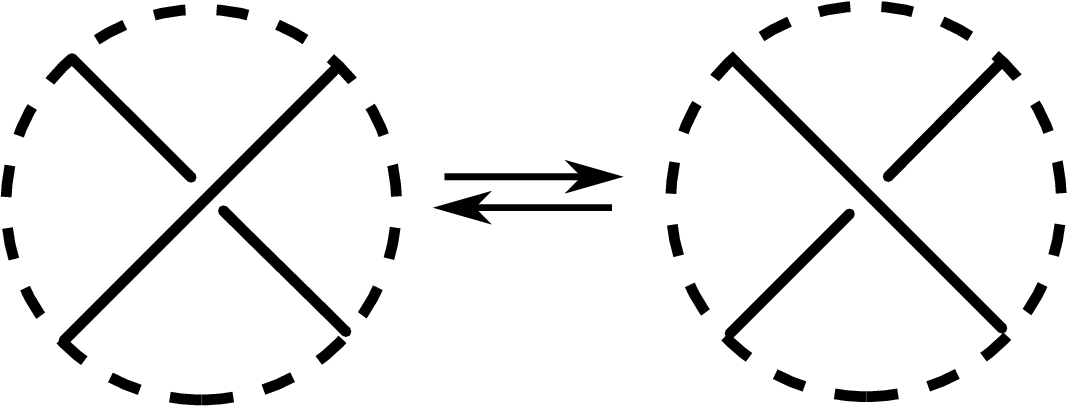}
\caption{Изменение перекрестка.
Плоские диаграммы совпадают вне кругов, ограниченных пунктирном.
Никакие другие ребра плоских диаграмм, кроме показанных на рисунке, не пересекаются с этими кругами.
(Те же соглашения приняты на рис. \ref{f:reid}, \ref{f:planiso}, \ref{f:skeino}, \ref{f:pass} и \ref{f:skein}.)}
\label{f:crossing}
\end{figure}

{\bf Изменение перекрестка} --- это замена перехода на проход или наоборот, см. рис. \ref{f:crossing}.

Очевидно, после любого изменения перекрестка в крайних слева диаграммах трилистника и восьмерки, показанных на рис. \ref{f:trefoil},
получится диаграмма узла, изотопного тривиальному.

\begin{lemma}\label{l:crossing} Любая диаграмма узла может быть преобразована изменениями перекрестков в диаграмму узла, изотопного тривиальному.\footnote{Эта простая лемма будет использована для построения инвариантов, использующего рекурсию ({\it скейн-соотношения}), см. теоремы \ref{con-lk2},   \ref{con-thma},\ref{con-lk2},\ref{con-thm},\ref{con-ale} и \ref{vas-kon}.}
\end{lemma}

\begin{figure}[h]\centering
\includegraphics[scale=0.45]{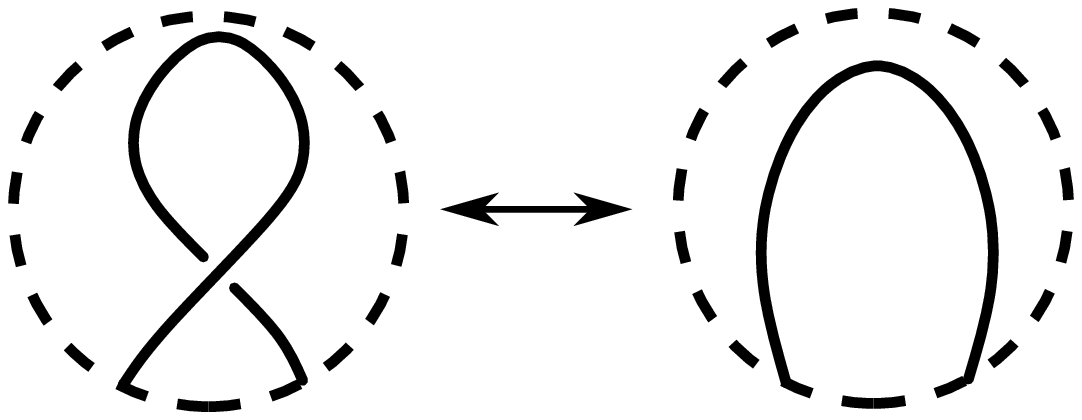}\qquad
\includegraphics[scale=0.45]{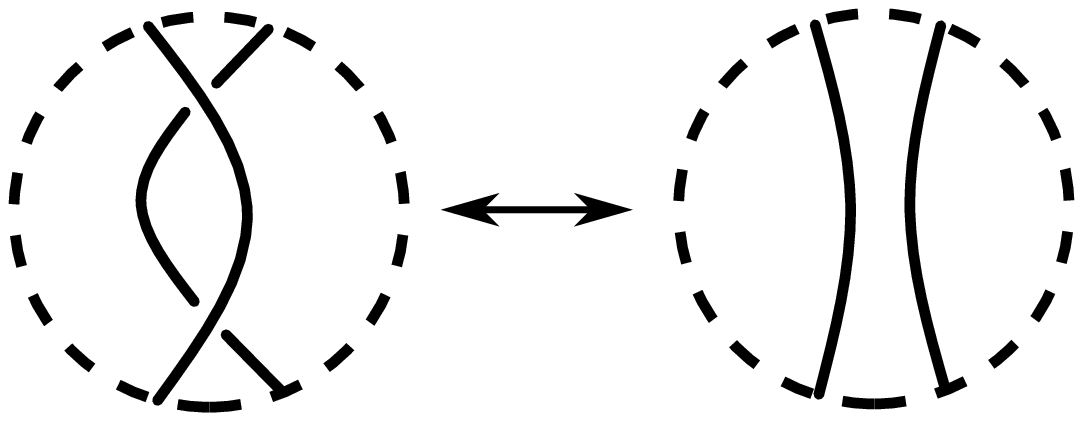}\qquad
\includegraphics[scale=0.45]{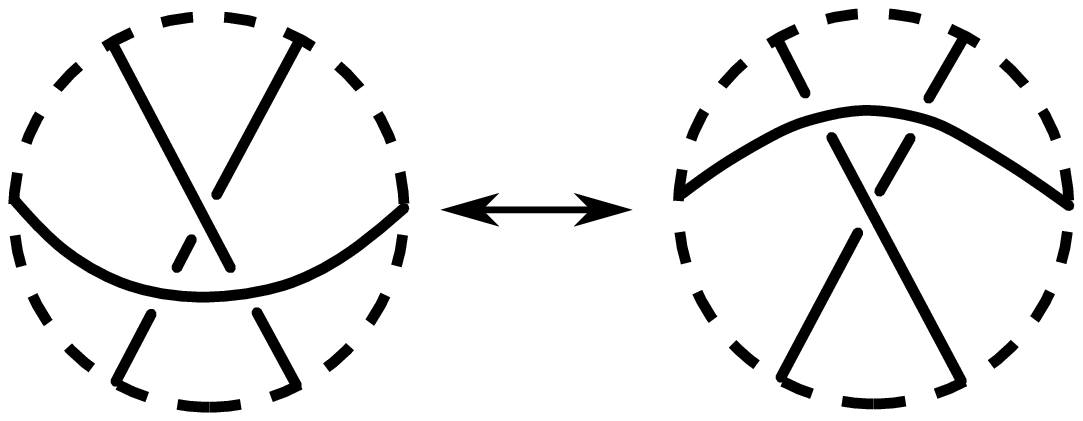}
\caption{Движения Рейдемейстера}
\label{f:reid}
\end{figure}

\begin{figure}[h]\centering
\includegraphics[scale=0.45]{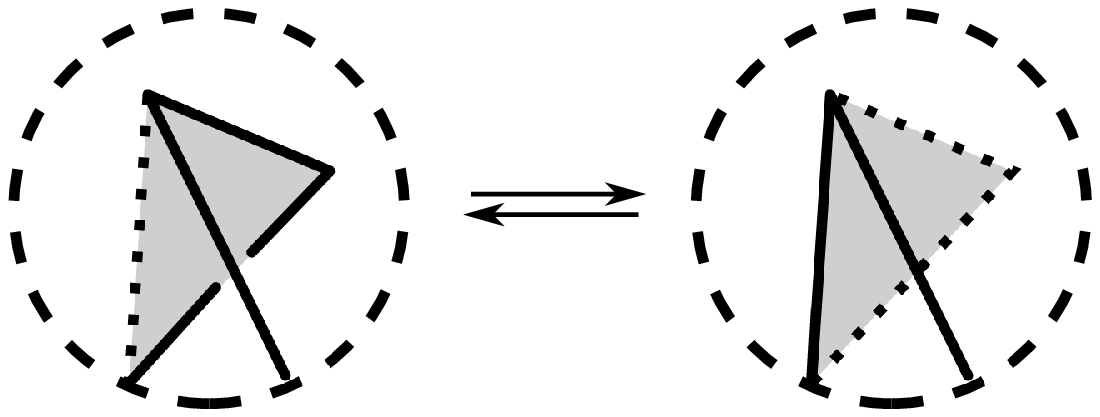}\qquad\qquad
\includegraphics[scale=0.45]{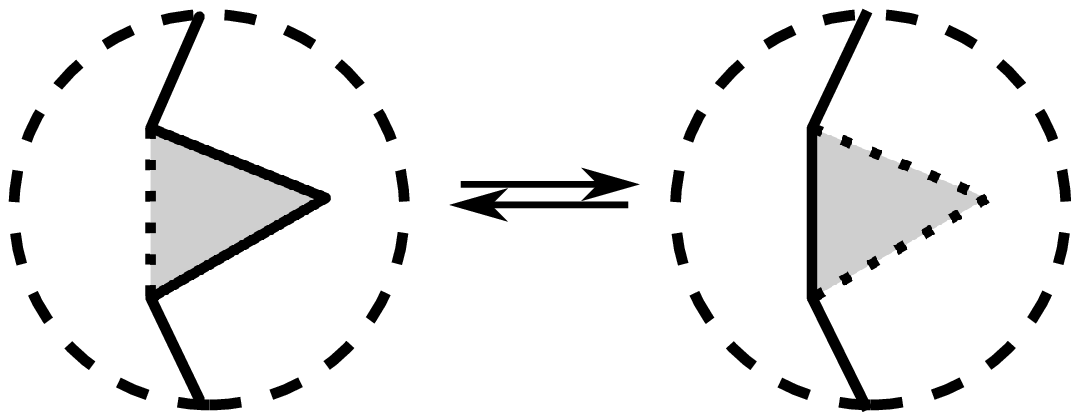}\quad
\includegraphics[scale=0.45]{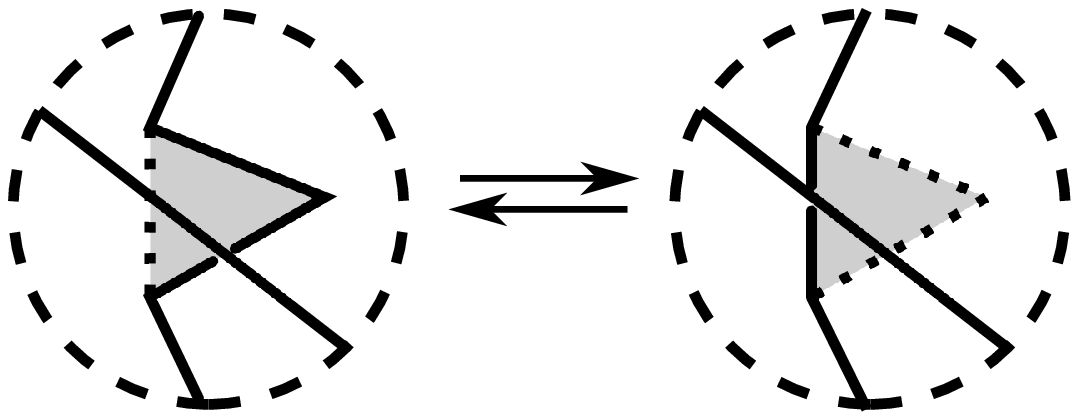}
\caption{Слева: к строгому определению первого движения Рейдемейстера
\newline
В середине, справа: движения изотопии плоскости}
\label{f:planiso}
\end{figure}

В этом тексте вместо изучения узлов с точностью до изотопии, мы будем изучать диаграммы узлов с точностью до эквивалентности, порожденной {\bf движениями Рейдемейстера} (см. рис. \ref{f:reid})\footnote{Строгое определение первого движения Рейдемейстера легко дается с помощью рис. \ref{f:planiso} слева. У других движений Рейдемейстера есть аналогичные строгие определения.
Вы можете использовать неформальное описание движений Рейдемейстера на рис. \ref{f:reid} и, таким образом, не пользоваться движениями изотопии плоскости.
См. сноску \ref{fo:reid}.}
и {\it движениями изотопии плоскости} (см. рис. \ref{f:planiso} в середине и справа).
Т.е. мы будем использовать без доказательства следующий результат.

\begin{theorem}[Рейдемейстер]\label{l:reid} Два узла изотопны тогда и только тогда, когда некоторая диаграмма первого узла может быть получена из некоторой диаграммы второго узла с помощью движений Рейдемейстера и движений изотопии плоскости.
\end{theorem}

См. \cite[\S1.7]{PS96}.\footnote{\label{fo:reid}
Так как в \cite[\S1.6]{PS96} нет такого же строгого определения для движений Рейдемейстера как для движений изотопии плоскости, рассуждения в \cite[\S1.7]{PS96} не являются строгим доказательством.
По-видимому, строгое доказательство может быть получено с использованием строгого определения движений Рейдемейстера.
\newline
Это также показывает, что наличие движений изотопии плоскости в утверждении из \cite[\S1.7]{PS96} не делает само утверждение строгим, поэтому их следует избегать.
Работая на интуитивном уровне, движениями изотопии плоскости лучше не пользоваться.
С помощью приведенного ниже альтернативного строгого определения движения изотопии плоскости выражаются через движения Рейдемейстера,
а потому не должны упоминаться в формулировке теоремы.
\newline
Приведем альтернативное строгое определение первого движения Рейдемейстера.
У других движений Рейдемейстера есть аналогичные строгие определения.
Возьмем на плоскости замкнутую несамопересекающуюся ломаную $L$, внутренность которой (см. теорему Жордана в замечании \ref{r:accur}) пересекает диаграмму $D$ по несамопересекающейся ломаной $M$, соединяющей две точки на $L$.
Обозначим через $N$ замкнутую несамопересекающуюся ломаную во внутренности ломаной $L$, для которой
$N\cap L=\emptyset$, $N\cap M$ --- одна точка, а $M\cup N$ --- (самопересекающаяся) ломаная общего положения.
{\it Первое движение Рейдемейстера} --- это замена $M$ на $M\cup N$ в $D$ с любой <<информацией>> в появившемся перекрестке.}
Справедливы также аналоги леммы \ref{l:crossing} и теоремы \ref{l:reid} для зацеплений.

%получается очень запутанно

\section{Гауссов коэффициент зацепления по модулю 2}\label{0rapl2}

%для плоских диаграмм

Предположим, что существует изотопия между двумя двухкомпонентными зацеплениями, и вторая компонента фиксирована в процессе изотопии.
Тогда след первой компоненты --- самопересекающийся цилиндр, не имеющий общих точек со второй компонентой.
Если после изотопии компоненты расцеплены, то цилиндр можно дополнить до самопересекающегося диска, не имеющего общих точек со второй компонентой.
Это наблюдение, вместе с \cite[лемма 4.3.2 о проекции]{Sk}, мотивирует следующее определение.

{\bf Коэффициентом зацепления по модулю два} $\lk_2$ плоской диаграммы двухкомпонентного зацепления называется взятое по модулю 2 число таких перекрестков диаграммы, для которых первая компонента проходит над второй компонентой.

\begin{pr}\label{a:lk2} Найдите коэффициент зацепления по модулю 2 для диаграмм на рис. \ref{f:hopf} и
для пар колец Борромео на рис. \ref{f:borwhitre}.b.
\end{pr}

\begin{lemma}\label{l:lk2} Коэффициент зацепления по модулю 2 сохраняется при движениях Рейдемейстера.
\end{lemma}

Эта лемма легко доказывается отдельно для каждого движения Рейдемейстера.

По лемме \ref{l:lk2} {\bf коэффициент зацепления по модулю 2} двухкомпонентного зацепления (или даже его изотопического класса) можно корректно определить как коэффициент зацепления по модулю 2 любой диаграммы этого зацепления.

Будем использовать без доказательства следующую {\it лемму о четности:} любые две замкнутые ломаные на плоскости, вершины которых находятся в общем положении, пересекаются в четном числе точек. Обсуждение и доказательство этой леммы можно найти в \S1.3 <<Число пересечений для ломаных на плоскости>> \cite{Sk18, Sk}.

\begin{assertion}\label{a:lk2-pr}
(a) Перестановка компонент двухкомпонентного зацепления сохраняет коэффициент зацепления по модулю 2.

(b) Существует двухкомпонентное зацепление, не изотопное тривиальному зацеплению,
но имеющее нулевой коэффициент зацепления по модулю 2.
\end{assertion}

П. (b) доказывается с помощью {\it целочисленного} коэффициента зацепления, см. \S\ref{0rapl}.

\begin{theorem}\label{con-lk2} Существует единственный изотопический инвариант $\lk_2$ двухкомпонентных зацеплений, принимающий значения в $\Z_2$, равный нулю на тривиальном зацеплении и такой, что
для любых зацеплений $K_+$ и $K_-$, диаграммы которых отличаются изменением перекрестка $A$,
$$\lk\phantom{}_2K_+-\lk\phantom{}_2K_-=\begin{cases}1 &A\text{ --- пересечение различных компонент;} \\
0 &A\text{ --- самопересечение одной из компонент.}\end{cases}$$
\end{theorem}

\begin{assertion}\label{a:lk2comp} Если коэффициент зацепления по модулю 2 двух (непересекающихся контуров) треугольников в пространстве равен нулю, то зацепление, образованное этими треугольниками, изотопно тривиальному.
\end{assertion}

Доказательство, по-видимому, не опубликовано, однако оно не очень сложное.
Мы призываем читателя опубликовать подробное рассуждение. Ср. \cite{Ko19}.

\section{Арф-инвариант}\label{0arf}

Рассмотрим плоскую диаграмму некоторого узла и ее точку $P$, не являющуюся перекрестком.
Назовем $P$ {\it отмеченной точкой}.
Неупорядоченная пара перекрестков $A$, $B$ называется {\bf скрещивающейся} (или $P$-скрещивающейся),
если, обходя диаграмму в каком-нибудь направлении, начиная с $P$, и отмечая только
перекрестки в $A$ и в $B$, мы сначала встретим переход в $A$, затем проход в $B$, затем проход в $A$
и, наконец, переход в $B$.
Далее, {\it $P$-Арф-инвариантом} $\arf_P$ плоской диаграммы называется четность числа ее $P$-скрещивающихся пар перекрестков.

\begin{pr}\label{e:con-cas} (a) Если $P$-Арф-инвариант плоской диаграммы не равен нулю, то $P$ не может быть точкой из утверждения \ref{a:arfmot}.

(b) Найдите $P$-Арф-инвариант (какой-нибудь плоской диаграммы) тривиального узла, трилистника и восьмерки
(при Вашем выборе отмеченной точки $P$).
\end{pr}

\begin{lemma}\label{p:con-cas-ab} (a) $P$-Арф-инвариант не зависит от выбора отмеченной точки $P$.

(b) Арф-инвариант плоской диаграммы сохраняется при движениях Рейдемейстера.
\end{lemma}

Согласно п.~(a) {\it Арф-инвариант} плоской диаграммы можно корректно определить как ее $P$-Арф-инвариант для произвольной отмеченной точки $P$.
Поэтому формулировка п.~(b) осмыслена.
Согласно п.~(b) {\it Арф-инвариант} $\arf$ узла (или даже изотопического класса узлов) можно корректно определить как Арф-инвариант произвольной плоской диаграммы этого узла.

\emph{Указания.} (a) Будем перемещать отмеченную точку вдоль плоской диаграммы. Тогда достаточно показать, что Арф-инвариант не изменяется, когда отмеченная точка проходит через перекресток.

(b) Докажите утверждение для каждого движения Рейдемейстера отдельно, всякий раз выбирая отмеченную точку наиболее удобным образом.

\begin{assertion}\label{a:arf-pr} Существует узел, не изотопный тривиальному, Арф-инвариант которого равен нулю.
\end{assertion}

Это доказывается с помощью \emph{инварианта Кэссона}, см.~\S\ref{0cas}.

\begin{figure}[h]\centering
\includegraphics[scale=0.6]{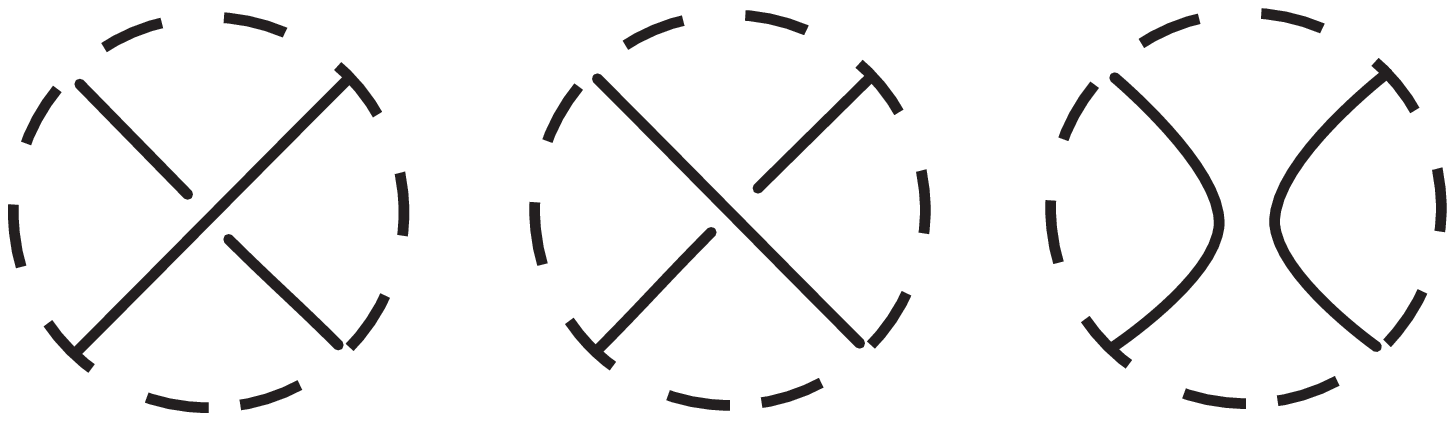}
\caption{Зацепления $K_+,K_-,K_0$}
\label{f:skeino}
\end{figure}

\begin{theorem}\label{con-thma} Существует единственный изотопический инвариант $\arf$ узлов, принимающий значения в $\Z_2$, равный нулю на тривиальном узле и такой, что
$$\arf K_+-\arf K_-=\lk\phantom{}_2K_0$$
для любых узлов $K_+$ и $K_-$, диаграммы которых отличаются как показано на рис.~\ref{f:skeino}, причем $K_0$ является двухкомпонентным зацеплением. (Последнее эквивалентно существованию ориентации, для которой   рис. \ref{f:skeino} превращается в рис. \ref{f:skein}.)
\end{theorem}

\begin{figure}[h]\centering
\includegraphics[scale=0.8]{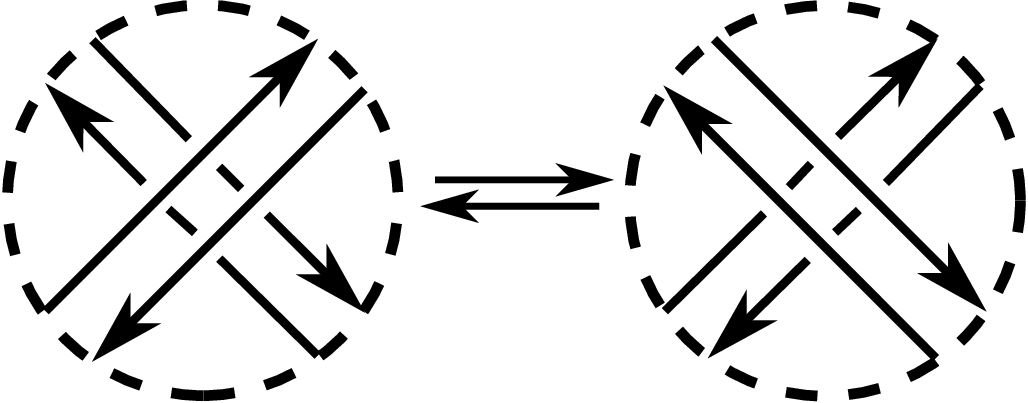}
\caption{Пронос}
\label{f:pass}
\end{figure}

\begin{assertion}\label{lk-pamo} Два узла называются \emph{Арф-эквивалентными}, если какая-нибудь диаграмма
первого узла (с некоторой ориентацией) может быть преобразована в какую-нибудь диаграмму второго узла
(с некоторой ориентацией) движениями Рейдемейстера и с помощью \emph{проносов}, показанныx на рис. \ref{f:pass}.

(a) Если два узла Арф-эквивалентны, то их Арф-инварианты совпадают.

(b) Восьмерка Арф-эквивалентна трилистнику.

(c) \cite[стр. 75--78]{Ka87} Если у двух узлов Арф-инварианты совпадают, то эти узлы Арф-эквивалентны.
\end{assertion}

\section{Приложение: правильные раскраски}\label{s:reco}

Этот раздел иллюстрирует связь с раскрасками в комбинаторике \cite{Ra20}.
Он использует только материал разделов \ref{0isot} и \ref{0isoli}.
Подробнее про правильные раскраски см. в \cite{Pr98}.

{\it Дугой} на плоской диаграмме (узла или зацепления) называется связный фрагмент ломаной, идущий от одного прохода до следующего.
{\bf Правильная раскраска} плоской диаграммы (узла или зацепления) --- это такая раскраска ее дуг в три цвета, что по крайней мере два цвета использовано, и в каждом перекрестке встречаются либо все три цвета, либо только один цвет.
Плоская диаграмма (узла или зацепления) называется {\bf раскрашиваемой в три цвета} если ее можно правильно раскрасить.

\begin{pr}\label{wha-cole} Для каждого из следующих узлов или зацеплений возьмите любую плоскую диаграмму и определите, является ли она раскрашиваемой в три цвета:

(a) тривиальный узел;
\quad
(b) трилистник;
\quad
(c) восьмерка.

(d-i) зацепления на рис. \ref{f:hopf} и \ref{f:borwhitre}.b.
\end{pr}

\begin{lemma}[{\cite[стр. 29-30, теорема 4.1]{Pr95}}]\label{wha-colre}
Раскрашиваемость в три цвета плоской диаграммы сохраняется при движениях Рейдемейстера.
\end{lemma}

\begin{theorem}\label{wha-colr} (a) Ни одно из зацеплений на рис. \ref{f:hopf} и \ref{f:borwhitre} (кроме тривиального зацепления) не изотопно тривиальному зацеплению.

(b) Узел $5_1$ не изотопен тривиальному узлу.
\end{theorem}

\begin{figure}[h]\centering
\includegraphics[scale=0.9]{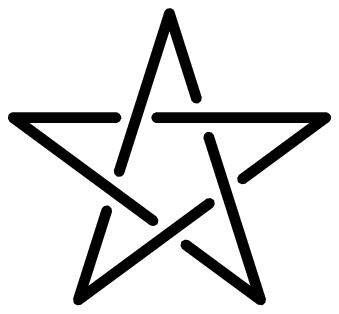}
\caption{Узел $5_1$}
\label{f:51}
\end{figure}

\section{Ориентированные узлы и зацепления; связные суммы}\label{s:ori}

\begin{figure}[h]\centering
\includegraphics[scale=0.5]{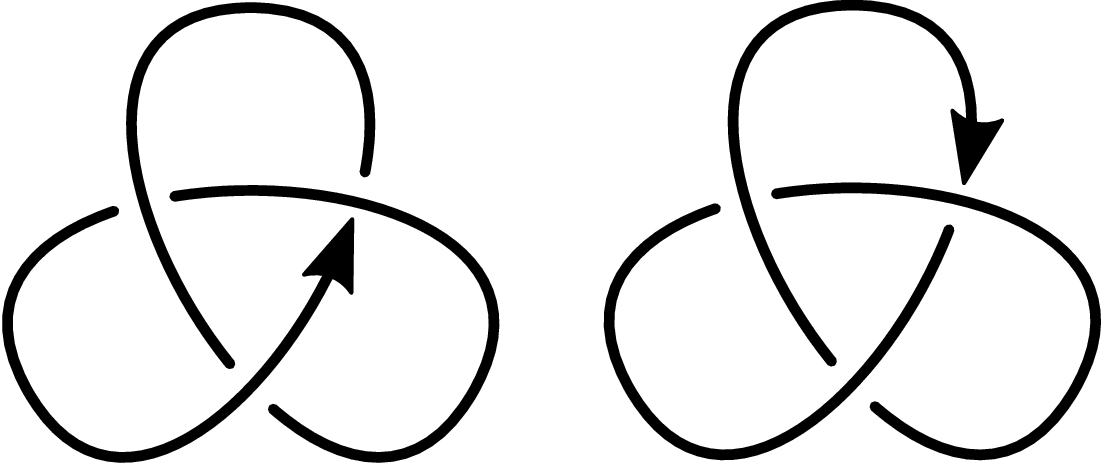}
\caption{Два трилистника с противоположными ориентациями}
\label{f:ori}
\end{figure}

Вы знаете, что такое ориентированная ломаная, так что вы также знаете, что такое ориентированный узел
(рис. \ref{f:ori}).

Неформальное понятие и строгое определение {\it изотопии} ориентированных узлов даются аналогично изотопии узлов.

\begin{assertion}\label{wha-chirsph} Изотопные ориентированные ломаные без самопересечений на плоскости и на сфере определяются аналогично изотопным ориентированным узлам в пространстве.

(a) Ориентированный сферический треугольник изотопен на сфере тому же треугольнику с противоположной ориентацией.

(b) Аналог пункта (a) для плоскости неверен.
\end{assertion}

\begin{assertion}\label{wha-inv} Следующие пары узлов с противоположными ориентациями изотопны: два тривиальных узла; два трилистника; две восьмерки.
\end{assertion}

\begin{theorem}[Троттер, 1964]\label{wha-ninv}
Существует ориентированный узел, не изотопный такому же узлу с противоположной ориентацией.
\end{theorem}

Этот факт доказывается с использованием {\it полиномов Джонса} \cite[\S3]{PS96}, \cite[\S2.4]{CDM}; доказательство выходит за рамки данного текста.

\begin{figure}[h]\centering
\includegraphics[scale=0.5]{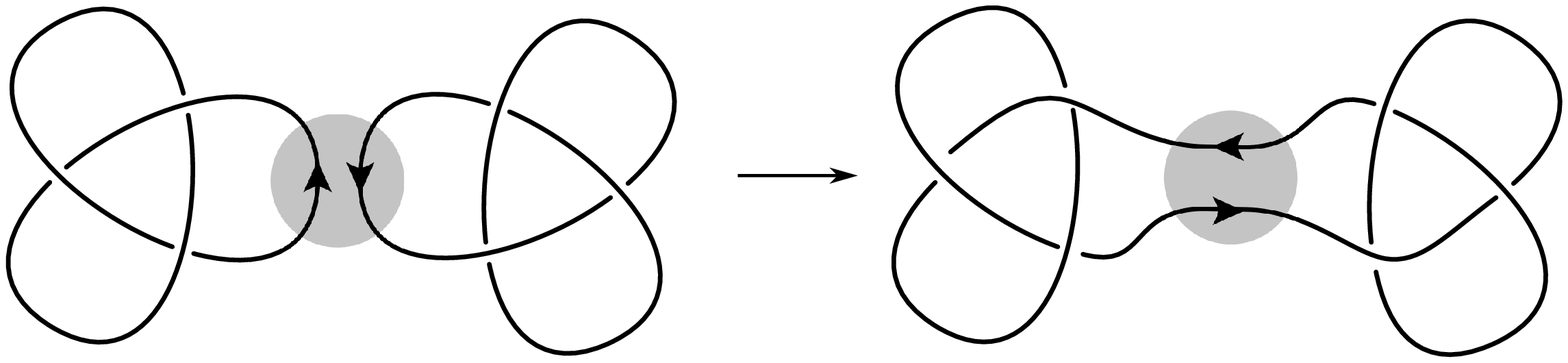}
\includegraphics[scale=0.4]{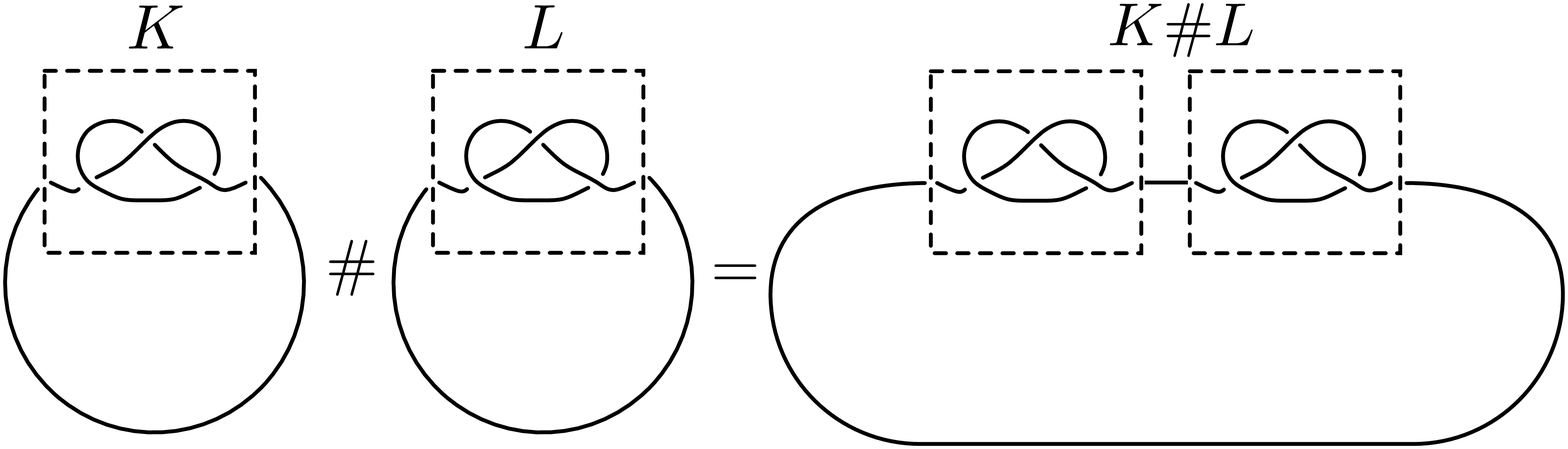}
\caption{Связная сумма узлов}
\label{f:consum}
\end{figure}

%{\bf Связная сумма} $\#$ {\it ориентированных} узлов определена на рис. \ref{f:consum}.
%Эта операция не является корректно определенной для ориентированных узлов.
%Так что мы обозначаем через $K\# L$ любую связную сумму узлов $K$ и $L$.

{\bf Связная сумма} $\#$ {\it ориентированных} узлов определена на рис. \ref{f:consum}.
{\bf Связная сумма} $\#$ (неориентированных) узлов определяется аналогично (забудьте про стрелки).
Ни та, ни другая операция не является корректно определенной на множестве узлов или ориентированных узлов соответственно; см. также замечание \ref{r:well}.
Так что мы обозначаем через $K\# L$ любую из связных сумм узлов или ориентированных узлов $K$ и $L$.

\begin{assertion}\label{wha-cons} Для любых ориентированных узлов $K,L,M$ и тривиального ориентированного узла $O$ имеем

(a) $K\#O=K$. \quad (b) $K\#L=L\#K$. \quad (c) $(K\#L)\#M=K\#(L\#M)$. \quad

(d) Для любых неориентированных узлов $K,L$ выполнено $\arf(K\#L)=\arf K+\arf L$.

(Строгий смысл пункта (a) --- <<существует связная сумма узлов $K$ и $O$, изотопная узлу $K$>>.
Строгий смысл пунктов (b,c) и (d) аналогичен.
См., впрочем, замечание \ref{r:well}.)
\end{assertion}

\begin{proof}[Набросок доказательства]
(a) См. рис. \ref{f:kok}.

\begin{figure}[h]\centering
\includegraphics[scale=0.7]{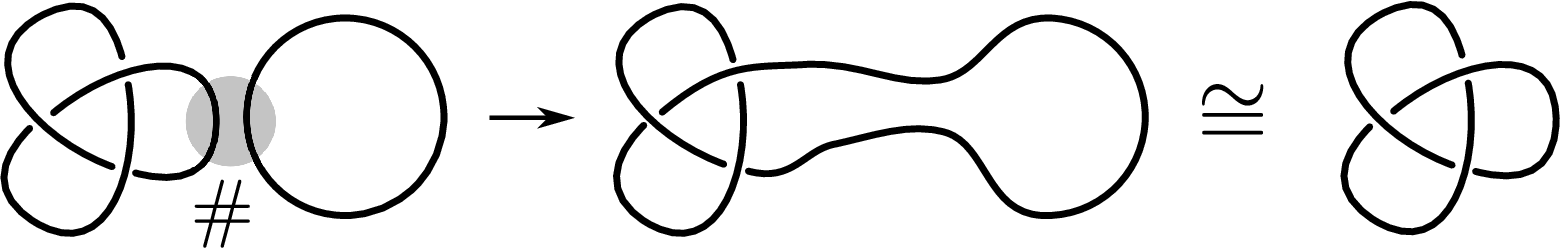}
\caption{Доказательство равенства $K\#O = K$}
\label{f:kok}
\end{figure}

\begin{figure}[h]\centering
\begin{tabular}{@{}cc@{}}
\includegraphics[scale=0.4]{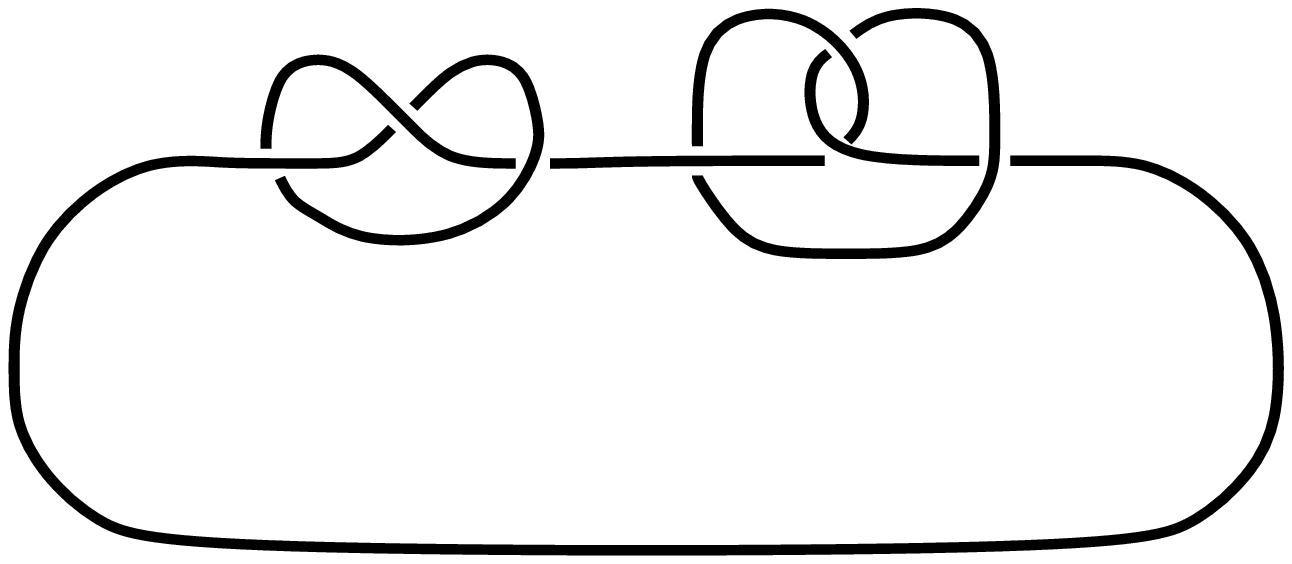}&\\
\includegraphics[scale=0.4]{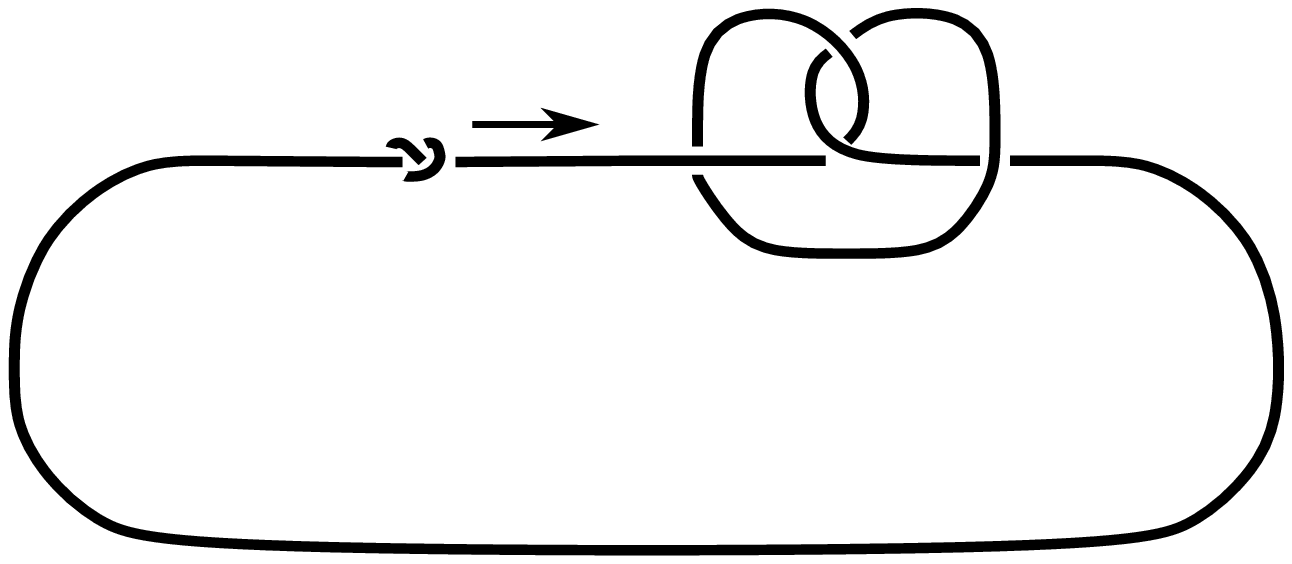}&\\
\includegraphics[scale=0.4]{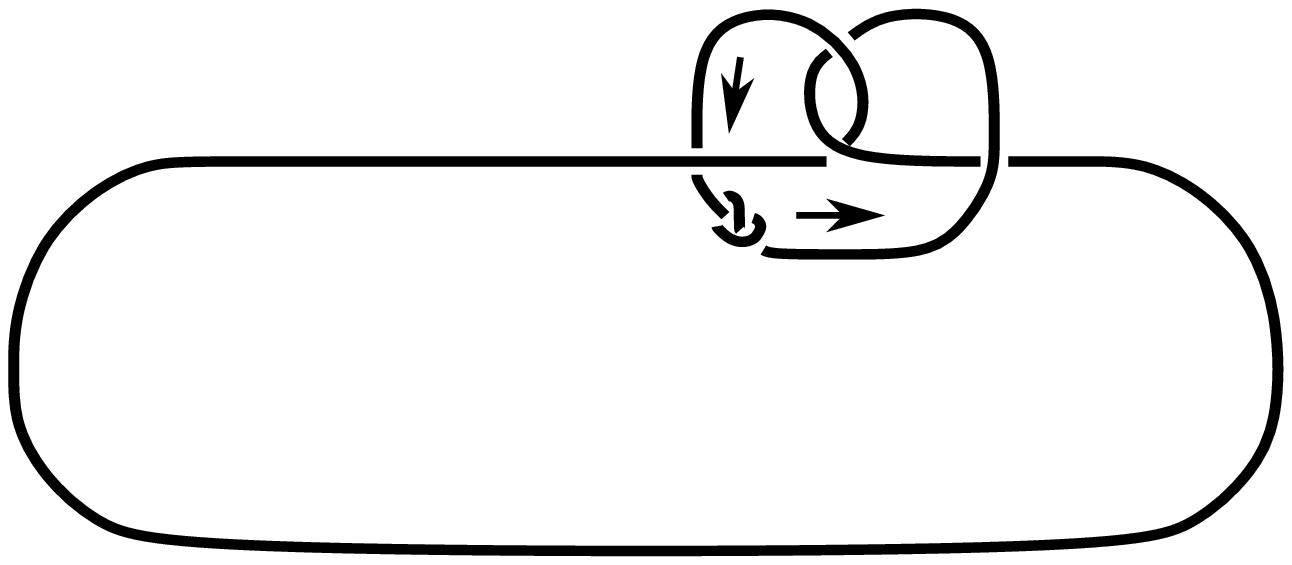}& \qquad \qquad \includegraphics[scale=0.4]{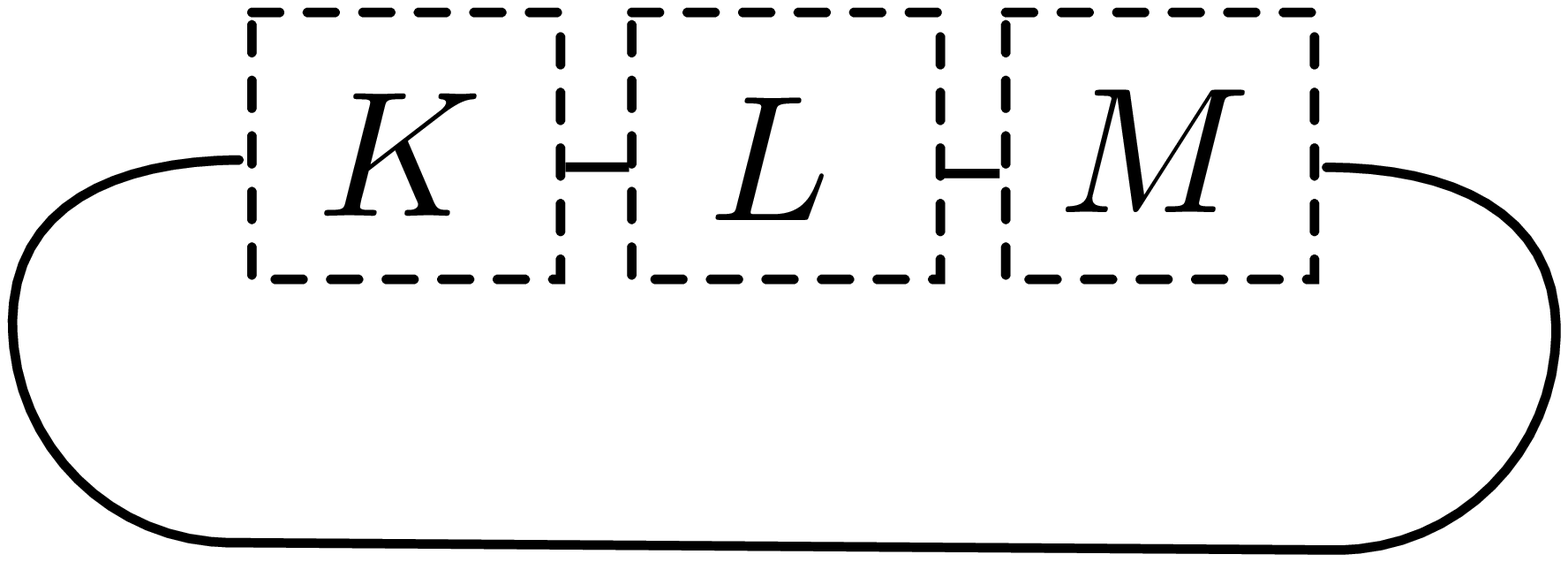}\\
\includegraphics[scale=0.4]{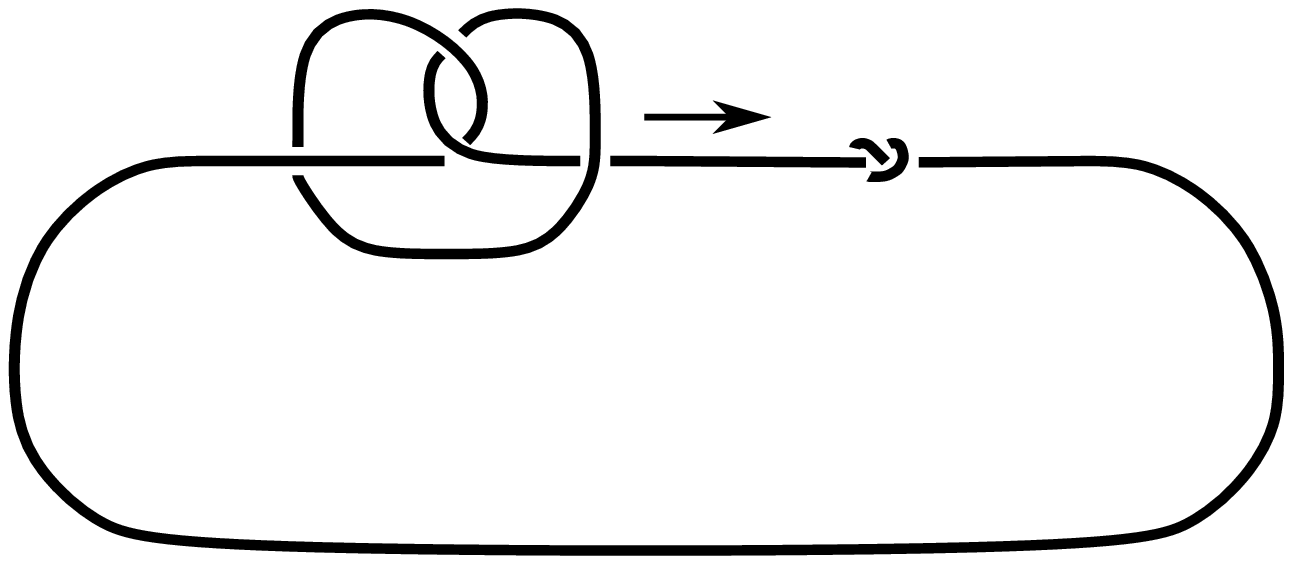}&\\
\includegraphics[scale=0.4]{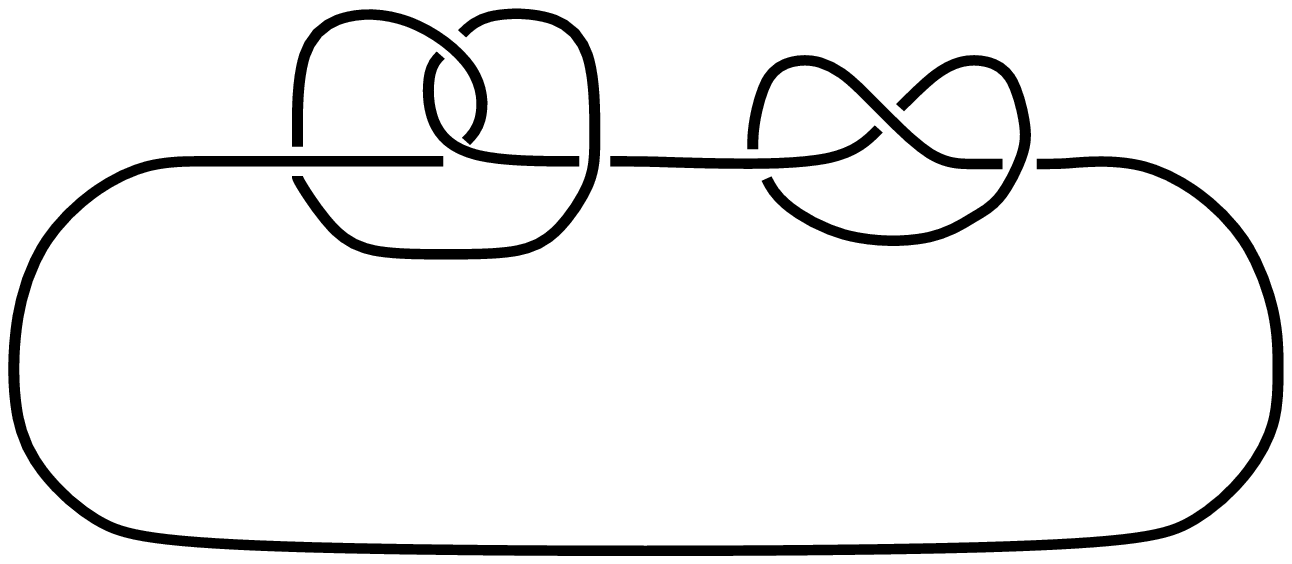}&\\
\end{tabular}
\caption{Доказательства равенств $K\#L = L\#K$ (слева) и $(K\#L)\#M = K\#(L\#M)$ (справа)}\label{f:comm}
\end{figure}

(b) Достаточно взять <<маленький>> узел, изотопный узлу $L$ и <<провести>> его через узел $K$, см. рис. \ref{f:comm} слева.

(c) Узлы из левой и правой частей равенства изотопны узлу на рис. \ref{f:comm} справа.

(d) Выберите отмеченную точку $P$ близко к <<месту соединения>> узлов.
Проверьте, что всякая $P$-скрещивающаяся пара перекрестков в $K\#L$ получается либо из скрещивающейся пары перекрестков в $K$, либо из скрещивающейся пары перекрестков в $L$.
\end{proof}

\begin{remark}\label{r:well} {\it Изотопический класс} узла --- это множество всех узлов, ему изотопных. Ориентированный изотопический класс $[K\#L]$ связной суммы двух ориентированных изотопических классов $[K],[L]$ ориентированных узлов $K,L$ не зависит от выборов, сделанных в процессе построения, и от выбора представителей $K,L$ классов $[K],[L]$.
Следовательно, связная сумма ориентированных изотопических классов ориентированных узлов корректно определена как $[K]\#[L]:=[K\#L]$, см. \cite[замечание 2.3.a]{Sk15}.
Для изотопических классов неориентированных узлов связная сумма не является корректно определенной операцией \cite{CSK}.
\end{remark}

\begin{theorem}\label{t:consum}
Пусть $K,L,M$ --- изотопические классы ориентированных узлов, а $O$ --- изотопический класс ориентированного тривиального узла. Тогда

(a) если $K\#L=O$, то $K=L=O$;

(b) если $K\#L=K\#M$, то $L=M$.
\end{theorem}

Доказательство выходит за рамки данного текста, см. \cite[теорема 1.5]{PS96}.
(В этом отрывке из \cite{PS96} необходимо всюду заменить <<узел>> на <<ориентированный узел>> ввиду замечания \ref{r:well}.)

\begin{figure}[h]\centering
\includegraphics[scale=0.9]{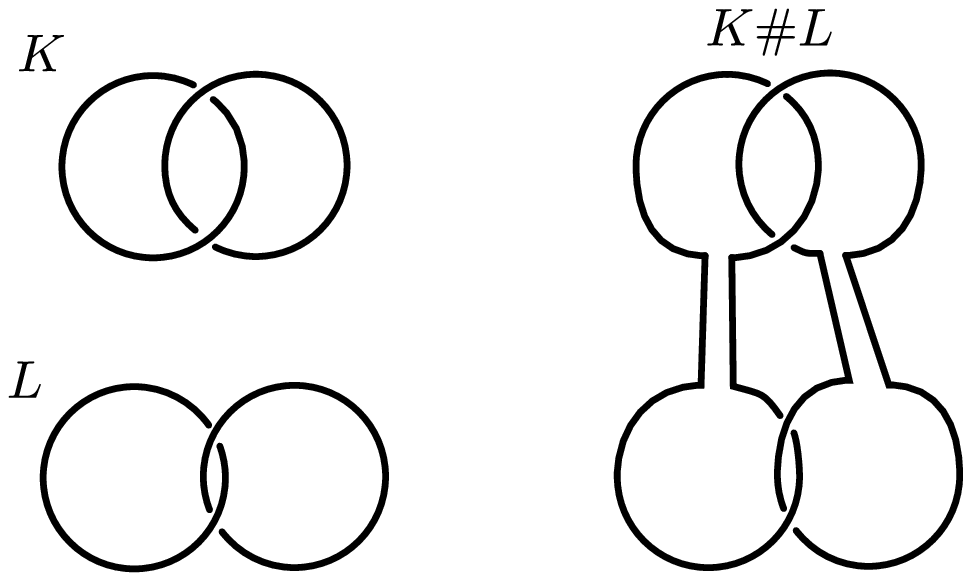}
\caption{Связная сумма зацеплений}
\label{f:consuml}
\end{figure}

Связная сумма $\#$ зацеплений (упорядоченных или нет, ориентированных или нет) определяется аналогично связной сумме узлов, см. рис. \ref{f:consuml}.
Эта операция не является корректно определенной для зацеплений, и замечание \ref{wha-csli} показывает, что эта операция также не дает корректно определенного понятия связной суммы изотопических классов зацеплений.
Так что мы обозначаем через $K\# L$ любую связную сумму зацеплений $K$ и $L$.

\begin{assertion}\label{wha-consli} (a,b,c) Докажите аналоги утверждений \ref{wha-cons}.a,b,c для зацеплений.

(d) Для любых неориентированных двухкомпонентных зацеплений $K,L$ выполнено $\lk_2(K\#L)=\lk_2K+\lk_2L$.
\end{assertion}

\begin{figure}[h]
\centering
\includegraphics[scale=0.5]{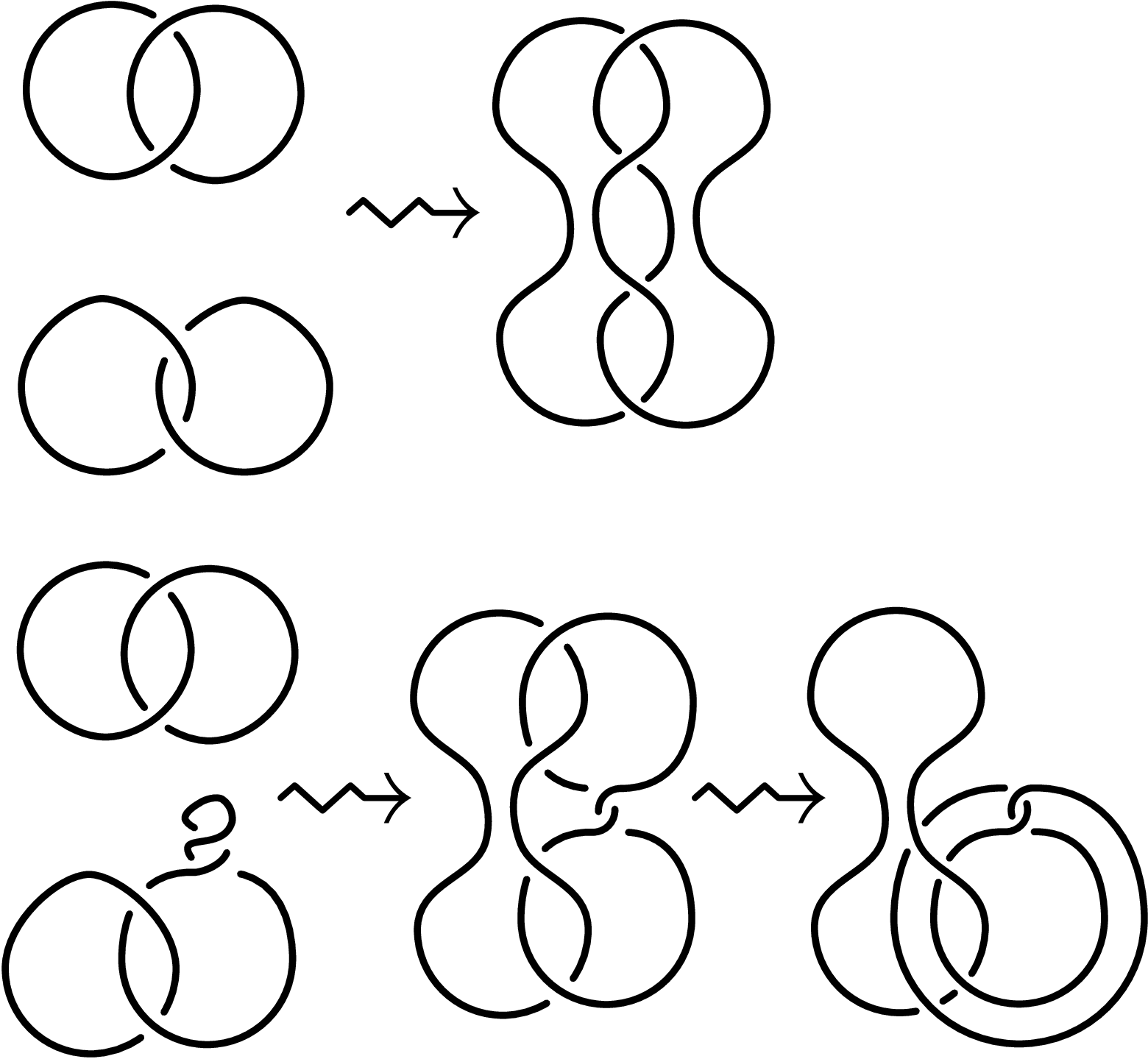}
\caption{Связная сумма изотопических классов зацеплений не определена корректно}
\label{f:ryabichev}
\end{figure}

\begin{remark}\label{wha-csli} Существуют две изотопные пары $(K,L)$ и $(K',L')$
двухкомпонентных зацеплений  такие, что некоторые связные суммы $K\#L$ и $K'\#L'$ не изотопны.
(Зацепления могут быть упорядочены или нет, ориентированы или нет --- всего получается 4 утверждения.)
Для случая неупорядоченных зацеплений пример дается парой одинаковых зацеплений, каждое из которых состоит из трилистника и тривиального узла, содержащихся в непересекающихся кубах, ср. \cite[рис. 3.16]{PS96}.
Пример для случая упорядоченных зацеплений содержится в \cite{As}.
На рис. \ref{f:ryabichev} показан другой пример, предложенный А.\,Рябичевым.
\end{remark}

\comment

\begin{figure}[h]\centering
\includegraphics[scale=0.4]{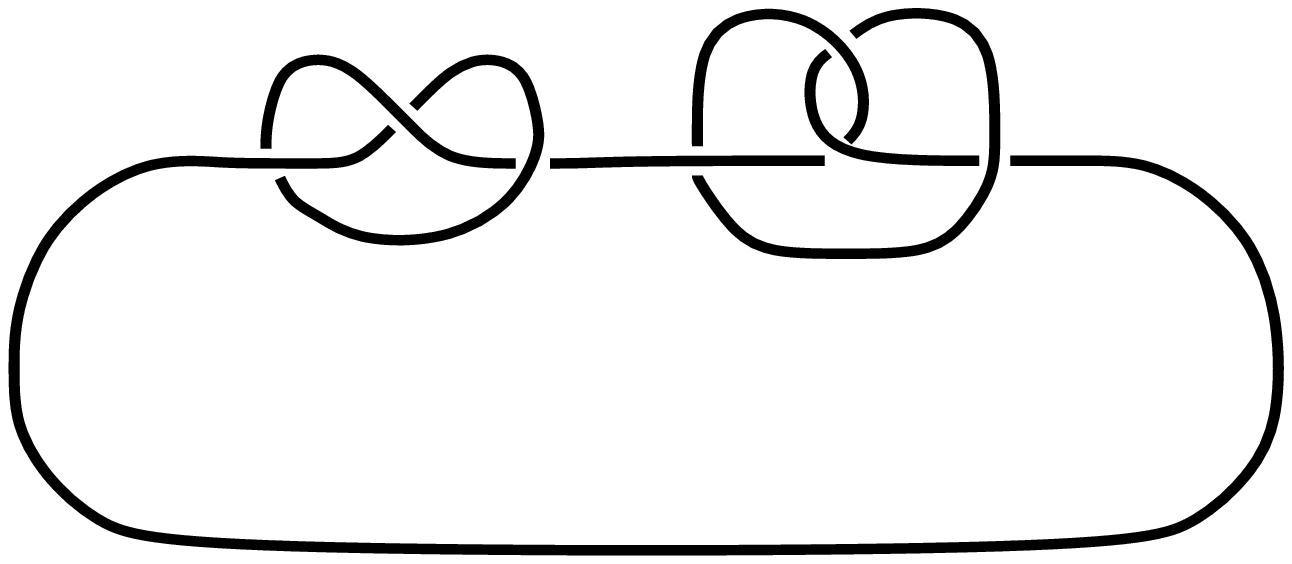}\\
\includegraphics[scale=0.4]{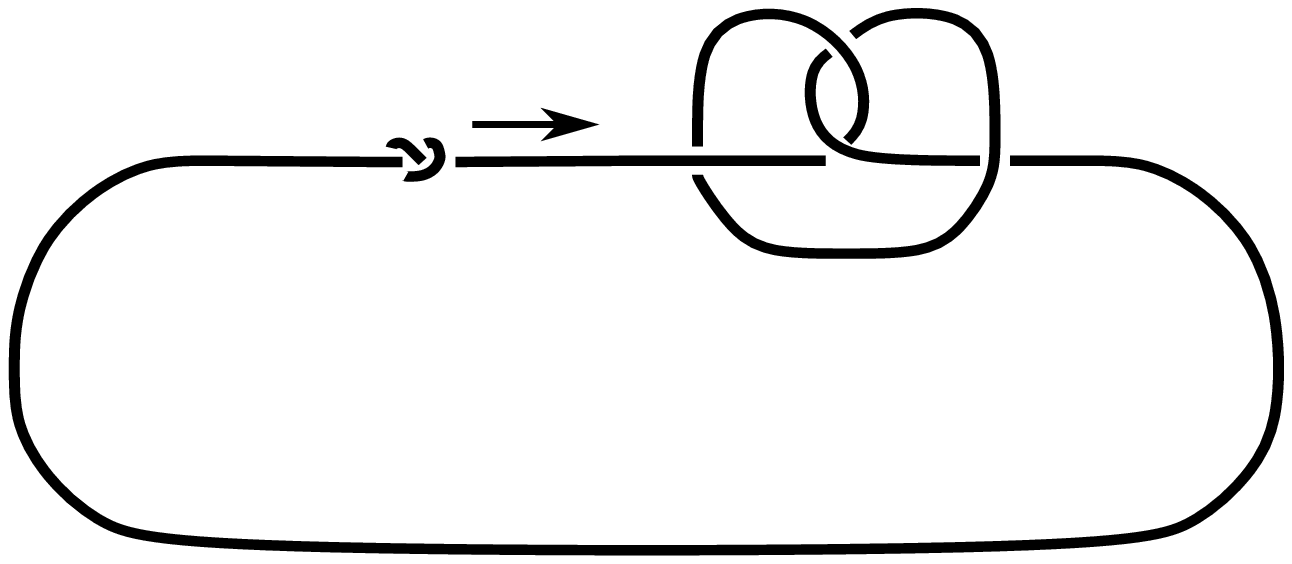}\\
\includegraphics[scale=0.4]{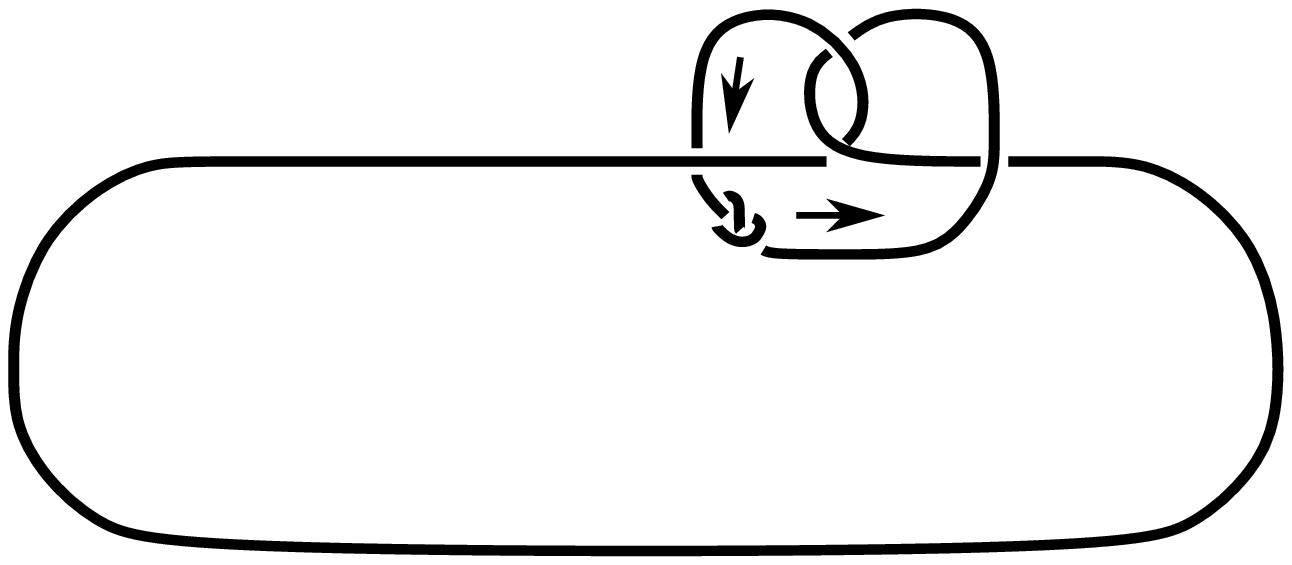}\\
\includegraphics[scale=0.4]{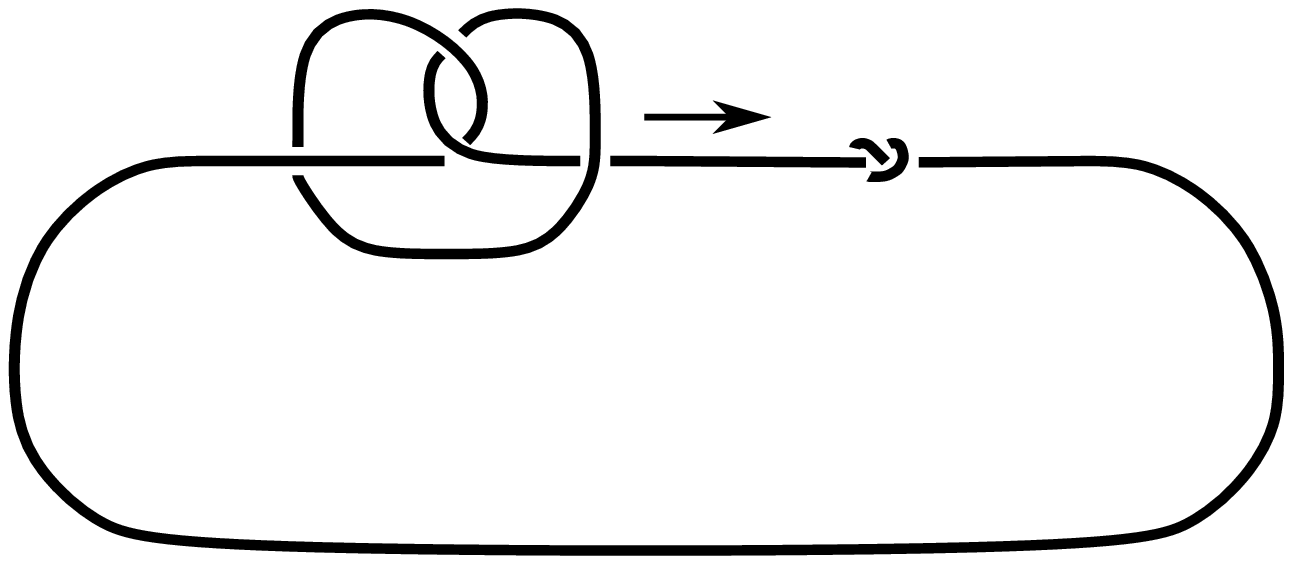}\\
\includegraphics[scale=0.4]{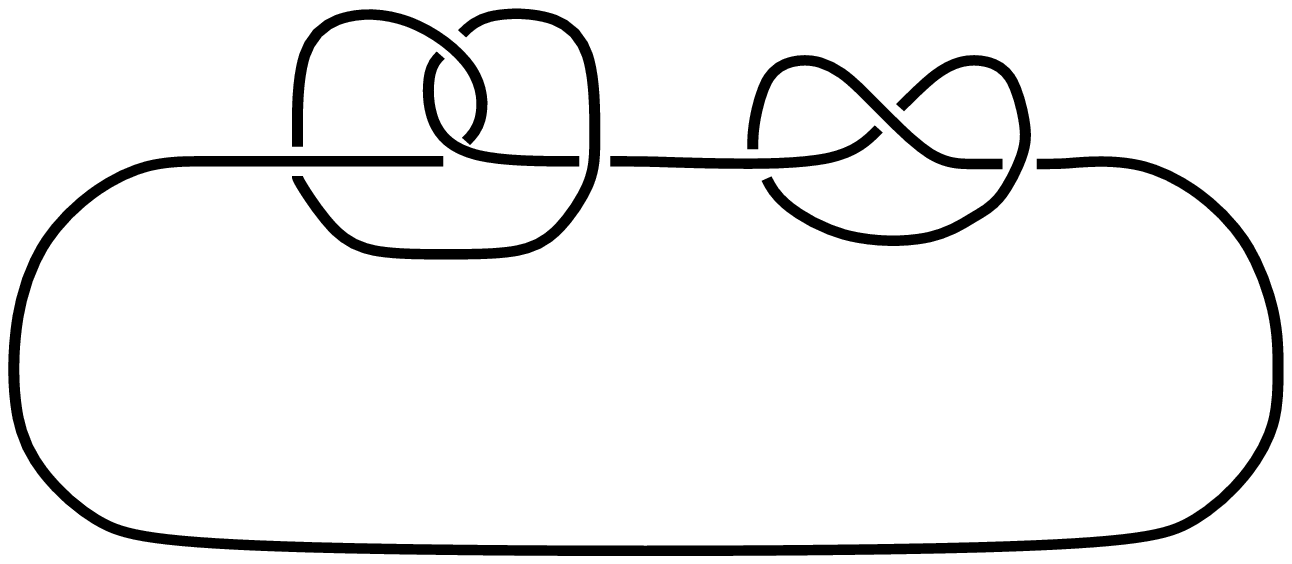}\\
\end{figure}

\begin{figure}[h]\centering
\includegraphics[scale=0.4]{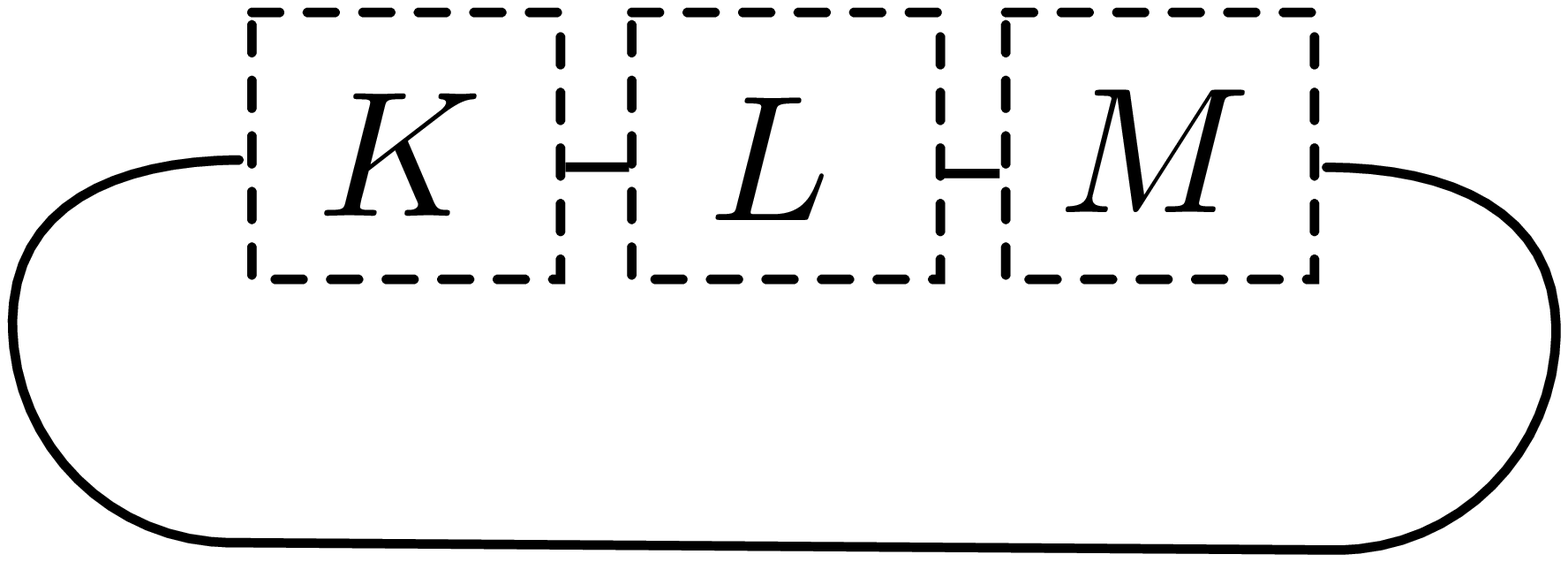}
\caption{Доказательство равенства $(K\#L)\#M = K\#(L\#M)$}
\label{f:asso}
\end{figure}

\endcomment

\section{Гауссов коэффициент зацепления}\label{0rapl}

%для плоских диаграмм

\begin{figure}[h]\centering
\includegraphics[scale=0.7]{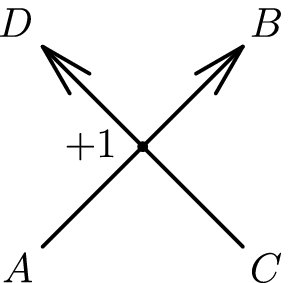}\quad \includegraphics[scale=0.7]{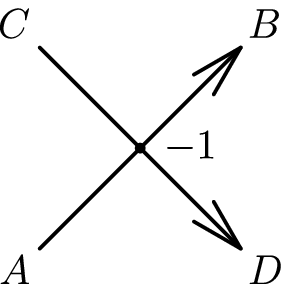}
\caption{Знак точки пересечения}
\label{f:sign}
\end{figure}

Пусть $(\overrightarrow{AB},\overrightarrow{CD})$ --- упорядоченная пара векторов (ориентированных отрезков) на плоскости, пересекающихся в точке $P$.
Определим {\bf знак} точки $P$ как $+1$, если $ABC$ ориентирован по часовой стрелке, и как $-1$ в противном случае
 (рис. \ref{f:sign}).

{\bf Коэффициент зацепления} $\lk$ плоской диаграммы ориентированного двухкомпонентного зацепления --- это сумма знаков всех точек пересечения на диаграмме, в которых первая компонента проходит над второй компонентой. В каждом перекрестке {\it первый} ({\it второй}) вектор --- это звено первой (второй) компоненты.

\begin{pr}\label{e:lk} Найдите коэффициент зацепления для (некоторой плоской диаграммы) зацепления Хопфа и пары колец Борромео, выбрав ориентацию компонент самостоятельно.
\end{pr}

\begin{lemma}\label{p:lk} Коэффициент зацепления сохраняется при движениях Рейдемейстера.
\end{lemma}

Доказательство аналогично доказательству леммы \ref{l:lk2}.
Достаточно проверить, что знаки всех перекрестков не изменяются.

По лемме \ref{p:lk} {\bf коэффициент зацепления} ориентированного двухкомпонентного зацепления
(или его класса изотопности) корректно определен как коэффициент зацепления любой плоской
диаграммы зацепления.

{\it Модуль коэффициента зацепления} (неориентированного) двухкомпонентного зацепления (или его класса изотопии) корректно определен вне зависимости от ориентации компонент.

Будем использовать без доказательства следующую {\it лемму о тривиальности:} для любых двух замкнутых ломаных в плоскости, вершины которых находятся в общем положении, сумма знаков их точек пересечения равна нулю.
Обсуждение и доказательство см. в \S1.3 <<Число пересечений для ломаных на плоскости>> в \cite{Sk18}, \cite{Sk}.

\begin{assertion}\label{p:lk-pr} (a) Перестановка компонент меняет знак коэффициента зацепления.

(b) Изменение ориентации любой компоненты меняет знак коэффициента зацепления.

(c) Существует ориентированное двухкомпонентное зацепление с коэффициентом зацепления $-5$.

(d) Для любой связной суммы $K\#L$ ориентированных двухкомпонентных зацеплений $K,L$ имеем $\lk(K\#L)=\lk K+\lk L$.

(e) Существует двухкомпонентное зацепление, не изотопное тривиальному зацеплению,
но имеющее нулевой коэффициент зацепления.
\end{assertion}

П. (e) доказывается с помощью \emph{многочлена Александера-Конвея}, см. \S\ref{0con}.

\begin{figure}[h]\centering
\includegraphics[scale=0.6]{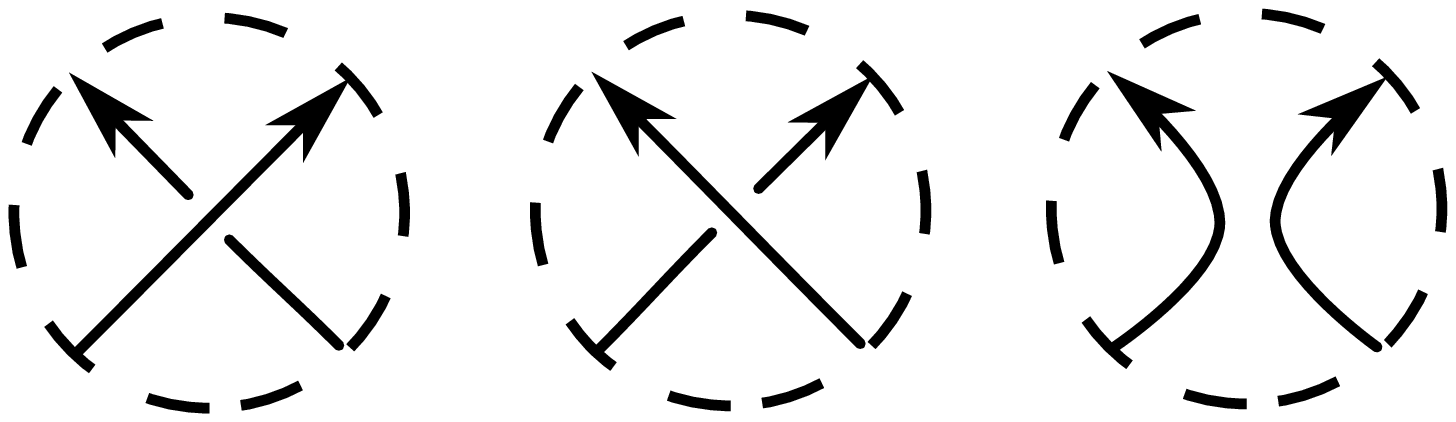}
\caption{Ориентированные зацепления $K_+,K_-,K_0$}
\label{f:skein}
\end{figure}

\begin{theorem}\label{con-lk}
Существует единственный целочисленный изотопический инвариант $\lk_2$ ориентированных двухкомпонентных зацеплений,  принимающий значение 0 на тривиальном зацеплении и такой, что для любых зацеплений $K_+$ и $K_-$, диаграммы которых отличаются изменением перекрестка $A$, как показано на рис.~\ref{f:skein},
$$\lk K_+-\lk K_-=\begin{cases}1 &A\text{ --- пересечение различных компонент;} \\
0 &A\text{ --- самопересечение одной из компонент.}\end{cases}$$
\end{theorem}

Доказательство аналогично доказательству теоремы \ref{con-lk2}.

\section{Инвариант Кэссона}\label{0cas}

{\bf Знак} перекрестка ориентированной плоской диаграммы узла определен после рисунка \ref{f:sign};
первый (второй) вектор --- это вектор перехода (прохода).
Очевидно, знак не зависит от ориентации диаграммы, и таким образом определен для неориентированных диаграмм.

{\bf Знак} $P$-скрещивающейся пары перекрестков на плоской диаграмме узла (для любой отмеченной точки $P$) --- это произведение знаков соответствующих перекрестков.

{\it $P$-инвариант Кэссона} плоской диаграммы --- это сумма знаков всех $P$-скрещивающихся пар перекрестков.

\begin{pr}\label{e:con-casca} (a) Аналог задачи \ref{e:con-cas}.b для инварианта Кэссона.

(b) Нарисуйте плоскую диаграмму узла и отмеченную точку $P$ такие, что $P$-инвариант Кэссона плоской диаграммы равен $-5$.
\end{pr}

\begin{lemma}\label{p:con-casca-ef} (a,b) Аналог леммы \ref{p:con-cas-ab}.a,b для инварианта Кэссона.
\end{lemma}

Поэтому {\bf инвариант} (число) {\bf Кэссона} $c_2$ плоской диаграммы узла (или даже изотопического класса узла) является корректно определенным инвариантом, если положить его равным $P$-инварианту Кэссона любой плоской диаграммы этого узла с любой отмеченной точкой $P$.

\begin{assertion}\label{con-casca} (a,b) Аналог утверждений \ref{wha-cons}.d и \ref{a:arf-pr} для инварианта Кэссона.
\end{assertion}

П. (b) доказывается с помощью \emph{многочлена Александера-Конвея}, см. \S\ref{0con}.

Обозначим через $D_+,D_-,D_0$ любые три диаграммы ориентированных (узлов или) зацеплений, отличающиеся, как показано на рис. \ref{f:skein}
(см. соглашения о рисунках в подписи к рис. \ref{f:crossing}).
Также обозначим через $K_+,K_-,K_0$ любые три зацепления с диаграммами $D_+,D_-,D_0$.

\begin{theorem}\label{con-thm}
Существует единственный целочисленный изотопический инвариант $c_2$ (неориентированных) узлов, который равен 0 для тривиального узла, и такой что
$$c_2(K_+)-c_2(K_-)=\lk K_0$$
для любых узлов $K_+$ и $K_-$, диаграммы которых отличаются как показано на рис.~\ref{f:skein}.
(Отметим что $K_0$ обязательно является двухкомпонентным зацеплением; число $\lk K_0$ определено корректно, так как изменение ориентации обеих компонент ориентированного зацепления не меняет коэффициента зацепления.)
\end{theorem}

\section{Полином Александера-Конвея}\label{0con}

Этот раздел использует только материал разделов \ref{0isot}, \ref{0isoli} и \ref{s:ori}
(кроме утверждений \ref{p:conzer}.b,c, использующих также разделы \ref{0rapl} и \ref{0cas}).

\begin{pr}\label{t:a3} (a) Существует единственный изотопический инвариант ориентированных трехкомпонентных зацеплений $\arf$ со значениями в $\mathbb{Z}_2$, принимающий значение 0 на тривиальном зацеплении и такой, что
$$\arf K_+ -\arf K_- =\begin{cases}\lk\phantom{}_2 K_0 &\text{в точке пересечения различных компонент;} \\
0 &\text{в точке самопересечения одной из компонент.}\end{cases}$$
(Здесь $\lk_2 K_0$ определено, так как $K_0$ обязательно является двухкомпонентным зацеплением.)\footnote{Это утверждение --- частный случай версии по модулю 2 теоремы~\ref{con-ale}; было бы интересно получить его прямое доказательство, поскольку оно могло бы проиллюстрировать идею доказательства теоремы~\ref{con-ale} на простейшем из нетривиальных случаев.
\newline
Теорема \ref{con-thma} является аналогом данного утверждения для узлов (которые являются однокомпонентными зацеплениями).
Определение инварианта $\arf$ из \S\ref{0arf} применимо только для узлов, и здесь мы пытаемся обобщить его на трехкомпонентные зацепления.}

(b) Предполагая существование инварианта $\arf$ из п.~(a), вычислите его для колец Борромео (при вашем выборе ориентаций компонент).
\end{pr}

\begin{theorem}\label{con-ale} Существует единственная бесконечная последовательность $c_{-1}=0,c_0,c_1,c_2,\ldots$ целочисленных инвариантов ориентированных неупорядоченных зацеплений такая, что

$\bullet$ на тривиальном узле $c_0=1$ и $c_1=c_2=\ldots=0$;

$\bullet$ для любого $n\ge0$ и зацеплений $K_+,K_-,K_0$ с рис.~\ref{f:skein} выполнено
$$c_n(K_+)-c_n(K_-)=c_{n-1}(K_0).$$
\end{theorem}

Доказательства существования в утверждении \ref{t:a3}.a и теореме \ref{con-ale} выходят за рамки данного текста.
%Мы будем пользоваться существованием без доказательства.
\footnote{Не вполне ясно, упорядоченные или неупорядоченные зацепления имеются в виду в \cite[\S2.3.1]{CDM}.
Поэтому мы приводимВывод сильной версии (для неупорядоченных зацеплений) из слабой (для упорядоченных зацеплений) см. в \S\ref{s:app}.}
Доказательство см. в \cite{Al28}, \cite[\S3-\S5]{Ka06'}, \cite{Ka06}, \cite[\S5.5 и \S5.6]{Ma18}, \cite{Ga20}.
Про связь с правильными раскрасками см. в \cite[\S6]{Ka06'}.

%Существование приведено в параграфе 5.6 ка следствие теоремы из параграфа 5.5.

Многочлен $C(K)(t):=c_0(K)+c_1(K)t+c_2(K)t^2+\ldots$ называется {\it многочленом Конвея}, см. утверждение \ref{p:conzer}.e.
Введение этого многочлена позволяет находить все инварианты $c_n$ так же быстро, как и один из них.
Формула из теоремы~\ref{con-ale} эквивалентна равенству
$$
C(K_+)-C(K_-)=tC(K_0).
$$

\begin{pr}\label{con-cal}
Найдите многочлены Конвея следующих зацеплений (выберите сами ориентацию для каждой компоненты).

(a) тривиальное зацепление с двумя компонентами; \quad

(b) тривиальное зацепление с $n$ компонентами; \quad

(c) зацепление Хопфа; \quad
(d) трилистник; \quad
(e) восьмерка; \quad

(f) зацепление Уайтхеда; \quad
(g) кольца Борромео; \quad
(h) узел $5_1$.
\end{pr}

\begin{assertion}\label{p:conzer} (a) Имеет место равенство $c_0(K)=1$ если $K$ узел, и $c_0(K)=0$ в противном случае (т.е. если $K$ имеет более одной компоненты).

(b) Если $K$ --- узел, то $c_{2j+1}(K)=0$, и $c_2$ --- инвариант Кэссона.

(c) Если $K$ --- двухкомпонентное зацепление, то $c_{2j}(K)=0$, и $c_1$ --- коэффициент зацепления.

(d) Если $K$ --- любое $k$-компонентное зацепление, то $c_j(K)=0$, если $j\le k-2$ или $j-k$ четно.

(e) Для каждого узла или зацепления все инварианты $c_n$, кроме конечного числа, обращаются в нуль.
\end{assertion}

\begin{assertion}\label{con-pro} (a) Изменение ориентации всех компонент зацепления (в частности, изменение ориентации узла) сохраняет многочлен Конвея.

(b) \cite[2.3.4]{CDM} Существует двухкомпонентное зацепление такое, что изменение ориентации одной его компоненты меняет степень многочлена Конвея (в частности, такое изменение ни сохраняет многочлен Конвея, ни меняет его знак).

(c) Для любой связной суммы $K\#L$ узлов $K,L$ верно $C(K\#L)=C(K)C(L)$.
\end{assertion}

\begin{assertion}\label{con-spl}
Зацепление называется \emph{разделяемым}, если оно изотопно зацеплению, компоненты которого содержатся в непересекающихся шарах.

(a) Зацепление Хопфа, зацепление Уайтхеда и кольца Борромео не являются разделяемыми зацеплениями.

(b) Коэффициент зацепления разделяемого зацепления равен нулю.

(c) Многочлен Конвея разделяемого зацепления нулевой.
\end{assertion}

\comment

{t:a3}
(b) Существует единственный изотопический инвариант ориентированных двухкомпонентных зацеплений $a_3$ со значениями в $\mathbb{Z}_2$, принимающий значение 0 на тривиальном зацеплении и такой, что
$$a_3(K_+)-a_3(K_-)=\arf K_0.$$
(Здесь $\arf K_0$ определено, так как $K_0$ либо узел, либо трехкомпонентное зацепленине.)

(c) Существует единственный изотопический инвариант ориентированных четырехкомпонентных зацеплений $a_3$ со значениями в $\mathbb{Z}_2$, принимающий значение 0 на тривиальном зацеплении и такой, что
$$a_3(K_+)-a_3(K_-)=\begin{cases}\arf K_0 &\text{в точке пересечения различных компонент;} \\
0 &\text{в точке самопересечения одной из компонент.}\end{cases}$$
(Здесь $\arf K_0$ определено, так как $K_0$ --- трехкомпонентное зацепление.)

\begin{pr}\label{con-cal2}
Вычислите (выбрав самостоятельно ориентации компонент)

(a) Инвариант $\arf$ колец Борромео;

(b,c,d) Инвариант $a_3$ зацепления Хопфа, зацепления Уайтхеда и {\it четырехкомпонентного зацепления Борромео},
т.е. такого четырехкомпонентного зацепления, для которого любое трехкомпонентное подзацепление изотопно тривиальному зацеплению, но само зацепление не изотопно тривиальному.
%вроде бы общепринятое название для этого `Brunnian link'=<<Брунново зацепление>>.
%добавить ссылку на рис. \ref{f:4-borromean}?
\end{pr}

\begin{figure}[h]\centering
\includegraphics[scale=0.6]{4-borromean.eps}
\caption{Четырехкомпонентное зацепление Борромео}
\label{f:4-borromean}
\end{figure}

\endcomment

\section{Инварианты Васильева-Гусарова}\label{0vas}

Этот раздел использует только материал разделов \ref{0isot}, \ref{0isoli} и \ref{s:ori}
(кроме задачи \ref{vas-pro}.2, использующей также раздел \ref{0cas}).

{\it Сингулярный} (ориентированный) {\it узел} --- замкнутая ориентированная ломаная в $\R^3$, все точки самопересечения которой являются простыми двойными (эти точки не являются вершинами ломаной).
Два сингулярных узла {\it изотопны}, если существует сохраняющий ориентацию кусочно-линейный гомеоморфизм $h\colon\R^3\to\R^3$, переводящий первый узел во второй с сохранением ориентации.
Обозначим через $\Sigma$ множество всех изотопических классов сингулярных узлов.

{\it Хордовая диаграмма} --- слово из $2n$ букв, среди которых ровно $n$ различных, причем каждая буква встречается дважды.
Слова, отличающиеся циклическим сдвигом, считаются одинаковыми.
Хордовую диаграмму обычно изображают в виде ориентированной окружности с набором хорд, ср. \cite[\S1.5]{Sk20}.
Для произвольного сингулярного узла $K$ обозначим через $\sigma(K)$ хордовую диаграмму, построенную следующим образом.
Будем равномерно двигаться по ориентированной окружности в направлении ориентации, и для каждой ее точки $A$ рассмотрим соответствующую точку $f(A)$ сингулярного узла $K$.
Соединим две точки окружности хордой, если они обе соответствуют некоторой точке самопересечения сингулярного узла $K$ \cite[4.8]{PS96}, \cite[3.4.1]{CDM}.\footnote{Другими словами, рассмотрим произвольное кусочно-линейное отображение $f:S^1\to\R^3$, образ которого есть $K$. Проведем хорду $XY$ для каждой пары точек $X,Y$ такой, что $f(X)=f(Y)$.
Хордовую диаграмму не следует путать с {\it диаграммой Гаусса} (проекции) узла $g:S^1\to\R^3$ (несингулярного), которая определяется как (произвольным образом ориентированная) хордовая диаграмма композиции проекции $\R^3\to\R^2$ с отображением $g$ \cite[4.8]{PS96} \cite[1.8.4]{CDM}.}

\begin{figure}[h]\centering
\includegraphics[scale=1]{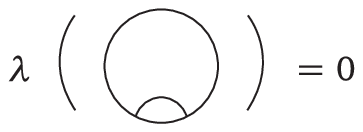}\qquad\qquad\includegraphics[scale=1]{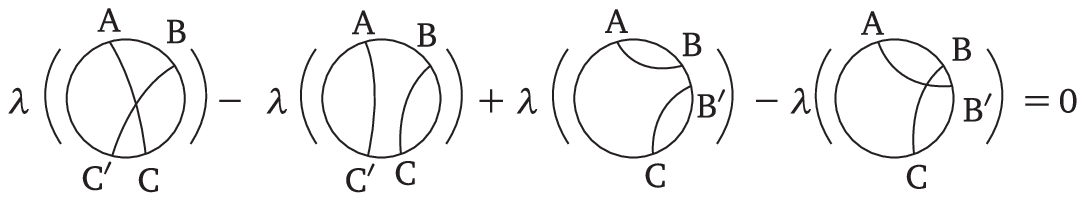}
\caption{Одночленное и четырехчленное соотношения}
\label{f:term}
\end{figure}

\begin{theorem}\label{vas-kon} Пусть $n\ge0$ --- целое и $\lambda\colon\delta_n\to\R$ --- функция из множества $\delta_n$ всех хордовых диаграмм с $n$ хордами.
Функция $\lambda$ удовлетворяет одночленному и четырехчленному соотношениям с рис. \ref{f:term} тогда и только тогда, когда существует отображение $v\colon\Sigma\to\R$ (т.е. инвариант сингулярных узлов), для которого

\begin{figure}[h]\centering
\includegraphics[scale=1]{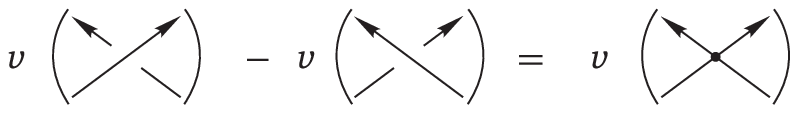}
\caption{Скейн-соотношение Васильева (обратите внимание на разницу с рис. \ref{f:skein})}
\label{f:vskein}
%<<скейн>>, а не <<скайн>> потому что так в Википедии
\end{figure}

(1) Выполнено скейн-соотношение Васильева с рис. \ref{f:vskein},

($2_n$) $v(K)=0$ для любого сингулярного узла $K$, имеющего более $n$ точек самопересечения,

(3) $v(K)=\lambda(\sigma(K))$ для любого сингулярного узла $K$, имеющего ровно $n$ точек самопересечения.
\end{theorem}

Доказательство выходит за рамки данного текста.

Насколько мне известно, теорема \ref{vas-kon} Васильева-Концевича никогда раньше не формулировалась в таком виде, коротком и удобном для вычисления инвариантов (хотя эта формулировка неявно использовалась в процессе вычислений).
%Поэтому
Я благодарен С.\,Чмутову за подтверждение того, что теорема \ref{vas-kon} верна и действительно эквивалентна теореме Васильева-Концевича в традиционной формулировке, см., например, \cite[теорема 4.2.1]{CDM}, ср. \cite[теорема 4.12]{PS96}.

Отображение $v\colon\Sigma\to\R$, для которой выполнено условие (1), называется {\it инвариантом Васильева-Гусарова}.
Если дополнительно выполнено условие ($2_n$), то $v$ называется {\it инвариантом порядка не выше $n$}.

Следующее утверждение означает, что функция $v$ из теоремы \ref{vas-kon} единственна с точностью до инварианта Васильева-Гусарова порядка не выше $n-1$.

\begin{assertion}[{\cite[предложение 3.4.2]{CDM}}]\label{vas-uni}
Разность $v-v'$ отображений $v,v'\colon\Sigma\to\R$, удовлетворяющих условиям (1), ($2_n$) и (3), удовлетворяет условиям (1) и ($2_{n-1}$).
\end{assertion}

\begin{pr}\label{vas-pro} (a) Докажите часть <<тогда>> теоремы \ref{vas-kon}.

(0),(1),(2) Докажите часть <<только тогда>> теоремы \ref{vas-kon} для $n=0,1,2$.

\emph{Указание.} В случае $n=2$ воспользуйтесь теоремой \ref{con-thm}.

(Часть <<только тогда>> теоремы \ref{vas-kon} для $n=3$ может быть доказана с использованием коэффициента при $h^3$ в $J(e^h)$, где $J$ --- $t$-параметризованный многочлен Джонса \cite[2.4.2, 2.4.3]{CDM} \cite[(4.6)]{PS96}.)
\end{pr}

В оставшихся задачах используйте без доказательства часть <<только тогда>> теоремы \ref{vas-kon}.
Далее утверждение <<$v(K)=x$ для любого особого узла $K$ с хордовой диаграммой $A$>> сокращенно записывается как <<$v(A)=x$>>.

\begin{pr}\label{vas-cal2}
(a) Существует единственный инвариант Васильева-Гусарова $v_2\colon\Sigma\to\R$ порядка не выше 2 такой, что $v_2(O)=0$ для тривиального узла $O$, и $v_2(1212)=1$.
(Здесь (1212) обозначает единственную <<нетривиальную>> хордовую диаграмму с двумя хордами, см. \cite[рис. 4.4]{PS96}, 3-я диаграмма в первом ряду.)

\emph{Указание.} Это следует из теоремы \ref{con-thm}, однако попробуйте вывести это из теоремы \ref{vas-kon}.

(b,b',c,d) Вычислите $v_2$ для (произвольным образом ориентированных) правого трилистника, левого трилистника, восьмерки и узла $5_1$.
\end{pr}

\begin{pr}\label{vas-cal3}
(a) Существует единственный инвариант Васильева-Гусарова $v_3\colon\Sigma\to\R$ порядка не выше 3 такой, что $v_3(O)=0$ для тривиального узла $O$ и левого трилистника $O$, и $v_3(123123)=1$.
(Здесь (123123) обозначает <<самую симметричную нетривиальную диаграмму с 3-я хордами>>, см. \cite[рис. 4.4]{PS96}, 5-ая диаграмма во втором ряду.)

(b,b',c,d) Аналог задачи \ref{vas-cal2}, для $v_3$.

\emph{Указание.} См. задачи 2, 3, 4ab и результаты/теоремы 11, 13, 14 из \cite[\S4]{PS96}.
\end{pr}

\begin{pr}\label{vas-cal4} (a) \cite[задача 4.4.b]{PS96} Существует единственный инвариант Васильева-Гусарова $v_4\colon\Sigma\to\R$ порядка не выше 4 такой, что

$\bullet$ $v_4(O)=0$ для тривиального узла $O$, левого трилистника $O$ и правого трилистника $O$,

$\bullet$ $v_4(12341234)=2$, $v_4(12341432)=3$, и $v_4(12341423)=5$.

(b,b',c,d) То же, что и в задаче \ref{vas-cal2}, для $v_4$.
\end{pr}

\section{Некоторые доказательства, указания  и решения}\label{s:app}

\begin{figure}[h]\centering
\includegraphics[scale=0.25]{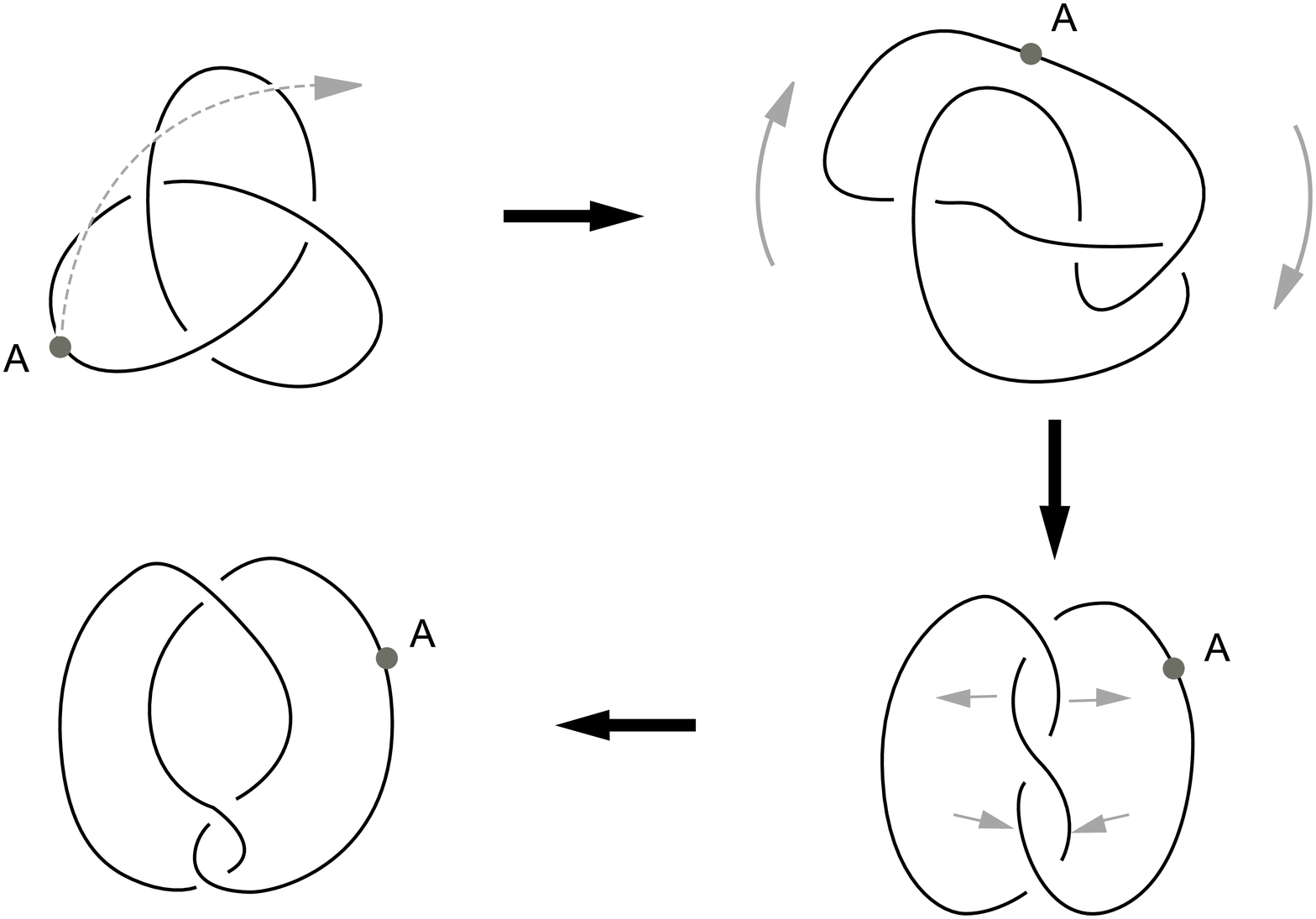}\qquad\includegraphics[scale=0.7]{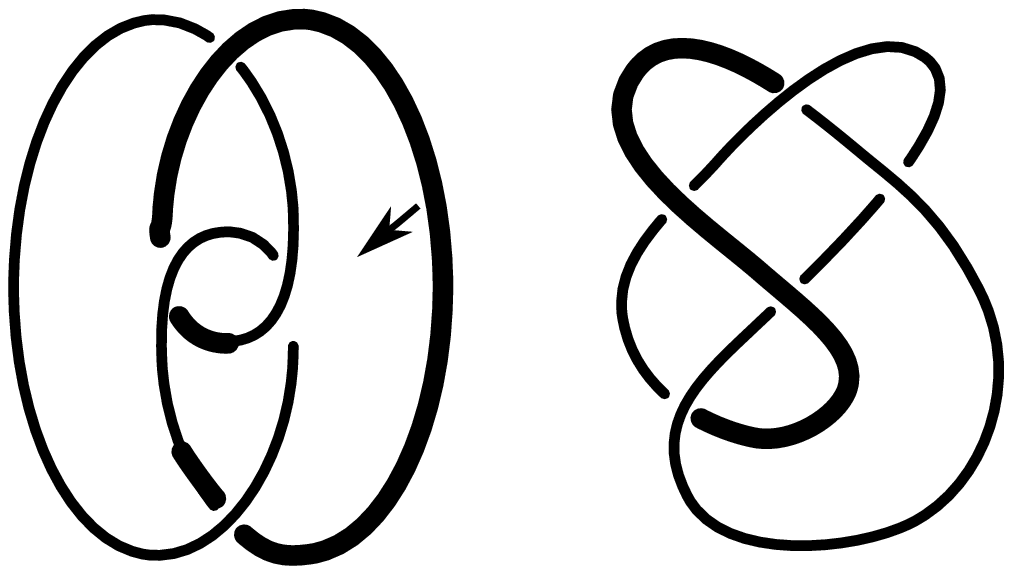}
\caption{Изотопии трилистника и восьмерки}
\label{f:itref}
\end{figure}

{\bf \ref{a:crossing}.} (a,b) <<Проще всего изготовить трилистник и восьмерку из веревки или шнурка, а затем попытаться получить из этих узлов все узлы, изображенные на рис. \ref{f:trefoil}. Выполнить некоторые преобразования трилистника и восьмерки вам поможет рис. \ref{f:itref}.>> \cite[\S2]{Pr95}
(Рис. \ref{f:itref} слева выполнен Д. Кроо.)

(e) Рассмотрим два узла с совпадающими плоскими диаграммами в <<горизонтальной>> плоскости $\pi$.
Для каждой точки $X$ в пространстве пусть $p(X)$ --- прямая, проходящая через $X$ и перпендикулярная $\pi$.
Пусть $h(X)$ --- высота точки $X$ относительно $\pi$, положительная ($h(X)>0$) если $X$ находится в верхнем полупространстве и отрицательная ($h(X)<0$) если $X$ находится в нижнем полупространстве.
Каждой точке $A$ первого узла поставим в соответствие точку $A’$ второго узла с помощью следующей процедуры.

\emph{Случай 1: Проекция точки $A$ на $\pi$ не является перекрестком на плоской диаграмме.}
В этом случае прямая $p(A)$ пересекает первый узел только в точке $A$.
Так как плоские диаграммы узлов совпадают, то прямая $p(A)$ пересекает второй узел также только в одной точке.
Обозначим эту точку через $A'$.

\emph{Случай 2: проекция точки $A$ на $\pi$ является перекрестком на плоской диаграмме.}
В этом случае прямая $p(A)$ пересекает первый узел еще в одной точке $B$.
Так как плоские диаграммы узлов совпадают, то прямая $p(A)$ пересекает второй узел тоже в двух точках: $C$ и $D$.
Без ограничения общности можно считать что $h(C)>h(D)$.
Тогда если $h(A)>h(B)$, то положим $A'=C$, иначе $A'=D$.

Для каждой точки $A$ первого узла и каждого числа $t\in[0,1]$ пусть $A(t)$ --- точка на прямой $p(A)$ на высоте $h(A(t))=(1-t)h(A)+t h(A')$.
По построению $A(0)=A$, $A(1)=A'$, и искомой изотопией является преобразование первого узла, двигающее точку $A(0)$ к точке $A(1)$ с постоянной скоростью так, что в момент времени $t$ она находится в положении $A(t)$, является искомой изотопией.

\begin{figure}[h]\centering
\includegraphics[scale=0.5]{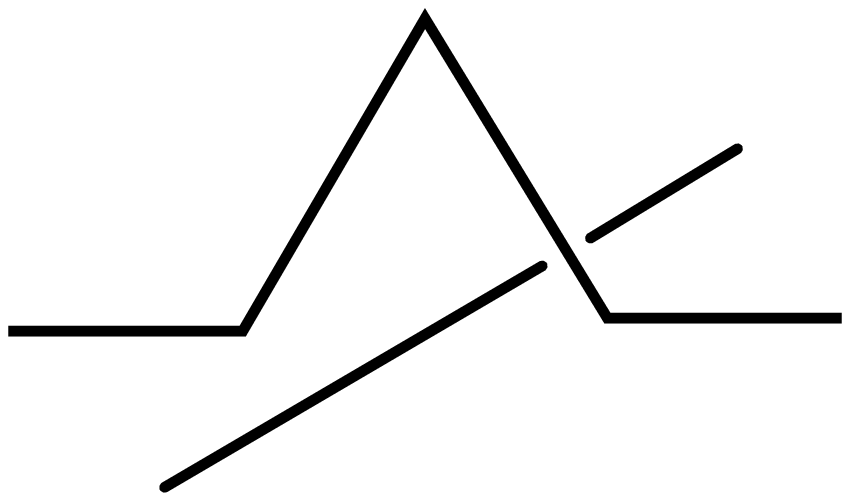}
\caption{Мостик через перекресток}
\label{f:bridge}
\end{figure}

\smallskip
{\bf \ref{p:diag}.} См. рис. \ref{f:bridge}.
Для каждого перекрестка на плоской диаграмме выберем на участке перекрестка, являющемся переходом, две точки, близкие к пересечению и находящиеся по разные стороны от него.
Заменим отрезок между двумя выбранными точками <<мостиком>>, соединяющим выбранные точки и <<подымающимся>> над плоской диаграммой.
После замены всех перекрестков такими <<мостиками>> мы получим искомый узел.

\smallskip
{\bf \ref{wha-tre}.} (a) Используйте утверждения \ref{e:con-cas}.b, \ref{e:con-casca}.a и леммы \ref{p:con-cas-ab}.a,b, \ref{p:con-casca-ef}.a,b.
Или можете воспользоваться утверждениями \ref{wha-cole}.a,b и леммой \ref{wha-colre}.

(b) Рассмотрим любую связную сумму $n$ трилистников.
Согласно утверждениям \ref{e:con-casca}.a и~\ref{con-casca}.a, инвариант Кэссона этого узла равен $n$.
Следовательно, по лемме \ref{p:con-casca-ef}.a,b соответствующие узлы не изотопны для различных значений $n$.

\smallskip
{\bf \ref{t:hopf}.} (a) Чтобы отличить зацепление Хопфа от остальных, воспользуйтесь утверждением \ref{a:lk2} и леммой \ref{l:lk2}.
Чтобы отличить зацепление Уайтхеда от тривиального зацепления, воспользуйтесь утверждением \ref{wha-cole} (или \ref{con-cal}) и леммой \ref{wha-colre} (или теоремой \ref{con-ale}).

(b) Используйте утверждение \ref{con-cal} и теорему \ref{con-ale}.

\smallskip
{\bf \ref{a:arfmot}.}
Выберем узел, проецирующийся на данную плоскую диаграмму тем же способом, что и в утверждении~\ref{p:diag}.
Будем считать, что все <<мостики>> лежат в верхнем полупространстве относительно плоскости проекции.
По предположению существуют точки $X$ и $Y$ узла, которые разбивают узел на две ломаные $p$ и $q$ такие, что

$\bullet$ $q$ лежит на плоскости проекции и проходит только через проходы;

$\bullet$ $p$ проецируется на ломаную $p'$, которая проходит только через переходы.

Возьмем точку $Z$ в верхнем полупространстве и точку $T$ в нижнем полупространстве.
Построим изотопию между заданным узлом и замкнутой ломаной $XZYT$, которая изотопна тривиальному узлу.
Построение изотопии состоит из трех шагов, в каждом из которых точки $X,Y$ остаются фиксированными.

{\it Шаг 1. Изотопия между $q$ и $XTY$.}
Предположим, что $q=A_0A_1\ldots A_n$, где $A_0=X$ и $A_n=Y$.
Тогда изотопия получается последовательностью элементарных движений
$$
A_0A_1\to A_0TA_1,\quad TA_1A_2\to TA_2,\quad TA_2A_3\to TA_3,\quad\ldots TA_{n-1}A_n\to TA_n.
$$

{\it Шаг 2. Изотопия между $p$ и $p'$.} Уберем все <<мостики>> с помощью элементарных движений.

{\it Шаг 3. Изотопия между $p'$ и $XZY$.} Строится аналогично шагу~1.

\smallskip
{\bf \ref{l:crossing}.} Следует из утверждения \ref{a:arfmot}.

{\it Другая идея доказательства (ср.~\cite[теорема 3.8]{PS96}).}
Обозначим через $\pi$ горизонтальную плоскость, содержащую диаграмму.
Для каждой точки $X$ в пространстве $p(X)$ и $h(X)$ определяются как в доказательстве утверждения~\ref{a:crossing}.c.
Пусть $l$ --- прямая на плоскости, которая проходит через вершину $A_0$ плоской диаграммы, и пусть вся диаграмма лежит по одну сторону от прямой $l$.
Пусть $A_0,A_1,\ldots, A_n$ --- все вершины плоской диаграммы в порядке их появления при движении вдоль диаграммы в некотором направлении.
Выберем точки $B_0,\ldots,B_n$ так, что $A_i\in p(B_i)$ для $i=1,\ldots,n$, и $h(B_i)<h(B_j)$ для $i<j$.
Пусть $B_{n+1}$ --- точка, проекция которой на $\pi$ близка к $A_0$, и $h(B_{n+1})>h(B_n)$.
Тогда узел $B_0\ldots B_nB_{n+1}$ изотопен тривиальному узлу.
Действительно, по построению прямой $l$, проекция узла на любую плоскость, перпендикулярную $l$, является замкнутой несамопересекающейся ломаной.
Остается только изменить перекрестки плоской диаграммы так, чтобы они согласовывались с проекцией построенного узла на плоскость $\pi$.

\smallskip
{\bf \ref{a:lk2}.} \textit{Ответ:} 1 для зацепления Хопфа и 0 для остальных зацеплений.

\smallskip
{\bf \ref{l:lk2}.} При движениях I и III число перекрестков, в которых первая компонента проходит над второй, не изменяется.
При движении II это число изменяется на 0 или $\pm 2$.
	
\smallskip
{\bf \ref{a:lk2-pr}.} (a) Рассмотрим плоскую диаграмму зацепления.
По лемме о четности (сформулированной перед утверждением \ref{a:lk2-pr}) число таких перекрестков, где первая компонента проходит над второй, имеет ту же четность что и число перекрестков, где вторая компонента проходит над первой.
Это в точности то, что нужно доказать.

(b) Примером является четвертое зацепление на рис. \ref{f:hopf}.
Используя целочисленный коэффициент зацепления, можно доказать, что это зацепление не изотопно тривиальному, см. \S\ref{0rapl}.

\smallskip
{\bf \ref{con-lk2}.} {\it Существование.} По лемме \ref{l:lk2} коэффициент зацепления по модулю 2 является изотопическим инвариантом.
Справедливость скейн-соотношения легко проверить непосредственно.

{\it Единственность.} Предположим, что $\mathrm f$ --- другой инвариант, отличный от $\lk_2$ и удовлетворяющий предположениям.
Тогда $\mathrm f-\lk_2$ является изотопическим инвариантом, равном нулю на тривиальном зацеплении и инвариантным относительно изменений перекрестков.
Согласно аналогу леммы \ref{l:crossing} для зацеплений любая плоская диаграмма зацепления может быть получена из диаграммы зацепления, изотопного тривиальному, при помощи изменений перекрестков.
Следовательно, $\mathrm f-\lk_2=0$.

\smallskip
{\bf \ref{e:con-cas}.} (a) Если $P$ --- такая точка на плоской диаграмме, как описано в утверждении~\ref{a:arfmot}, то не существует $P$-скрещивающихся пар перекрестков.
Следовательно, $P$-Арф-инвариант равен нулю.

\smallskip
{\bf \ref{p:con-cas-ab}.}
(a) Пусть $P_1$ и $P_2$ --- две отмеченные точки такие, что отрезок $P_1P_2$ содержит в точности один перекресток $X$.

{\it Случай 1: $P_1P_2$ участвует в проходе}.
Тогда $X$ не образует ни $P_1$-скрещивающейся, ни $P_2$-скрещивающейся пары ни с каким другим перекрестком.
Поэтому $P_1$- и $P_2$-Арф-инварианты диаграммы равны.

{\it Случай 2: $P_1P_2$ участвует в переходе}.
Тогда $X$ делит диаграмму на две замкнутые ломаные $q_1$ и $q_2$ такие, что $P_1$ лежит на $q_1$, а $P_2$ лежит на $q_2$.
Обозначим через $n_1$ (соответственно, $n_2$) число всех таких перекрестков на $q_1$ и $q_2$, в которых $q_1$ проходит над $q_2$ (соответственно, $q_2$ проходит над $q_1$).
Обозначим через $N_1$ количество $P_1$-скрещивающихся пар, образованных $X$ и некоторым пересечением ломаных $q_1$ и $q_2$.
Обозначим $P_1$-Арф-инвариант диаграммы $D$ через $\arf_{P_1}D$.
Аналогичные обозначения будем использовать и для отмеченной точки $P_2$.
Тогда
$$
{\arf}_{P_1} D-{\arf}_{P_2} D=N_1-N_2=n_1-n_2\underset{2}\equiv n_1+n_2\underset{2}\equiv0,
$$
где $D$ --- данная плоская диаграмма. В последней формуле

$\bullet$ первое равенство выполнено поскольку пара перекрестков в $D$ является $P_1$-скрещивающейся или $P_2$-скрещивающейся (но не одновременно!) тогда и только тогда когда она образована перекрестком $X$ и некоторым пересечением ломаных $q_1$ и $q_2$;

$\bullet$ второе равенство выполнено поскольку $N_1=n_1$ и $N_2=n_2$; действительно, точка пересечения ломаных $q_1$ и $q_2$ образует $P_1$-скрещивающуюся (соответственно, $P_2$-скрещивающуюся) пару с $X$ тогда и только тогда когда в этой точке $q_1$ проходит над (соответственно, под) $q_2$;

$\bullet$ $\underset{2}\equiv$ обозначает сравнение по модулю 2;

$\bullet$ последнее сравнение выполнено ввиду леммы о четности, примененной к ломаным $q_1$ и $q_2$.

\begin{figure}[h]\centering
\includegraphics[scale=0.45]{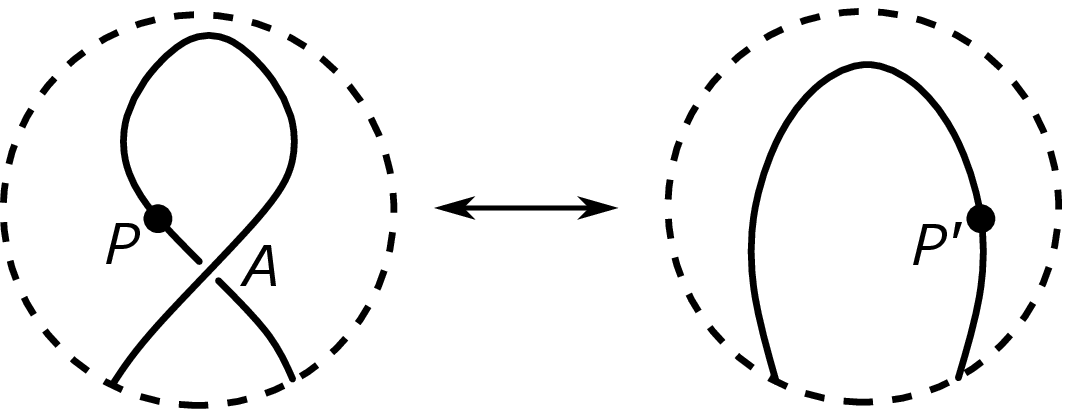}\qquad
\includegraphics[scale=0.45]{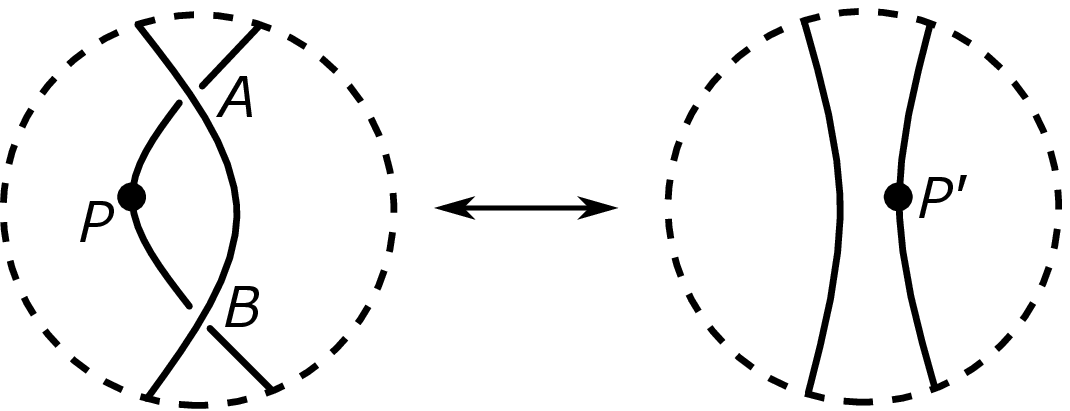}\qquad
\includegraphics[scale=0.45]{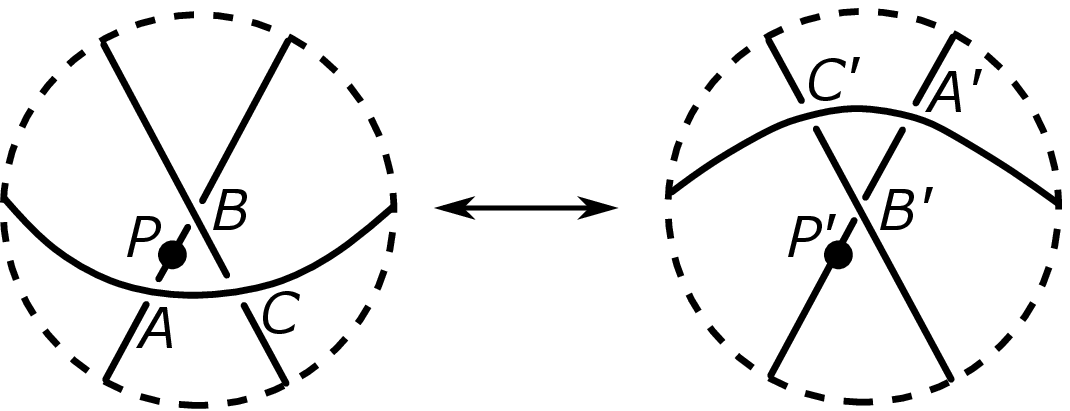}
\caption{Арф-инвариант не меняется при движениях Рейдемейстера}
\label{f:arf-reid}
\end{figure}

(b) {\it Движение типа $I$.} Выберем отмеченные точки до и после движения как показано на рис.~\ref{f:arf-reid} (слева).
Проверьте, что перекресток $A$ не образует $P$-скрещивающуюся пару ни с одним другим перекрестком.

{\it Движение типа $II$.} Выберем отмеченные точки до и после движения, как показано на рис.~\ref{f:arf-reid} (в центре).
Проверьте, что ни один из перекрестков $A$, $B$ не образует $P$-скрещивающуюся пару ни с каким другим перекрестком.

{\it Движение типа $III$.} Выберем отмеченные точки до и после движения, как показано на рис.~\ref{f:arf-reid} (справа).
Проверьте, что ни один из перекрестков $A$, $B$ не образует $P$-скрещивающуюся пару ни с каким другим перекрестком, и что ни один из перекрестков $A'$, $B'$ не образует $P'$-скрещивающуюся пару ни с каким другим перекрестком.
Затем проверьте, что произвольный перекресток $X$, отличный от $A$, $B$, $C$, образует $P$-скрещивающуюся пару с $C$, если и только если $X$ образует $P'$-скрещивающуюся пару с $C'$.

\smallskip
{\bf \ref{a:arf-pr}.} Рассмотрим любую связную сумму $K$ двух трилистников.
Из утверждения~\ref{wha-cons}.d следует, что $\arf K=0$.
В то же время по результату задачи~\ref{e:con-cas}.b и утверждению \ref{con-casca}.a $c_2(K)\ne0$.
Поэтому $K$ не изотопен тривиальному узлу.

\begin{figure}[h]\centering
\includegraphics[scale=0.6]{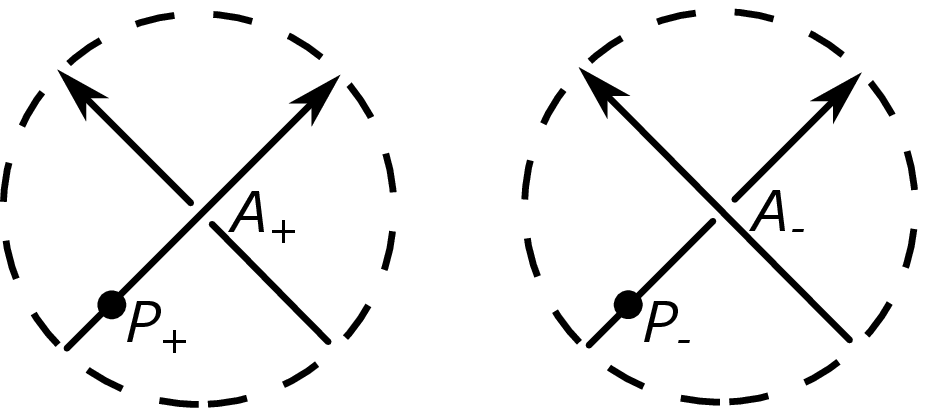}
\caption{К доказательству скейн-соотношения для Арф-инварианта}
\label{f:arf-skein}
\end{figure}

\smallskip
{\bf \ref{con-thma}.} {\it Существование.} По лемме \ref{p:con-cas-ab}, Арф-инвариант является изотопическим инвариантом.
Приведем указания к проверке скейн-соотношения.
Выберем отмеченные точки $P_+$, $P_-$, как показано на рис.~\ref{f:arf-skein}.
Проверьте, что перекресток $A_-$ не образует $P_-$-скрещивающуюся пару ни с каким другим перекрестком в $K_-$.
Затем проверьте, что число таких перекрестков, которые образуют $P_+$-скрещивающуюся пару с $A_+$ в $K_+$, сравнимо с $\lk\phantom{}_2K_0$ по модулю 2.

{\it Единственность.} Доказательство аналогично доказательству теоремы~\ref{con-lk2}.
Используйте саму лемму~\ref{l:crossing} вместо ее аналога для зацеплений.

\begin{figure}[h]\centering
\includegraphics[scale=0.7]{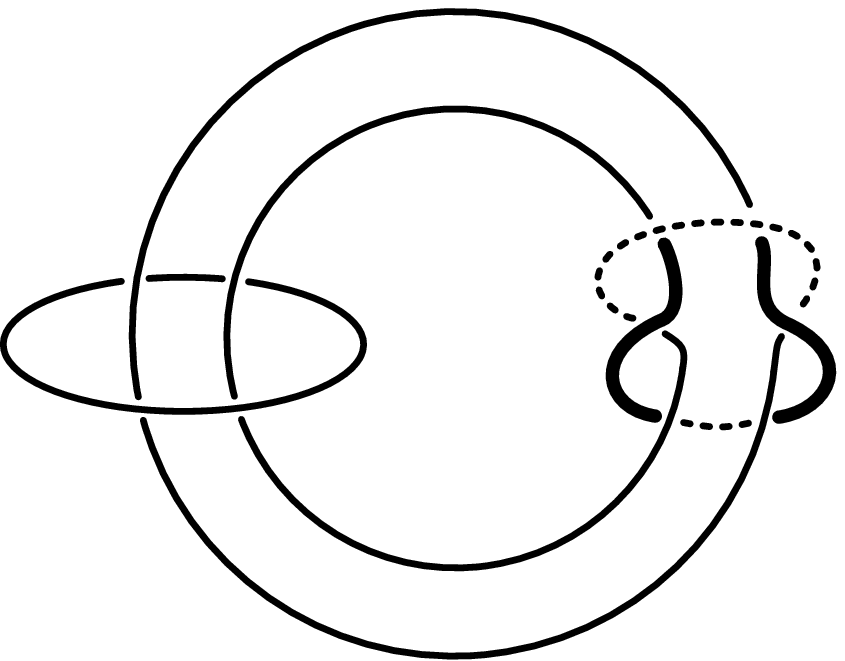}\qquad\includegraphics[scale=0.9]{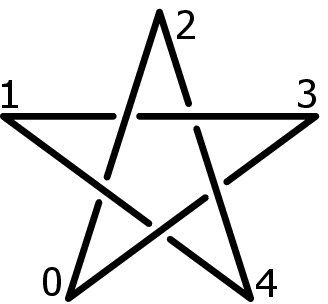}
\caption{Раскраска зацепления в 3 цвета и раскраска узла $5_1$ в 5 цветов}
\label{f:51-colored}
\end{figure}

\smallskip
{\bf \ref{wha-cole}.} {\it Ответы:} b,e,h --- 3-раскрашиваемы, a,c,d,f,g,i --- не раскрашиваемы.
Пример раскрашенного в три цвета трилистника см. в
\cite[стр. 30, рис. 4.3]{Pr95}.
Правильную раскраску последней диаграммы с рис. \ref{f:hopf} см. на рис. \ref{f:51-colored} слева.
(Эта диаграмма была ошибочно названа нераскрашиваемой в три цвета в \cite[\S4]{Pr95}. Эта небольшая ошибка была найдена Л. Баннором, С. Зотовой и Л. Кравцовой.)

\smallskip
{\bf \ref{wha-colr}.}
(a) Большая часть п. (a) следует из леммы \ref{wha-colre} и утверждений~\ref{wha-cole}.d-h
(см. \cite[стр. 30]{Pr95}).
Последняя диаграмма на рис. \ref{f:hopf} отличается от тривиального зацепления числом правильных раскрасок плоской диаграммы.
Докажите, что это число сохраняется при движениях Рейдемейстера.

(b) Плоская диаграмма называется {\it раскрашиваемой в 5 цветов}, если существует раскраска ее дуг в цвета $0,1,2,3,4$ такая, что хотя бы два цвета использованы, и в каждой точке пересечения если верхняя дуга имеет цвет $a$, а две нижние дуги --- цвета $b$ и $c$, то $2a\equiv b+c\pmod{5}$.
Аналогично лемме~\ref{wha-colre} можно показать, что раскрашиваемость в 5 цветов сохраняется при движениях Рейдемейстера.
Узел $5_1$ является раскрашиваемым в 5 цветов (см. рис.~\ref{f:51-colored} справа), а тривиальный узел --- нет.
Поэтому они не изотопны.

\smallskip
{\bf \ref{wha-chirsph}.} (b) {\it Первое решение.}
Назовем ориентированную ломаную {\it положительной}, если ограниченная ею часть плоскости лежит слева от любого ее звена (здесь используется теорема Жордана, сформулированная в замечании~\ref{r:accur}). Покажите, что положительность ломаной сохраняется при элементарных движениях.

{\it Указание ко второму решению.} Положительность ломаной из предыдущего решения можно эквивалентно определить следующим образом.
Будем говорить, что ориентированная ломаная $A_1\ldots A_n$ {\it положительна}, если для любой ее внутренней точки $O$ сумма ориентированных углов
$\angle A_1OA_2+\angle A_2OA_3+\ldots+\angle A_{n-1}OA_n+\angle A_nOA_1$
положительна (т.~е. {\it индекс} всякой внутренней точки ломаной относительно самой ломаной положителен).

\smallskip
{\bf \ref{wha-inv}.} Каждый из трех указанных ориентированных узлов переводится в ориентированный узел с противоположной ориентацией вращением на угол $\pi$ относительно <<вертикальной>> оси, проходящей через <<самую верхнюю>> точку узла
(см. самую левую диаграмму на рис. \ref{f:trefoil} в первом и втором ряду для трилистника и восьмерки соответственно). Это вращение входит в непрерывное семейство вращений на угол $\pi t,t\in[0,1]$ относительно той же прямой. Это и есть требуемая изотопия.

\smallskip
{\bf \ref{wha-consli}.} (d) Проверьте, что всякая точка пересечения разных компонент в $K\#L$ получается из точки пересечения разных компонент в $K$ или точки пересечения разных компонент в $L$.

\smallskip
{\bf \ref{e:lk}.} 	\textit{Ответы:} $\pm1$; 0.

\smallskip
{\bf \ref{p:lk-pr}.} (a) Доказательство аналогично доказательству утверждения \ref{a:lk2-pr}.a.
Рассмотрим плоскую диаграмму зацепления.
Из леммы о тривиальности (сформулированной перед утверждением \ref{p:lk-pr}) следует, что сумма знаков точек пересечения, в которых первая компонента выше второй, противоположна сумме знаков точек пересечения, в которых вторая компонента выше первой.
Поскольку перестановка компонент зацепления меняет знаки всех точек пересечения компонент, то коэффициент зацепления не меняется.

(b) Смена ориентации любой компоненты меняет знаки всех точек пересечения.

(c) Достаточно рассмотреть любую связную сумму пяти зацеплений Хопфа, у которых ориентации на компонентах выбраны таким образом, чтобы их коэффициент зацепления был равен~$-1$.

(d) Доказательство аналогично \ref{wha-cons}.d.
Множество перекрестков в плоской диаграмме зацепления $K\#L$ является объединением множеств перекрестков в плоских диаграммах зацеплений $K$ и $L$, причем знак каждого перекрестка в $K\#L$ совпадает с знаком этого перекрестка в соответствующем ему зацеплении $K$ или $L$.

(e) Примером является зацепление Уайтхеда (оно не изотопно тривиальному по теореме~\ref{t:hopf}.a).

\smallskip
{\bf \ref{e:con-casca}.} (a) \textit{Ответы:} 0, 1 и $-1$.

Тривиальный узел не имеет перекрестков, а значит в нем нет и скрещивающихся пар перекрестков. Поэтому инвариант Кэссона этого узла равен нулю.

Все три перекрестка трилистника имеют одинаковый знак.
Так как трилистник имеет только одну скрещивающуюся пару перекрестков (вне зависимости от выбора отмеченной точки), то инвариант Кэссона этого узла равен 1.

(b) Рассмотрим любую связную сумму пяти восьмерок.
По п.~(a) и утверждению~\ref{con-casca}.a ниже, инвариант Кэссона этого узла равен $-5$.

\smallskip
{\bf \ref{p:con-casca-ef}, \ref{con-casca}.a, \ref{con-thm}.} доказательства аналогичны доказательствам леммы~\ref{p:con-cas-ab}, утверждения~\ref{a:arf-pr} и теоремы~\ref{con-thma} соответственно.
Важно следить за знаками точек пересечения.
Для доказательства леммы \ref{p:con-casca-ef}.a воспользуйтесь леммой о тривиальности, сформулированной перед утверждением~\ref{p:lk-pr}.

\smallskip
{\bf \ref{con-casca}.} (b) Рассмотрим любую связную сумму трилистника и восьмерки.
По п.~(a) и утверждению \ref{e:con-casca}.a инвариант Кэссона этого узла равен нулю.
Однако по утверждениям~\ref{con-cal}.d,e и~\ref{con-pro}.c, многочлен Конвея этого узла равен $(1+t^2)(1-t^2)\ne1$.
Значит, этот узел не изотопен тривиальному узлу.

\smallskip
{\bf \ref{t:a3}.} (b) {\it Ответ:} 0.

{\it Замечание.} Инвариант $\arf=c_2\bmod2$
для зацеплений, вообще говоря, зависит от ориентации компонент (см. \cite[2.3.4]{CDM} для инварианта $c_3\bmod2$).

\smallskip
Пусть $D$ --- плоская диаграмма некоторого зацепления.
Обозначим через $\mathrm{cr}D$ количество перекрестков в $D$, а через $\mathrm{ch}D$ минимальное число изменений перекрестков, необходимое для получения из $D$ диаграммы зацепления, изотопного тривиальному (такая последовательность изменений перекрестков существует по аналогу леммы \ref{l:crossing} для зацеплений).

\smallskip
{\bf \ref{con-ale}.} Едиственность доказывается аналогично теоремам \ref{con-lk}, \ref{con-thm}; сначала решите задачу \ref{con-cal}.

{\it Вывод сильной версии (для неупорядоченных зацеплений) из слабой (для упорядоченных зацеплений).}
Достаточно доказать, что все инварианты $c_n$, определенные для упорядоченных зацеплений, сохраняются при изменении порядка компонент.

Пусть $D$ --- плоская диаграмма некоторого зацепления с двумя или более компонентами и пусть $D'$ --- плоская диаграмма, полученная из $D$ изменением порядка компонент.
Доказательство проведем индукцией по $\mathrm{cr}D$.
Если $\mathrm{cr}\,D=0$, то $D$ --- диаграмма зацепления, изотопного тривиальному, и по утверждению~\ref{con-cal}.b имеем $C(D)=0$ для любого порядка компонент.
Предположим, что $\mathrm{cr}\,D>0$; тогда продолжим доказательство индукцией по $\mathrm{ch}D$.
Если $\mathrm{ch}D=0$, то $D$ --- диаграмма зацепления, изотопного тривиальному; этот случай рассмотрен выше.
Предположим, что $\mathrm{ch}D>0$.
Пусть $D_*$ --- зацепление, полученное из $D$ таким изменением перекрестка, что $\mathrm{ch}D_*<\mathrm{ch}D$.
Обозначим через $D'_*$ зацепление, полученное из $D'$ изменением того же перекрестка. Тогда
$$
\pm(C(D)-C(D_*))=C(D_0)\quad\mbox{and}\quad\pm(C(D')-C(D'_*))=C(D'_0),
$$
где $D_0$ --- диаграмма зацепления $K_0$ (с некоторым порядком компонент) на рис.~\ref{f:skein}, $D$ и $D_*$ --- это $D_+$ и $D_-$ в некотором порядке, а $D'_0$ обозначает то же самое для $D'$, $D'_*$.
%эта фраза, по-моему, непонятна; но придумать, как сделать лучше, трудно
Заметим, что диаграммы $D_*$ и $D'_*$ совпадают с точностью до порядка компонент.
То же верно для диаграмм $D_0$ и $D'_0$.
Так как $\mathrm{ch}D_*<\mathrm{ch}D$ и $\mathrm{cr}\,D_0<\mathrm{cr}\,D$, то по предположению индукции имеем $C(D_*)=C(D'_*)$ и $C(D_0)=C(D'_0)$.
Тогда $C(D)=C(D')$.

\smallskip
{\bf \ref{con-cal}.} {\it Ответы:} (a, b) 0; (c) $\pm t$; (d) $1+t^2$; (e) $1-t^2$; (f) $\pm t^3$; (g) $\pm t^4$; (h) $1+3t^2+t^4$.

{\it Замечание.} Знаки в ответах к пп.~(c), (f), (g) зависят от ориентации компонент.

{\it Указание.} Примеры вычислений для пп.~(a), (c), и (d) см. в~\cite[2.3.2]{CDM}.

\smallskip
{\bf \ref{p:conzer}.} Пусть $D$ --- плоская диаграмма данного зацепления $K$.

(a) Для любой диаграмы $D_*$, полученной из $D$ изменением перекрестка, имеем $c_0(D)-c_0(D_*)=0$,
т.~е. $c_0$ не меняется при изменениях перекрестков.
По аналогу леммы~\ref{l:crossing} для зацеплений диаграмма $D$ может быть получена изменениями перекрестков из диаграммы зацепления, изотопного тривиальному.
Утверждение следует из определения $c_0$ на тривиальном узле и утверждения~\ref{con-cal}.b.

(b,c) Первые части обоих утверждений являются частными случаями п.~(d).
Вторые части обоих утверждений следуют из определения $c_1,c_2$ и теорем~\ref{con-lk}, \ref{con-thm}.

(d) Доказательство индукцией по $\mathrm{cr}\,D$.
Если $\mathrm{cr}\,D=0$, то $K$ изотопно тривиальному зацеплению.
Если $K$ --- узел, то $C(D)=1$.
В противном случае, $C(D)=0$ по утверждению~\ref{con-cal}.b.
Предположим, что $\mathrm{cr}\,D>0$; тогда продолжим доказательство индукцией по $\mathrm{ch}D$.
Если $\mathrm{ch}D=0$, то $K$ изотопно тривиальному зацеплению; этот случай рассмотрен выше.
Предположим, что $\mathrm{ch}D>0$.
Пусть $D_*$ --- зацепление, полученное из $D$ изменениями перекрестков и такое, что $\mathrm{ch}D_*<\mathrm{ch}D$.
Тогда $\pm(c_j(D)-c_j(D_*))=c_{j-1}(D_0)$, где $D_0$ --- диаграмма с рис.~\ref{f:skein}, $D$ и $D_*$ --- это $D_+$, $D_-$ в некотором порядке.
Заметим, что зацепление $D_*$ состоит из $k$ компонент, а зацепление $D_0$ состоит из $k'=k\pm1$ компонент. Поэтому если $j\le k-2$, то $j-1\le k'$, и если $j-k$ четно, то и $(j-1)-k'$ четно.
Так как $\mathrm{ch}D_*<\mathrm{ch}D$ и $\mathrm{cr}\,D_0<\mathrm{cr}\,D$, то по предположению индукции имеем $c_j(D_*)=c_{j-1}(D_0)=0$.
Тогда $c_j(D)=0$.

(e) Докажите аналогично п.~(d) что {\it $c_j(D)=0$ для любой плоской диаграммы $D$ и $j>\mathrm{cr}\,D$}.

\smallskip
{\bf \ref{con-pro}.} (a) Доказательство аналогично доказательству утверждения~\ref{p:conzer}.d.

(c) Пусть $D$ и $E$ --- плоские диаграммы зацеплений $K$ и $L$.
Аналогично утверждению~\ref{p:conzer}.d покажите. что $C(D\#E)=C(D)C(E)$ индукцией по $\mathrm{cr}\,D$ при фиксированном $E$.

\begin{figure}[h]\centering
\includegraphics[scale=0.6]{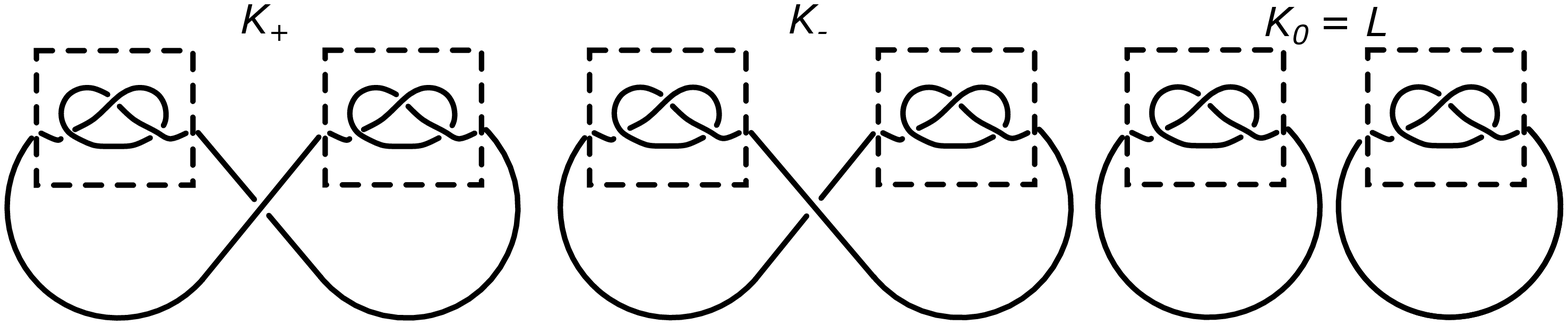}
\caption{Доказательство равенства $C(\mbox{разделяемое зацепление})=0$}
\label{f:con-split}
\end{figure}

\smallskip
{\bf \ref{con-spl}.} (a) Следует из утверждений~\ref{con-cal}.c,f,g и пп.~(b,c).

(c) Если $L$ --- разделяемое зацепление, то существуют зацепления $K_+$, $K_-$, $K_0$ такие, что

$\bullet$ их плоские диаграммы отличаются как на рис.~\ref{f:skein};

$\bullet$ зацепления $K_+$ и $K_-$ изотопны;

$\bullet$ зацепление $K_0$ изотопно $L$.

Получаем $C(L)=C(K_0)=\frac{1}{t}(C(K_+)-C(K_-))=0$, см. рис. \ref{f:con-split}.

%\newpage

{\it В этом списке звездочками отмечены книги, обзоры и популярные статьи.}